\documentclass[10pt,hidelinks]{article}
\usepackage{authblk}

\author[1]{Steen Ryom-Hansen\thanks{supported in part by FONDECYT grant 1221112
and EPSRC grant EP/W007509/1}}

\affil[1]{Instituto de Matem\'aticas, Universidad de Talca, Chile}

\usepackage{cancel}
\usepackage{pstricks-add}
\usepackage{yfonts}
\usepackage{pstricks-add}
\usepackage{enumitem}
\usepackage{calc}

\usepackage{pst-func}
\usepackage{stackrel}
\usepackage{bbold}
\usepackage[normalem]{ulem}
\usepackage{etex}
\usepackage[left=2cm,top=2cm,right=2cm,bottom=2cm]{geometry}
\usepackage{amsmath, amsthm, amssymb}
\allowdisplaybreaks
\usepackage{pst-all}
\usepackage[textwidth=0.5in]{todonotes}
\usepackage[rightcaption]{sidecap}

\usepackage{relsize}
\usepackage{float}
\usepackage[caption = false]{subfig}

\usepackage[hypertexnames=false]{hyperref}

\usepackage{multicol,lipsum}
\usepackage[scaled=1]{helvet}
\usepackage{titlesec}
\usepackage{multido}
\usepackage{graphicx}
\usepackage{multirow}
\usepackage{blindtext}
\usepackage{color}
\usepackage{lipsum}
\usepackage{mathtools}
\usepackage[vcentermath]{youngtab}
\usepackage{young}
\usepackage[all,ps,dvips,graph]{xy}
\usepackage[misc,geometry]{ifsym}

\usepackage{titlesec}
\usepackage{lineno,hyperref}

\let\OLDthebibliography\thebibliography
\renewcommand\thebibliography[1]{
  \OLDthebibliography{#1}
  \setlength{\parskip}{0pt}
  \setlength{\itemsep}{0pt plus 0.3ex}
}

\titleformat{\section} {\normalfont\scshape \large \centering}{ \thesection}{1em}{}
\definecolor{morado}{rgb}{0.5,0,0.5}

\newcommand{\F}{  \Phi} 

\newcommand{\C}{  \hat{c}} 
\newcommand{\logq}{  {\log}} 
\newcommand{\JM}{{\bf JM} }

\newcommand{\QQ}{\mathbb Q}

\newcommand{\ord}{{\rm{ord}}}

\newcommand{\pos}{{\rm{pos}}}

\newcommand{\ii}{\boldsymbol{i}}
\newcommand{\jj}{\boldsymbol{j}}

\newcommand{\II}{\mathcal{I}_n}
\newcommand{\JJ}{\mathcal{J}}
\newcommand{\JJn}{\mathcal{J}_n}
\newcommand{\KK}{\mathcal{K}}
\newcommand{\KKn}{\mathcal{K}_n}

\newcommand{\KKii}{\mathcal{K}_{ \ii}}
\newcommand{\KKjj}{\mathcal{K}_{ \jj}}

\newcommand{\Kdomain}{\Bbbk}

\newcommand{\TT}{{\mathcal T }_{q,t}}

\newcommand{\SSS}{{\mathcal S }_{q,t}}

\newcommand{\nrelbar}{\,\not\!\!\!\!\;\text{---}\,}

\newcommand{\FF}{ { \mathbb F}}

\newcommand{\OO}{\mathcal{O}}
\newcommand{\kk}{\boldsymbol{k}}

\newcommand{\A}{\mathcal{A}}
\newcommand{\Z}{\mathbb{Z}}
\newcommand{\e}{ \overline{e}}

\newcommand{\res}{ \textrm{res} }

\newcommand{\Par}{{\rm Par}   }
\newcommand{\Comp}{{\rm Comp}   }

\newcommand{\spa}{{\rm span}}

\newcommand\bs{\mathbf{s}}
\newcommand\bt{\mathbf{t}}
\newcommand\bu{\mathbf{u}}

\newcommand{\s}{\mathfrak{s}}

\newcommand{\q}{\hat{q}}

\newcommand{\T}{  \mathfrak{t}}

\newcommand{\Si}{\mathfrak{S}}

\newcommand{\std}{{\rm Std}}
\newcommand{\Rstd}{{\rm RStd}}
\newcommand{\Cstd}{{\rm CStd}}
\newcommand{\tab}{{\rm Tab}}

\newcommand{\Hecken}{ {\mathcal{H}_n(q)}}

\newcommand{\YYresiduefield}{ \mathcal{Y}^{\, \Kdomain}_{d, n}(q)}
\newcommand{\YYZ}{ \mathcal{Y}^{\, \Z}_{d, n}(q)}
\newcommand{\YYfieldext}{ \mathcal{Y}^{\, \Kdomain^{\prime}}_{d, n}(q)}
\newcommand{\YYOO}{ \mathcal{Y}^{\, \OO}_{d, n}(\q)}
\newcommand{\YYfractionfield}{ \mathcal{Y}^{\, \Kdomain(\q)}_{d, n}(\q)}
\newcommand{\YYfractionfielda}{ \mathcal{Y}^{\, a}_{d, n}(\q)}

\newcommand{\Eord}{ {\mathcal E}^{\ord}_n(q)}
\newcommand{\Eordfield}{ {\mathcal E}^{\ord, \Kdomain}_n(q)}
\newcommand{\Efield}{ {\mathcal E}^{\, \Kdomain}_n(q)}
\newcommand{\Efieldext}{ {\mathcal E}^{\, \Kdomain^{\prime}}_n(q)}

\newcommand{\Efractionfield}{ {\mathcal E}^{\, \Kdomain(\q)}_n(\q)}

\newcommand{\E}{ {\mathcal E}_n(q)}
\newcommand{\Eprime}{ {\mathcal E}^{\prime}_n(q)}

\newcommand\blambda{{\boldsymbol\lambda}}

\newcommand\bnu{{\boldsymbol\nu}}

\newcommand\bmu{{\boldsymbol\mu}}

\newcommand{\YY}{ \mathcal{Y}_{d, n}(q)}
\newcommand{\RR}{ {\mathcal R}_n(\Gamma_{e,d})}
\newcommand{\REs}{ {\mathcal E}^{sub}_n(\Gamma_{e})}
\newcommand{\REords}{ {\mathcal E}^{\ord, sub}_n(\Gamma_{e,d})}
\newcommand{\RE}{ {\mathcal E}_n(\Gamma_{e})}
\newcommand{\REord}{ {\mathcal E}^{\ord}_n(\Gamma_{e})}

\newcommand{\YYdeqone}{ \mathcal{Y}_{1, n}(q)}
\newcommand{\RRdeqone}{ {\mathcal R}_n(\Gamma_{e})}

\newcommand{\bT}{\pmb{\mathfrak{t}}}

\newcommand{\SetPar}{ \operatorname{SetPar}}

\newcommand{\OrdSetPar}{ \operatorname{OrdSetPar}}

\modulolinenumbers[5]

\newtheorem{theorem}{Theorem}
\newtheorem{lemma}[theorem]{Lemma}

\newtheorem{definition}[theorem]{Definition}
\newtheorem{corollary}[theorem]{Corollary}

\newtheorem{remark}[theorem]{Remark}

\newenvironment{dem}{\noindent \textit{Proof:} }{\quad \hfill $\square$}

\numberwithin{equation}{section}

\newgray{plomoclaro}{.90}
\newgray{plomooscuro}{.70}
\newgray{blanco}{1}

\begin{document}
\Yvcentermath1
\sidecaptionvpos{figure}{lc}


\title{A KLR-like presentation for the bt-algebra }

\date{\vspace{-5ex}}
\maketitle
\begin{abstract}
  We consider the bt-algebra ${  \mathcal E}_n(q)$ of knot theory, defined over an
  arbitrary field $ \Bbbk$. 
    We find a KLR-like presentation for $ {\mathcal E}_n(q) $ showing
      that it is a $ \mathbb Z$-graded algebra if 
  $ q \in \Bbbk^{\times} \setminus \{1 \} $ admits
  a square root in $ \Bbbk $. We introduce the ordered bt-algebra
 $ \Eord$ and show that it also has a KLR-like presentation, 
  without restriction on $ q \in \Bbbk^{\times} \setminus \{1 \} $. In
  particular, $ {\mathcal E}^{\ord}_n(q)$ is a $ \mathbb Z$-graded algebra for all
  $ q \in \Bbbk^{\times} \setminus \{1 \} $. 
\end{abstract}

\section{Introduction}
In this article we study the representation theory of the Yokonuma-Hecke algebra 
$ \YY $ of type $A$ and the related braids and ties algebra $ \E$, also denoted bt-algebra for short.
The Yokonuma-Hecke algebra was
introduced in the sixties by Yokonuma, see \cite{Yokonuma}, as a generalization of the Iwahori-Hecke algebra,
whereas the bt-algebra was introduced in the noughts by Aicardi and Juyumaya, see
\cite{AicardiJuyumaya1}, via a deframization process involving $ \YY$. 

\medskip
Much of the current interest in $ \YY $ and $ \E $ is derived from their applications in
knot theory, see for example
\cite{AicardiJuyumaya}, \cite{ArcisJuyu}, \cite{ChJuKaLa} and \cite{MarceloFlores}, 
where Jones'
Markov trace formalism for the HOMFLYPT
polynomial is suitably adapted to the $ \YY $ and $ \E $ settings. In turns out that
the associated link invariants 
are different from the HOMFLYPT polynomial, and in fact stronger than
the HOMFLYPT polynomial
on links that are not knots.

\medskip
The interest in both algebras, however, goes much further than their applications in knot
theory,
and in fact both algebras have been studied from a wide range of perspectives in recent years.
In particular, it has been shown that $ \YY $ is isomorphic to Shoji's modified Ariki-Koike algebra, 
whereas 
$ \E $ has connections to Jones and Martin's partition algebra,
Clifford theory, Coxeter group theory, geometric representation theory and much more, see
for example \cite{Banjo}, \cite{ERH}, \cite{ KimKim}, \cite{LaPo}, \cite{Lu1}, \cite{Lu2},
\cite{Marin}, \cite{Martin}, \cite{Shoji}.

\medskip
For generic choices of the parameter $ q $, the algebras $ \YY $ and $ \E $ are semisimple,
but our main focus is here 
the root of unity case where they are non-semisimple, in general.

\medskip
Since Soergel's seminal papers in the nineties,  see \cite{Soergel}, \cite{Soergel2}, 
$\Z$-gradings on non-semisimple algebras 
have been a  
recurring theme in representation theory of objects of Lie type, 
see for example  \cite{AJS}, \cite{EhStro}, \cite{HaziMartinParker}, \cite{LehrerLyu},
\cite{LoPlaRy}, \cite{PlaRy}, \cite{Ro2}, \cite{Catharina}.
The methods for obtaining these gradings are very disperse, but 
for the present paper, 
Brundan, Kleshchev and Rouquier's isomorphism Theorem $ \Hecken \cong 
R_n(\Gamma_e)  $, with $ \Hecken $ denoting the Iwahori-Hecke algebra for 
$ q $ an $ e$'th primitive root of unity
and $ R_n(\Gamma_e)  $ the KLR-algebra 
for $ \Gamma_e $ a cyclic quiver of order $ e$, is especially relevant.
Indeed, $ R_n(\Gamma_e)$ is $\Z$-graded, 
and so $ \Hecken $ becomes $\Z$-graded as well via this isomorphism.

\medskip
The starting point of our work is 
Rostam's paper \cite{Ro}, in which Brundan, Kleshchev and Rouquier's isomorphism 
is extended to an isomorphism $ \YY \cong  \RR $
where $ q $ still is an $e$'th primitive root of unity but $ R_n(\Gamma_{e, d} )$ is the KLR-algebra 
for the disjoint union $ \Gamma_{e,d} $ of
$ d $ cyclic quivers each of order $ e $.
Using his results,   
we obtain in part $   { \bf a )} $ of our main Theorem 
\ref{twomaintheoremsX} 
a non-trivial $ \Z$-grading on $ \E $ through an isomorphism 
$ \E \cong \RE $, where $ \RE $ is an algebra given by a KLR-like presentation. 
On the other hand, it should be pointed out that, strictly speaking, $ \RE $  is not a KLR-algebra,
even in the general sense of Rouquier,
see \cite{Rouq}.

\medskip
As suggested by the name, the bt-algebra $ \E $ is a diagram algebra with diagram basis consisting of \lq braids\rq\ and \lq ties\rq\, as follows 
\begin{equation}\label{btdiagramIntro}
D=  \raisebox{-.45\height}{\includegraphics[scale=0.7]{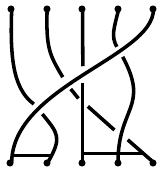}} 
\end{equation}
with the ties of $ D $ defining a set partition $ A_D $ on $ \{1,2, \ldots, n \}$, for example 
$ A_D = \{ \{ 1,2 \},\{ 3,4,5 \} \} $ in \eqref{btdiagramIntro}. 
In particular, the dimension of $ \E $ is $ n! b_n $, where $ b_n$ is the $n$'th Bell number, that
is the cardinality of $ \SetPar_n $. 

\medskip
This dependency on $\SetPar_n$ is reflected in the 
definition of $ \RE$ which is given by generators
\begin{equation} 
 \{ \e(\ii,A)  | \ii \in (\Z/ e \Z)^n,  A \in \SetPar_n \}
 \cup \{ y_1, y_2, \ldots, y_n \} \cup \{ \psi_1, \psi_2, \ldots, \psi_{n-1} \}, 
\end{equation}
subject to a long list of relations, see
Definition \ref{quiverE}, in the spirit of 
the usual KLR-relations for $ R_n(\Gamma_{e} )$. 
Similarly to $ R_n(\Gamma_{e} )$ there is also a diagram calculus 
associated with $\RE$, for example for $ e = 3, n= 4, \ii = (1,2,0,1), A = \{ \{1,2\}, \{3,4 \} \} $
we have 
\begin{equation}\label{720}
  \e(\ii,A) \mapsto  \raisebox{-.5\height}{\includegraphics[scale=0.7]{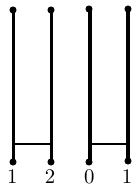}} \raisebox{-20\height}{,}
  \, \,  \, \,   
y_3 e(\ii, A) \mapsto   \raisebox{-.5\height}{\includegraphics[scale=0.7]{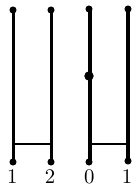}} \raisebox{-20\height}{,}  \, \,  \, \, 
\psi_2 e(\ii,A) \mapsto   \raisebox{-.5\height}{\includegraphics[scale=0.7]{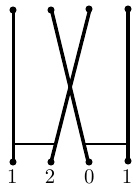}} \, \, \, 
 \raisebox{-20\height}{,}
\end{equation}
and via this calculus the relations for $ \RE $ may be rewritten in a diagrammatic form, that we give in 
\eqref{REdiag1A}-\eqref{REdiag1Alast}.

\medskip
In \cite{ERH} it was shown that the bt-algebra $ \E $ is a cellular algebra
in the sense of \cite{GL}. It is 
endowed with a family $ {\cal L}_n = \{ X_1, X_2, \ldots, X_n \} $ of
Jucys-Murphy elements 
whose generalized eigenspaces and intertwiners between them
play a prominent role throughout the paper. On the other hand, 
$ {\cal L}_n $ is not a separating family of $\JM$-elements for $ \E $ in the sense of
\cite{Mat-So}, and this fact causes several technical problems. 
For the purpose of this paper, however, 
we can solve these problems by embedding $ \E = \Efield $ in an appropriate $ \YYfieldext $,
on which we can apply the full $ \JM$-theory, see section \ref{JM for YY}. 

\medskip
Our isomorphism $ \E \cong \RE $
in  $   { \bf a )} $ of Theorem 
\ref{twomaintheoremsX}
holds for any ground field $ \Kdomain $. For brevity of the exposition we do not consider the
degenerate case $q =1 $, but 
unfortunately we also need
to impose the technical condition that $ q \in \Kdomain  $
admit a square root $ q^{\frac{1}{2}} $ in $ \Kdomain $. 
(For example, if $ \Kdomain = \FF_q $ this condition can be checked using Legendre symbols and 
quadratic reciprocity). 
The condition on the existence of $ q^{\frac{1}{2}} $ is related to the lack
of a natural total order on the blocks of a general set partition $ A \in \SetPar_n $.
This leads us to introduce {\it the ordered bt-algebra $ \Eord$},
see Definition
\ref{orderedBT}, 
a diagram algebra similar to $ \E $, but with a diagram basis consisting of braids and {\it ordered ties} as follows
\begin{equation}\label{btdiagramIntroA}
D=  \raisebox{-.45\height}{\includegraphics[scale=0.7]{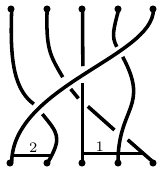}} 
\end{equation}
with the ties of $ D $ defining an {\it ordered set partition} $A_D^{\ord} $, for
example $ A_D^{\ord} = \{ \{3, 4,5 \}, \{1,2\} \} $ in
\eqref{btdiagramIntroA}. 
In particular, $ \dim \Eord = n! f_n $ where $ f_n$ is
the $n$'th Fubini number, the cardinality of ordered set partitions on $ n$. 

\medskip
Similarly, in Definition \ref{quiverEord} we introduce the algebra $ \REord $, an ordered 
set partition version of
$ \RE$ also given by a KLR-like presentation.  
With these definitions, in part $   { \bf b)} $ of Theorem
\ref{twomaintheoremsX}, we establish our second main algebra isomorphism $ \Eord \cong  \REord$,
that this time holds without restriction on $ q \in \Kdomain^{\times} \setminus \{1\}$. 
Since $ \REord $ is a $ \Z$-graded algebra, we conclude that $ \Eord $ is a $ \Z$-graded algebra
for any choice of $ q \in \Kdomain^{\times} \setminus \{1\}$. 

\medskip
Let us give a brief overview of the contents of the paper.

\medskip
In section \ref{Partitions and tableaux}, we fix the notation concerning
partitions, compositions, tableaux and so on, that we shall adhere to throughout the paper. 
Next, in section \ref{sectionYH} we fix the basic notation concerning the Yokonuma-Hecke algebra
$\YY $ of type $A$.  
There are several slightly different presentations
of $ \YY$ in the literature, but we use the presentation 
introduced in \cite{ChlouPouchin} which is 
the best adapted to KLR-theory.
In section \ref{sectionYH} we also fix some notation concerning set partitions $ \SetPar_n $, 
and ordered set partitions $ \OrdSetPar_n $, on $ \{1,2,\ldots, n \} $.

\medskip
In section \ref{idempotents in YY} we recall the Jucys-Murphy
$   {\mathcal L}_{d,n}=  \{X_1, \ldots, X_n, t_1,\ldots, t_n\}  $ for $ \YY$ that were
introduced \cite{ChPoAn} and were shown in \cite{ERH} to be a separating family of $\JM$-elements
for $ \YY$. We use them to construct idempotents in $ \YY $, associated with set partitions.

\medskip
Section \ref{JM for YYsection} is a key section of our paper. 
We develop the theory of 
$ \JM$-elements $   {\mathcal L}_{d,n} $
within the framework of \cite{Mat-So}, including 
contents $ c_{\bs} $, idempotents $ E_{\bs} $, the seminormal basis $ \{ f_{\bs \bt} \} $ and so on, for 
the generic Yokonuma-Hecke algebra $ \YYfractionfield$. In the non-semisimple case we consider
class idempotents $ E_{ \TT} $ and show that they coincide with the generalized eigenspace projectors
for $   {\mathcal L}_{d,n} $, generalizing a result of Hu-Mathas for the cyclotomic Hecke algebra, see
\cite{hu-mathas}. We recall the intertwiners $ \F_a $ introduced in \cite{Ro},
and show in Theorem \ref{commutation intertwiner} and 
Theorem \ref{commutation intertwiner setpar} a series of commutation formulas involving them. 
These properties were needed for 
the proof of $ \Hecken \cong R_n(\Gamma_{e} ) $ when $ d=1$ and for 
$ \YY \cong  \RR$, 
but they were not all proved in 
\cite{brundan-klesc} and \cite{Ro}.  
It is one of the purposes of the present paper to fill this gap which we do  
relying on the seminormal basis. 

\medskip
In section \ref{bt-algebra} we first recall the bt-algebra $ \E $. Next, 
in Theorem \ref{ortoBT} we use M\"obius inversion on the poset $ \SetPar_n$ to give
a new presentation for $ \E$, with associated diagram calculus that we present in
\eqref{illustrate1}--\eqref{illustrate3}. Partly inspired by Theorem \ref{ortoBT} we
then introduce the ordered bt-algebra $\Eord$, with associated diagram calculus.

\medskip
In section \ref{qHeckealgebras} we recall the KLR-algebra
$ R_n(\Gamma_{e, d}) $, considered in \cite{Ro}, and in section 
\ref{qhecke RE} we introduce the algebras $ \RE $ and $ \REord$ with associated
diagram calculi.
In the final section \ref{isomorphismProof} we prove our main Theorem
\ref{twomaintheoremsX} containing the isomorphisms
$  \E \cong \RE $ and $ \Eord \cong \REord$.

\medskip
\medskip
\noindent
    {\bf Acknowledgements}
The author wishes to express his gratitude to the 
organizers of
VI Encuentro de \'Algebra y Teor\'ia de Nudos
for giving him 
the opportunity to present this work.
He also wishes to thank D. Lobos for very useful conversations
related to it.

\section{Partitions and tableaux}\label{Partitions and tableaux}
In this section we fix some basic notions concerning the combinatorics of
partitions and tableaux.

\medskip
Let $ n $ be a non-negative integer.
A {\it composition} of $n $ of {\it length} $l$ is a finite sequence
of non-negative integers 
$ \lambda = (\lambda_1, \lambda_2, \ldots, \lambda_l) $ 
with sum $ n$. 
It is called a 
{\it partition} of $ n $
if 
$ \lambda_1 \ge \lambda_2 \ge \cdots \ge \lambda_l \ge 1   $.
We denote by $ \Comp_n^{l} $ the set of compositions of $n$ of length $  l $ and
by $ \Par_n $ the set of  partitions of $n$. 
For $ n= 0 $ we use the convention that $ \Par_0 = \{\emptyset\} $
and define 
$ \Par = \bigcup_{n=0}^{\infty} \Par_n $. For $ \lambda \in \Comp_n^l $ we define the {\it order}
$ | \lambda | $ of
$ \lambda $
via $ | \lambda | = n$.

\medskip
We identify compositions of $ n $ with their {\it Young diagrams}. For  
$ \lambda = (\lambda_1, \lambda_2, \ldots, \lambda_l) \in \Comp_n^l $ the Young diagram of $ \lambda $ consists
of $ l $ left aligned rows
of {\it nodes }(boxes), with the topmost row containing $ \lambda_1 $ nodes, the second topmost row
containing $ \lambda _2 $ nodes, and so on. For example, if $ \lambda = (5,3,2) \in \Par_{10}
\subseteq \Comp_{10}^3$ 
we have
\begin{equation}\label{firstpartition}
\lambda =   \raisebox{-.45\height}{\includegraphics[scale=0.7]{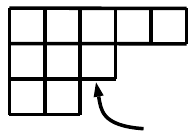}} \raisebox{-12\height}{.}
\end{equation}
We use matrix convention to label the nodes of $ \lambda$. For
example, for $ \lambda $ as in \eqref{firstpartition}, the node
indicated with an arrow is labeled $ \lambda[2,3]$.

\medskip
Let $ d $ be a positive integer. Then a {\it $ d$-multipartition} of $ n $,
or simply a {\it multipartition} if confusion is not possible, 
is a $ d $-tuple $ \blambda=
(\lambda^{(1)}, \lambda^{(2)}, \ldots, \lambda^{(d)}) $
such that $ \lambda^{(i)} \in \Par_{n_i } $ and such that $ n_1 +n_2 + \ldots + n_d = n $. 
We call $ \lambda^{(i)} $ the $i$'th {\it component} of $ \blambda$.
The set of $ d$-multipartitions of $ n $ is denoted by
$ \Par_{d, n} $. A $ d$-multipartition 
$ \blambda= (\lambda^{(1)}, \lambda^{(2)}, \ldots, \lambda^{(d)}) $ in 
$ \Par_{d, n} $ is identified with its {\it Young diagram} which by definition is the corresponding
tuple of Young diagrams for the $ \lambda^{(i)}$'s.
For example, for $ \blambda = ( \lambda^{(1)}, \lambda^{(2)}, \lambda^{(2)}) \in
\Par_{3,11} $ with $  \lambda^{(1)} = (3,2) $, $  \lambda^{(2)} = \emptyset $ and
$  \lambda^{(3)} = (2,2,2) $ we have
\begin{equation}\label{multipar}
\blambda =   \left(\raisebox{-.45\height}{\includegraphics[scale=0.7]{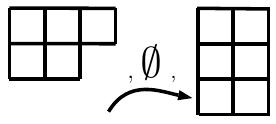}}\right)
 \raisebox{-12\height}{.}
\end{equation}
The nodes of $ \blambda$ are labeled by triples, via $ \blambda[r,c, p] 
= \lambda^{(p)}[r, c] $. For example, for $ \blambda $ as in 
\eqref{multipar}, the node indicated with an arrow has label $ \blambda[3,1,3]$.
For $ \blambda =( \lambda^{(1)}, \lambda^{(2)}, \ldots , \lambda^{(d)}) 
\in \Par_{d,n} $ we define the {\it order} $ ||\blambda ||$ of $ \blambda $ via
$ ||\blambda || = (| \lambda^{(1)}| , | \lambda^{(2)}| , \ldots, | \lambda^{(d)}| ) $. 
Note  that $  ||\blambda || \in \Comp_n^d$.

\medskip
For $ \lambda \in \Comp_n^l $ a $\lambda$-{\it tableau} $ \T $ is a filling
of the nodes of $ \lambda $ using the numbers $ \{1,2,\ldots, n \} $, each number once. 
A $\lambda $-tableau $ \T $ is called {\it row standard (resp. column standard)} 
if the numbers $ \{1,2,\ldots, n \} $ appear increasingly along the rows (resp. columns) of $ \T$
and it is called {\it standard} if it is both row and column standard. 
The set of (row standard, column standard, standard) $ \lambda $-tableaux is denoted 
$ \tab(\lambda) $ ($ \Rstd(\lambda)$, $ \Cstd(\lambda)$, $ \std(\lambda)$).
For example, for $ \lambda = (5,0,2,2) $ and $ \mu =(5,3,2) $ we have 
\begin{equation}
  \raisebox{-.45\height}{\includegraphics[scale=0.7]{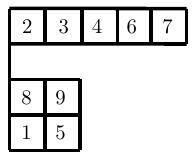}} \in \Rstd(\lambda), \, \, \, \, 
     \raisebox{-.45\height}{\includegraphics[scale=0.7]{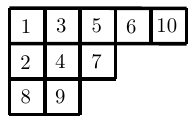}} \in \std(\mu). 
\end{equation}
The set of (row standard, column standard, standard) tableaux of {\bf partitions} of $ n $
is denoted $ \tab_n$ ($\Rstd_n$, $\Cstd_n$, $ \std_n$). We
set $ \tab = \bigcup_{n=0}^{\infty} \tab_n $ ($ \Rstd = \bigcup_{n=0}^{\infty} \Rstd_n $, 
$ \Cstd = \bigcup_{n=0}^{\infty} \Cstd_n $, $ \std = \bigcup_{n=0}^{\infty} \std_n $).
The set of (row standard, column standard, standard) tableaux of {\bf compositions} of $ n $
of length $ l $
is denoted $ \tab_n^{\Comp,l}$ ($\Rstd_n^{\Comp,l}$, $\Cstd_n^{\Comp,l}$, $ \std_n^{\Comp,l}$).

\medskip
For $ \blambda \in \Par_{d,n}$, a $ \blambda$-multitableau
$ \bT  $ is a filling
of the nodes of the components of $ \blambda $ using the numbers $ \{1,2,\ldots, n \} $, each number once.
The restriction of $ \bT $ to $ \lambda^{(i)} $ is the $i$'th {\it component} $ \T^{(i)} $ of $ \bT$, 
and therefore $ \bT $ may be viewed as the $d$-tuple $ \bt= (\T^{(1)}, \T^{(2)}, \ldots, \T^{(d)}) $. If $ a $ appears in $  \T^{(i)} $ we
write $ \pos_{\bT}(a) = i$. 

\medskip
A $\blambda $-tableau $ \bT $ is called {\it row standard (column standard)} 
if the numbers $ \{1,2,\ldots, n \} $ appear increasingly along the rows (columns) of
each component $ \T^{(i)}$ of $ \bT$, 
and it is called {\it standard} if it is both row and column standard. 
The set of (row standard, column standard, standard) $ \blambda $-tableaux is denoted 
$ \tab(\blambda) $ ($ \Rstd(\blambda)$, $ \Cstd(\blambda)$, $ \std(\blambda)$).
For example, for $ \blambda $ as in \eqref{multipar}
we have 
\begin{equation}\label{1.4}
\bT = \left(  \raisebox{-.4\height}{\includegraphics[scale=0.7]{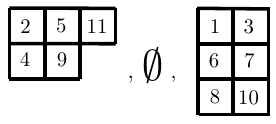}}\right) \in \std(\blambda)
\end{equation}
and $ \pos_{\bT}(5) = 1, \,  \pos_{\bT}(3) = 3 $.
The set of (row standard, column standard, standard) $d$-multitableaux 
of $d$-multipartitions of $ n$ is denoted $ \tab_{d,n}$, $\Rstd_{d,n}$, $\Cstd_{d,n}$, $ \std_{d,n}$.

\medskip
Suppose that $ \lambda \in \Par_n $ and $  \mu \in \Par_m$ with
$ \lambda = (\lambda_1, \lambda_2, \ldots, \lambda_{k} )$ and
$ \mu = (\mu_1, \mu_2, \ldots, \mu_l) $. 
Then we write $ \lambda \unlhd \mu $ if $ \sum_{i=1}^a \lambda_i  \le \sum_{i=1}^a \mu_i  $
for $ a = 1, 2, \ldots, {\rm max } (k,l) $
where we extend $ \lambda $ or $ \mu $ to the right with $ \lambda_i =0 $ or $ \mu_i =0 $
so that the sums on both sides make sense. 
This is the {\it dominance order} on $ \Par$.
It gives rise to the dominance order on $ \std$ as follows. For $ \T \in \tab(\lambda) $ 
define $ shape(\T) = \lambda$.
For $ \T \in \std_n $ and $ a \in \{1,2,\ldots, n\} $ let 
$ \T_{ \le a } \in \std_a $ be the standard tableau obtained from $ \T $ by deleting
the nodes containing numbers strictly larger than $ a $. For $ \s, \T \in \std $ we 
then define $ \s \unlhd \T $ if $ shape(\s_{ \le a }) \unlhd  shape(\T_{ \le a })  $ for all $ a $. 

\medskip
Suppose that $ \blambda, \bmu \in \Par_{d,n} $, with
$ \blambda = (\lambda^{(1)},\lambda^{(2)} , \ldots, \lambda^{(d)} )$ and
$ \bmu = (\mu^{(1)},\mu^{(2)} , \ldots, \mu^{(d)} )$. Then
the dominance order on $ \Par_{d, l} $ 
is given by $ \blambda \unlhd \bmu $ if 
$ \lambda^{(i)} \unlhd \mu^{(i)} $ for all $ i $. This is the dominance order on
$ \Par_{d, l} $ that we shall use in the present paper, but there are other interesting
orders on $ \Par_{d, l} $ that are also referred to as \lq dominance order\rq\ $\!\!$ in the literature. 

\medskip
For $ \bT \in \tab(\blambda)  $
we define $ shape(\bT) = \blambda $, and for $ \bT \in \std_{d, n} $ 
and $ 1 \le a \le n $ we let $ \bT_{\le a } \in \std_{d, a } $ be the multitableau obtained from 
$ \bT$ by deleting all nodes containing numbers strictly greater than $ a $.
Then the dominance order on $ \std_{d,n} $ is given by $ \bs \unlhd \bT $ if
$ shape(\bs_{ \le a} ) \unlhd shape(\bT_{ \le a} ) $ for all $ 1 \le a \le n $.

\medskip
Let $ \Si_n $ be the symmetric group consisting of the bijections of
$ \{ 1,2, \ldots, n \} $.
It acts on the left on 
$ \{ 1,2, \ldots, n \} $ via $ (\sigma, i ) \mapsto \sigma(i) $, where $ \sigma \in \Si_n $
and $ i \in \{ 1,2, \ldots, n \} $.
It is a Coxeter group on $ S = \{ \sigma_a | a =1,2,\ldots, n-1 \} $
where $ \sigma_a $ is the simple transposition $  \sigma_a = (a, a+1)$.
The $\Si_n$-action on $ \{ 1,2, \ldots, n \} $ induces a left $ \Si_n $-action on $ \tab_n$, 
written $ (\sigma, \T) \mapsto \sigma \T $, and it also induces a left $ \Si_n $-action on $ \tab_{d,n}$, 
written $ (\sigma, \bt) \mapsto \sigma \bt $, 
by 
permuting entries.

\section{The Yokonuma-Hecke algebra $ \YY  $. }\label{sectionYH}
In this section we recall the 
Yokonuma-Hecke algebra $ \YY $ and some of its fundamental properties.

\medskip
The Yokonuma-Hecke algebra was introduced in the sixties as a generalization of
the Iwahori-Hecke algebra of type $ A$. There are several slightly different
presentations of $ \YY$, but we shall use the
following one from \cite{ChlouPouchin} and  \cite{Ro}.

\begin{definition}
  \normalfont
Let $ \Kdomain $ be a domain and let $ d, n$ be positive integers, such that $ d $ is coprime to $ {\rm char }\,  \Kdomain$.
Let $ q \in \Kdomain^{\times}$. 
Then the Yokonuma-Hecke algebra $\YY =  \YYresiduefield $ is the associative $ \Kdomain$-algebra generated by
  $ \{g_1, g_2,\ldots, g_{n-1} \} $ and $ \{ t_1, t_2, \ldots, t_n\} $ subject to the following relations
	\begin{align}
\label{eq one}	t_a^d	& =  1   & &      \\
\label{eq two}	t_a t_b  & = t_b t_a &     &     \\
\label{eq three} t_a g_b & = g_b t_{ \sigma_b(a)}  &  &        \\
\label{eq four} g_a^2  & = q +  (q-1) e_a g_a   &  &        \\
\label{eq five} g_a g_b g_a   & = g_b g_a g_b   &   &    \mbox{if }   | a-b | = 1   \\
\label{eq six} g_a g_b    & = g_b g_a       &   &  \mbox{if } | a-b | > 1   
	\end{align}
where $ e_a = \frac{1}{d}\sum_{r=0}^{d-1} t_{a}^{r} t_{a+1}^{-r} $. 
\end{definition}

\medskip
    {\it Throughout the paper, unless stated otherwise we shall always assume that $ q \neq 1 $, corresponding to
      the non-degenerate case, although some of our results also hold for $ q =1$.}

\medskip    
Let $ \Z_{(d)}  \subseteq \QQ $ be the
localization of $ \Z $ at $ d $, that is $ \Z_{(d)} = S^{-1} \Z $ where $ S=\{1,d, d^2, \ldots \}$,
and define $ \YYZ  = \mathcal{Y}^{\, \Z_{(d)}[ \q,\q^{-1}]}_{d, n}(\q) $ where $ \q $ is an indeterminate. Then 
any 
$ \YYresiduefield $ can be obtained from $  \YYZ $ via base change.

\medskip
It follows from the relations \eqref{eq one} -- \eqref{eq six} 
that $ e_a $ is an idempotent in $ \YY $ and that 
$ g_a $ is invertible, with inverse $ g_a^{-1} = 
q^{-1} g_a +(q^{-1}-1) e_a $. For $1 \le a < b \le n$ we then define $ e_{ab} \in \YY$ via
\begin{equation}\label{e_{ab}}
e_{ab} = g_{b-1} g_{b-2}\cdots g_{a+1}e_a g_{a+1}^{-1}\cdots g_{b-2}^{-1} g_{b-1}^{-1}
\end{equation}
where by convention $ e_{a,a+1} = e_a  $.

\medskip
For $1 \le a < b \le n$ we set
$e_{ba} = e_{ab} $. 
Then for all $ a, b $ we have 
\begin{equation}\label{e_{ab}T}
e_{ab}=  \frac{1}{d}\sum_{r=0}^{d-1} t_{a}^{r} t_{b}^{-r}.
\end{equation}
It follows from 
\eqref{eq one}--\eqref{eq six}
that the $e_{ab}$'s are commuting idempotents.

\medskip
The combinatorics of set partitions shall be of importance throughout the present paper. 
Let $ \SetPar_n $ be the {\it set partitions} on $ \{ 1,2 \ldots, n \} $, 
that is the set of equivalence relations
on $ \{ 1,2 \ldots, n \}$.
To be precise, an element of $ \SetPar_n $ is a set of the form $ A=\{A_1, A_2, \ldots, A_k \} $
where the $ A_i$'s are non-empty, disjoint subsets of $ \{ 1,2 \ldots, n \} $
such that $ \bigcup_i A_i = \{ 1,2 \ldots, n \} $. For $ A = \{A_1, A_2, \ldots, A_k \} \in \SetPar_n$, 
the $ A_i$'s are called the {\it blocks} of $ A$. 
For example, for $ n=3 $ we have $  A = \{ \{1,2\}, \{3\} \} \in \SetPar_3 $ with
blocks $ A_1 = \{1,2\} $ and $ A_2 = \{3\} $.

\medskip
Let $ b_n = | \SetPar_n | $ be the number of set partitions of $ n $, that is $ b_n $ is 
the $n$'th {\it Bell number}. For example $ b_3 = |  \SetPar_3 | =5 $ since
\begin{equation}
  \SetPar_3 =
  \big\{ \{ \{1,2,3\}\},  \{ \{1,2\}, \{3\} \}, \{ \{1,3\}, \{2\}\}, \{ \{2,3\}, \{1\}\},  \{ \{1\}, \{2\},\{3\} \} \big\}.
\end{equation}  
The sequence $ \{ b_1, b_2, b_3, \ldots, \} $
is A000110 in OEIS.

\medskip
For $A \in \SetPar_n $ we denote by $ \sim_A $ the equivalence relation associated with $A$, 
that is $ a \sim_A  b $ if and only if $ a $ and $ b $ belong to the same block of $ A$.
We then define $ e(A) \in \YY $ via 
\begin{equation}\label{e_A}
e(A) = \prod_{\substack{1\le a,b\le n \ \\ a \, \sim_A \,b}} e_{ab}.
\end{equation}
The product in \eqref{e_A} can be taken in any order, since the $ e_{ab}$'s
commute. 
Moreover, the $ e(A)$'s are commuting idempotents in $ \YY$,
although they are not orthogonal. 

\medskip
The set of {\it ordered set partitions}, denoted $ \OrdSetPar_n $, is also of importance in
the present paper. By definition, an element of $  \OrdSetPar_n $ is a tuple
$ A^{\ord} = (A_1, A_2, \ldots, A_k ) $ such that 
$  A=  \{ A_1, A_2, \ldots, A_k \}  \in \SetPar_n$. 
The $ A_i$'s are
called the 
{\it blocks} for $ A^{\ord}$. For example we have 
\begin{equation}\label{ordered set partitions}
  \OrdSetPar_3 =
\begin{array}{l}
  \big\{ ( \{1,2,3 \}) ,  ( \{1,2\}, \{3\} ), ( \{3\}, \{1,2\} ),  ( \{1,3\}, \{2\}), ( \{2\}, \{1,3\}),  \\ 
 ( \{2,3\}, \{1\}), ( \{1\}, \{2,3\}),  ( \{1\}, \{2\},\{3\} ), 
  ( \{1\}, \{3\},\{2\} ),    \\   ( \{2\}, \{1\},\{3\} ),
    ( \{2\}, \{3\},\{1\} ), ( \{3\}, \{1\},\{2\} ), ( \{3\}, \{2\},\{1\} ) \big\}.
\end{array}
\end{equation}  

\medskip

Let $ f_n =   | \OrdSetPar_n  | $. Then $f_n $ is the $n$'th {\it Fubini number}, or 
the $n$'th {\it ordered Bell number}.
The sequence $ \{f_0,  f_1, f_2, f_3, \ldots \} $ is A000670 in OEIS.

\section{Idempotents in $ \YY$.}\label{idempotents in YY}
Suppose that $ \Kdomain $ is a field for which $ \YY = \YYresiduefield  $ exists and 
assume that there exists a $d$'th primitive root of unity $ \xi \in \Kdomain $. 
In this section we shall show how one can use the $ e(A)$'s to 
obtain orthogonal idempotents in $ \YY $.

\medskip
Let $  {\rm char }_q \, \Kdomain $ be the {\it quantum characteristic} of $q$, that is 
$ {\rm char }_q (\Kdomain) = {\rm min} \, A_q $ where 
\begin{equation}\label{quantum}
A_q =  \{   e \in {\mathbb N}  \, | \, 1 + q +q^2 + \ldots + q^{e-1} =0  \} 
\end{equation}
with the convention that $ {\rm min} \, \emptyset = 0 $.
Thus, by our standing hypothesis $ q \neq 1 $ we have
\begin{equation}
{\rm char }_q (\Kdomain)= {\rm min} \{   e \in {\mathbb N}  \, | \, q^e =1  \}.
\end{equation}

\medskip
Define $ I = {\mathbb Z }/ {\rm char }_q (\Kdomain)   {\mathbb Z }  $, $ J = {\mathbb Z }/d  {\mathbb Z } $ and $ K = I \times J$. 
Let $ {\II } = I^n $, $ {\JJn } = J^n $ and $ {\KKn } = K^n $.
We shall write elements of $ \II $ and $ \JJn $ in the form $ \ii=(i_1, i_2, \ldots, i_n) $
and $ \jj = (j_1, j_2, \ldots, j_n) $, whereas elements of $ \KKn $ shall be written in the form 
$\kk = (k_1, k_2, \ldots, k_n)  = ( (i_1, j_1), (i_2, j_2), \ldots, (i_n, j_n) ) $.
We shall sometimes write $ \ii_a $ for $ i_a $, and $ \jj_a $ for $ j_a $. 
We shall also sometimes use the notation $   k_{a, i}  $ and $   k_{a, j}  $
for the two coordinates of $ k_a$, that
is $ k_a = ( k_{a, i},  k_{a, j}) $. 
We set $ \ii_{\kk}=(k_{1,i}, k_{2, i}, \ldots, k_{n, i}) $ and 
$ \jj_{\kk}=(k_{1, j}, k_{2, j}, \ldots, k_{n, j}) $. For example, if $ \kk= ( (1,2), (2,3), (0,1), (2,2)) $
we have $  \ii_{\kk} = (1,2,0,2) $ and $  \jj_{\kk} = (2,3,1,2) $.

\medskip
There is a natural action of $ \Si_n $ on $ \II $, $ \JJn $ and $ \KKn $ by place permutation, given by 
\begin{equation}\label{placeperm}
 \begin{aligned}
   \sigma_a \ii &= (i_1, i_2, \ldots, i_{a+1}, i_a, \ldots, i_n),  & \sigma_a \jj = (j_1, j_2, \ldots, j_{a+1}, j_a, \ldots, j_n) \\
 \sigma_a \kk &= (k_1, k_2, \ldots, k_{a+1}, j_a, \ldots, k_n) & 
\end{aligned}
\end{equation}
for $ 1 \le a < n$.   

\medskip
Following \cite{ChPoAn}, \cite{WeidengCuiJinkuiWan} and \cite{ERH}
we define
$ \{X_1, X_2, \ldots, X_n \} \subseteq  \YY $ via $ X_1 = 1 $ and
recursively $ X_{a+1} =  q^{-1} g_a X_a g_a $. Then 
\begin{equation}\label{JM for YY}
  {\mathcal L}_{d,n}=  \{X_1, \ldots, X_n, t_1,\ldots, t_n\} \subseteq
  \YY
\end{equation}
is the set of {\it Jucys-Murphy elements} for
$ \YY $. 

\medskip
They form a commuting family of elements in $ \YYresiduefield $. Moreover, the 
eigenvalues for the action of the $ X_a$'s 
in $ \YYresiduefield $ by left multiplication are of the form $ q^i $ for some $ i \in I$,
see  Corollary \ref{eigenvalues} or \cite{WeidengCuiJinkuiWan},
and by \eqref{eq one} and \eqref{eq two} the actions of the $ t_a $'s in $ \YY $
are simultaneously diagonalizable. 
Hence there exist unique eigenspace projector idempotents 
$ e(\kk) \in   \YYresiduefield $ such that 
\begin{equation}\label{eigenspace}
e(\kk) \YYresiduefield=  \{ x \in \YYresiduefield \, | \, (X_a-q^{k_{a,i}})^N x = (t_a-\xi^{k_{a,j}})x = 0 \, \mbox{ for } N >>0 \} .
\end{equation}
They are orthogonal and complete, that is
\begin{equation}\label{regularmodule}
  e(\kk) e(\kk^{\prime}) = \delta_{ \kk \kk^{\prime} } e(\kk) \, \, \, \, \, \, \, \, \, \, \, \, \, \, \,
  \mbox{and}  \, \, \, \, \, \, \, \, \, \, \, \, \, \, \,
  1 = \sum_{ \kk \in \KKn } e(\kk)
\end{equation}
where $  \delta_{\kk \kk^{\prime}} $ is the Kronecker delta, although they are 
not all non-zero in general.


\medskip
For $ \ii \in \II $ and $ \jj \in \JJn$ we define $ \KKii, \KKjj  \subseteq \KKn $ via 
\begin{equation}
  \KKii  = \{ \kk \in \KKn \,  | \, \ii_{\kk} = \ii  \}, \, \, \,
\KKjj  = \{ \kk \in \KKn \,  |   \, \jj_{\kk} = \jj  \}
  \end{equation}
and so we obtain idempotents $ e(\ii), e(\jj) \in \YY $ as follows
\begin{equation}\label{def e(i)}
e(\ii) = \sum_{ \kk \in \KKii } e(\kk), \, \, \, 
e(\jj) = \sum_{ \kk \in \KKjj } e(\kk). \, \, \, 
  \end{equation}

By definition
\begin{equation}\label{eigenspaceIJ}
\begin{aligned}
  e(\ii) \YYresiduefield &=  \{ x \in \YYresiduefield \, | \, (X_a-q^{i_{a}})^N x  = 0 \, \mbox{ for } N >>0 \},  \\
  e(\jj) \YYresiduefield &=  \{ x \in \YYresiduefield \, | \, (t_a-\xi^{j_{a}})x = 0  \} \\
  \end{aligned}
\end{equation}
and each of the sets $ \{ e(\ii) \,| \, \ii \in \II \} $ and $ \{ e(\jj) \,| \, \jj \in \JJn \} $ is a complete set of orthogonal
idempotents.

\medskip

We have the following Lemma.
\begin{lemma}\label{description}
Let $ e_{ab} \in  \YY$ be the element from \eqref{e_{ab}}.
Then  
\begin{equation}\label{thefollowingdescription}
e_{ab} = \sum_{\substack{ \jj\in \JJn \\ j_{a} = j_{b}}} e(\jj) .
\end{equation}  
\end{lemma}
\begin{dem}
Since $ 1 = \sum_{ \jj \in \JJn} e(\jj) $ we get from \eqref{regularmodule} that 
    \begin{equation}
t_a = \sum_{ \jj \in \JJn} t_a e(\jj) = \sum_{ \jj \in \JJn} \xi^{ j_a}    e(\jj).
  \end{equation}
From this we deduce via \eqref{e_{ab}T} that 
  \begin{equation}
 e_{ab} = \frac{1}{d}\sum_{r=0}^{d-1} t_{a}^{r} t_{b}^{-r}
= \sum_{ \jj \in \JJn}\frac{1}{d} \sum_{r=0}^{d-1}\xi^{ r(j_{a}- j_{b}) }e(\jj).
  \end{equation}
  On the other hand, since $ \xi $ is a primitive $d$'th root of unity we have
  \begin{equation}
  \sum_{r=0}^{d-1}\xi^{ r(j_{a}- j_{b}) } =
 \left\{\begin{aligned}
 0  &  & \mbox{if } & j_{a} \neq j_{b}   \\
  d  &  & \mbox{if } & j_{a}= j_{b}        
\end{aligned} \right. 
\end{equation}  
  and so we arrive
  at $e_{ab} = \mathlarger{\sum}\limits_{\substack{ \jj \in \JJn \\ j_{a} = j_{b}}   } e(\jj) $ as claimed.
\end{dem}

\medskip

\begin{corollary}\label{description1}
Let $ \{ e(A) \, | \, A \in \SetPar_n \} \subseteq \YY$ be the elements defined in 
\eqref{e_A}. Then 
\begin{equation}\label{description1A}
e(A) = \sum_{a \sim_A b \, \Longrightarrow \,  j_{a} = j_{b}  } e(\jj) .
\end{equation}  
\end{corollary}
\begin{dem}
  This follows immediately from the definition of $ e(A) $ in \eqref{e_A} and Lemma \ref{description}. 
\end{dem}

\medskip
Note that \lq$\Longrightarrow \! $\rq\ in the expression for
$ e(A) $ in \eqref{description1A} cannot be replaced by \lq$ \Longleftrightarrow \!$\rq\ .
This motivates the following definitions.
Let $ A \in \SetPar_n$. Then we define $ \JJ_{A}  \subseteq \JJn $ via
\begin{equation}
\JJ_{A}  = \{ \jj \in \JJn \,  | \, j_{a} = j_{b} \Longleftrightarrow 
  a \sim_A b   \}.
  \end{equation}

\medskip
Note that in general it is possible that $\JJ_{A} = \emptyset $, but 
for $ d \ge n $ we have 
$ \JJ_{A} \neq \emptyset $ for any $ A \in \SetPar_n $. For example, for $ d=n=3 $
and  $A=\{ \{1, 2\}, \{3\} \} $
we have 
\begin{equation}
  \JJ_{A} = \{  (0,0,1), (0,0,2), (1,1,0), (1,1,2),(2,2,0), (2,2,1)   \}.
\end{equation}

\medskip
For $ A \in \SetPar $ we then define
\begin{equation}\label{overline e(A)first}
  \e( A) = \sum_{ \jj \in   \JJ_{A} } e(\jj)  \in \YY .
\end{equation}

Recall next that there is a natural poset structure $ (\SetPar_n, \subseteq) $. For 
$ A = \{ A_1, \ldots, A_k \}, B = \{ B_1, \ldots, B_l \}  \in \SetPar_n $, it is defined by $ A \subseteq B $ if
and only each $ B_i $ is a union of certain $ A_i$'s.
Note that 
\begin{equation}\label{beforeMoeb}
e(A) e(B) = e(B) \, \, \, \mbox{ if   } \, A \subseteq B 
\end{equation}
as one sees from \eqref{e_A}. 

\medskip
The M\"obius inversion function $ \mu $ associated 
with $ (\SetPar_n, \subseteq)$ is given by
\begin{equation}\label{MobiusFirst}
\mu(A,B) =  \left\{\begin{aligned}
 & (-1)^{r-s}\prod_{i=1}^{r-1} {i!}^{r_{i+1}}   &  & \mbox{if }  A \subseteq B   \\
 &  0  &  & \mbox{otherwise }         
\end{aligned} \right. 
\end{equation}
where $ r$  and $ s$  are the number of blocks of $A$ and $B$ respectively, and where $r_i$ is the
number of blocks of $B$ that contain exactly $i$ blocks of $A$.

\medskip
Let 
$ {\mathcal B }_{ n} $ be the subalgebra of $ \YY $ generated by $ \{ e(A)  \, | \, A \in \SetPar_n \} $. 
The following is the main Theorem of this section.
\begin{theorem}\phantomsection\label{teorem 6A}
  \begin{description}
  \item[a)] $ \{ \e(A)  | \,  A \in \SetPar_n \} $ is a set
    of orthogonal idempotents for $ \YY$. It is complete, in other words $ \sum_{ A \in \SetPar_n } \e(A) = 1 $. 
\item[b)] $ \{ \e(A) \neq 0 \, | \, A \in \SetPar_n \} $ is a basis for $  {\mathcal B }_n $. 
\end{description}  
  \end{theorem}
\begin{dem}
  We first show $   { \bf a)} $. Surely $ \e(A) $ is an idempotent since it is a 
  sum of orthogonal idempotents. Moreover, if $ A \neq A_1 $ the
  sets $ \JJ_{A} $ and $ \JJ_{A_1} $ are
  disjoint, and therefore $ \e( A) $ and $ \e( A_1) $ are orthogonal idempotents. From this we deduce that
  the $ \e(A) $'s are orthogonal idempotents, and they are complete in view of
\eqref{def e(i)}.

\medskip  
To show $   { \bf b)} $, we let $ \overline{{\mathcal B}}_n $ be the subalgebra generated by
the $ \e(A)$'s. We claim that $ {\mathcal B}_n =\overline{{\mathcal B}}_n $
from which 
$   { \bf b)} $ follows via $   { \bf a)} $. 
Here the inclusion $ {\mathcal B }_n \subseteq \overline{{\mathcal B }}_n$
follows from the identity    
\begin{equation}\label{Mobius}
 e(A) =  \sum_{\substack{  B \in \SetPar_n \\  A \subseteq B}} \e(B)  
\end{equation}  
which is a consequence of Lemma \ref{description} and \eqref{overline e(A)first}.

\medskip
To show the other inclusion we use M\"obius inversion on \eqref{Mobius} to arrive at 
\begin{equation}\label{firstMobius inversion}
 \e(A) =  \sum_{\substack{  B \in \SetPar_n \\  A \subseteq B}} \mu(A,B) e(B)
\end{equation}
where $ \mu $ is the M\"obius inversion function given in \eqref{MobiusFirst}. 
Since $ \mu(A,B) \in \Z$, 
the other inclusion $ {\mathcal B }_n \supseteq \overline{{\mathcal B }}_n$ now follows.  
\end{dem}

\medskip
The left $ \Si_n $-action on $ \{1,2 \ldots, n \} $ extends canonically to a left $ \Si_n $-action
on $ \SetPar_n$, denoted $ (\sigma, A) \mapsto \sigma A$. Using it we get the following Lemma 
that we shall need later on.
\begin{lemma}\label{afirstlemmaYY}
For $ 1\le a \le n-1 $ and $ A \in \SetPar_n $ the following identities hold in $ \YY$.
  \begin{description}
\item[a)] $g_a e(A) = e(\sigma_a  A) g_a$
\item[b)] $ g_a \e(A) = \e( \sigma_a  A) g_a$
\end{description}
\end{lemma}
\begin{dem}
  $   { \bf a)} $ follows from \eqref{e_{ab}T} and the definition of $ e(A) $ in
  \eqref{e_A} and $   { \bf b)} $ follows from $   { \bf a)} $
  via \eqref{firstMobius inversion}. 
\end{dem}

\section{Jucys-Murphy elements for $ \YY $}\label{JM for YYsection}
In this section we return to the Jucys-Murphy elements $ {\mathcal L}_{d,n}=  \{X_1, \ldots, X_n, t_1,\ldots, t_n\} $ for
$ \YY $, that were introduced in \eqref{JM for YY}.
Set
\begin{equation}
L_a = \left\{ \begin{array}{lll} X_a  & \mbox{ if } & 1 \le a \le n \\ t_{a-n }  &\mbox{ if } & n+1 \le a \le 2n
 \end{array} \right.
\end{equation}
such that $ {\mathcal L}_{d,n}  = \{L_1, L_2, \ldots, L_{2n} \}$.
In \cite{ERH} it was shown that   
$ {\mathcal L}_{d,n} $ is a family of $\JM$-elements for $ \YY $ in the sense of Mathas, see \cite{Mat-So}.
This result is important for us, so let us explain its precise meaning.

\medskip
We start out by recalling the definition from \cite{GL} of a cellular algebra.

\begin{definition}\label{cellular} Let $\Kdomain$ be a domain. Suppose that $\A$ is a
$\Kdomain$-algebra which is free as a $\Kdomain$-module. Suppose that $(\Lambda,
\geq)$ is a poset and that for each $\lambda\in \Lambda$ there is a
finite indexing set $T(\lambda)$ and elements $c_{st}^{\lambda}\in
\A$ such that
\begin{equation}
  \mathcal{C}=\{c_{st}^{\lambda}\mid \lambda\in \Lambda \mbox{ and } s,t\in T(\lambda)\}
\end{equation}  
is a $\Kdomain$-basis for $\A$. Then the pair $(\mathcal{C},\Lambda)$ is called a 
{\it cellular basis} for $\A$ if
\begin{enumerate}\renewcommand{\labelenumi}{\textbf{(\roman{enumi})}}
\item The $\Kdomain$-linear map $*:\A\to \A$ given by
$(\mathop{c_{st}^{\lambda}})^*=c_{ts}^{\lambda}$ for all
$\lambda\in\textcolor{black}{\Lambda}$ and all $s,t\in T(\lambda)$
is an algebra anti-automorphism of $\A$.

\item For any $\lambda\in \Lambda,\; t\in T(\lambda)$ and $a\in \A$
there exists $r_{  a s  u} \in \Kdomain $ such that for all $s\in T(\lambda)$
\begin{equation}
 a  c_{st}^{\lambda}  \equiv \sum_{u \in T(\lambda)} r_{a s  u}c_{u  t}^{\lambda} \mod{\A^{\lambda}}
\end{equation}  
where $\A^\lambda$ is the $\Kdomain$-submodule of $\A$ with basis
$\{c_{uv}^{\mu}\mid \mu\in \Lambda,\mu>\lambda \mbox{ and }
u,v\in T(\mu)\}$.
\end{enumerate}
If $\A$ has a cellular basis we say that $\A$ is a \textit{cellular
algebra}.
\end{definition}

\medskip
In Theorem 33 of \cite{ERH} 
it was shown that $ \YY $ is a cellular algebra in the sense of Definition
\ref{cellular}.
Here 
$ \Lambda =  \Par_{d, n} $ is endowed with the dominance order explained in section
\ref{Partitions and tableaux}, and 
$T(\blambda) = \std(\blambda) $ for $ \blambda \in \Lambda $. 
At this point, however, we do not need the details of the construction of the cellular basis element
$ c_{\bs \bt} $ itself, 
for $ \bs, \bt \in  \std(\blambda) $.

\medskip
The following definition was first given by 
Mathas in \cite{Mat-So}.

\begin{definition}{\label{JMdef}}
Suppose that $ \A $ is a $ \Kdomain $-algebra which is cellular with cellular basis $\mathcal{C}=\{c_{st} \, | \,
 \lambda\in \Lambda,s,t\in T(\lambda)\}$. Suppose moreover that each $ T(\lambda) $ is
endowed with a poset structure with order relation $ \rhd_{\lambda} $.
Then we say that a commuting set
 $\mathcal{L}=\{L_1,\ldots,L_M\}\subseteq \A $ is a family of $\JM$-elements for $\A$
if it satisfies $ L_i^{\ast} = L_i $ for all $ i $ and if
 there exists a set
 $\{c_{s}(i) \, | \, t\in T(\lambda), 1\leq i\leq M\} \subseteq \Kdomain$, called the contents of $ \lambda $, such
 that for all $\lambda\in \Lambda$ and $t\in T(\lambda)$ we have 
\begin{equation}{\label{contents}}
  L_i a_{s t}=c_{s}(i)a_{s t}+
  \sum_{\substack{u\in T(\lambda) \\ u \rhd_{\lambda} s}}r_{i u }a_{u t} \mod A^{\lambda}
\end{equation}
$ \mbox{ for some }r_{i u }\in  \Kdomain$.
\end{definition}

\medskip
In Proposition 38 of \cite{ERH} it was shown that $ {\mathcal L}_{d,n} = 
 \{L_1, L_2, \ldots, L_{2n} \} = \{X_1, X_2, \ldots, X_n, t_1, t_2, \ldots, t_n \}$
is
a family of $ \JM $-elements for $ \YY $ in the sense of Definition {\ref{JMdef}}. 
Here the poset structure on $ \std(\blambda) $ is the dominance order from section
\ref{Partitions and tableaux}, whereas 
for $ \bs \in \std(\blambda) $ the content function $ c_{\bs}(i) $ is given by the following formula
\begin{equation}\label{contentfunction}
  c_{\bs}(a) =  \left\{\begin{aligned} & q^{c-r} && \mbox{if} && 1 \le a \le n \\
& \xi^{p} && \mbox{if} &&  n+1 \le a \le 2n
  \end{aligned}  \right.
\end{equation}
where $ a $ is positioned in the node $\blambda[r,c,p] $ of $ \bs$.

\medskip
We next give a series of Lemmas involving commutation rules for the $ X_a$'s
that shall be useful later on.

\begin{lemma}\label{firstusefullemma}
Let $ 1 \le a \le n $ and $ A \in \SetPar_n $. Then the following identities hold in $ \YY$. 
  \begin{description}
  \item[a)]     $   X_a e(A) = e(A) X_a. $ 
  \item[b)]  $    X_a \e(A) = \e(A) X_a. $
   \item[c)] Suppose $ 1 \le a < n $, $ 1 \le b \le n $ and that $ b \neq a, a+1$. Then $    g_a X_b  =  X_b g_a .$
  \end{description}
\end{lemma}
\begin{dem}
  Using $ X_a = q^{1-a} g_{a-1} g_{a-2} \cdots g_2 g_1^2 g_2 \cdots g_{a-2} g_{a-1} $ we get 
  $ {\bf a)} $ from Lemma \ref{afirstlemmaYY}. 
  To show $ {\bf b)} $ we use $ {\bf a)} $ together with \eqref{firstMobius inversion}.
  Finally, $ {\bf c)} $ is a consequence of the braid relations
\eqref{eq five} and \eqref{eq six}.
\end{dem}

\begin{lemma}\label{veryfirstusefullemma}
  Suppose that $ A \in \SetPar_n $ and that $a \sim_A (a+1) $. Then the following identities hold in 
  $ \YY $. 
  \begin{description}
  \item[a)]     $  g_a^2 \e(A) =  (q + (q-1) g_a) \e(A). $
\item[b)]   $     g_a X_{a+1} \e(A) = X_{a}( g_a  + (q-1) )                     \e(A). $
  \item[c)]  $     g_a X_{a} \e(A) = X_{a+1}( g_a  + (1-q) )                     \e(A). $
  \end{description}
\end{lemma}
\begin{dem}
Let $ A $ and $ a $ be as in the announcement and define 
$ A_a=\{ \{1\}, \{2\}, \ldots, \{a,a+1\}, \ldots, \{n-1\},  \{n\} \} \in \SetPar_n$. 
Then $ A_a \subseteq A $ and $ e_a =  e(A_a) $ and 
so \eqref{beforeMoeb} 
together with 
\eqref{Mobius}, gives us 
$ e_a \e(A) = \e(A)$. Hence $   { \bf a)} $ follows from the relation 
\eqref{eq four}.  But $   { \bf a)} $ immediately implies $   { \bf b)} $ and  $   { \bf c)} $ via 
the definition of $ X_a $.
\end{dem}

\begin{lemma}\label{usefullemma}
  Suppose that  $ a \not\sim_A (a+1) $. Then we have that
  \begin{description}
  \item[a)]     $  g_a^2 \e(A) = q \e(A). $
  \item[b)]  $     g_a X_{a+1} \e(A) = X_a g_a \e(A). $
\item[c)]  $ g_a X_a \e(A) =     X_{a+1} g_a \e(A)   . $
  \end{description}
\end{lemma}
\begin{dem}
  For $   { \bf a)} $ we use once again
  $ A_a=\{ \{1\}, \{2\}, \ldots, \{a,a+1\}, \ldots, \{n-1\},  \{n\} \} \in \SetPar_n$, such
  that $ e_a =  e(A_a) $. This time we have $ A_a \nsubseteq A$ and so
${\bf a)} $ of Theorem \ref{teorem 6A} and 
\eqref{Mobius} 
imply that $ e_a \e(A) = 0$. Hence, using relation 
\eqref{eq four}
we get 
  \begin{equation}
g_a^2 \e(A) = \big(q+(q-1)e_a g_a\big) \e(A)  = \big(q+(q-1)g_a e_a\big) \e(A)  = q \e(A)
  \end{equation}
  which is $   { \bf a)} $. To show $   { \bf b)} $ we combine $   { \bf a)} $ with ${ \bf b)}$ of Lemma 
  \ref{firstusefullemma} to arrive at 
\begin{equation}
  g_a X_{a+1} \e(A) =   g_a q^{-1} g_a X_a g_a  \e(A) =  g_a q^{-1} g_a   \e(\sigma_aA)  X_a g_a=
  q^{-1} g_a^2   \e(\sigma_aA)  X_a g_a=  \e(\sigma_aA)  X_a g_a = X_a g_a \e(A)
  \end{equation}
which is $   { \bf b)} $. Finally $   { \bf c)} $ is proved the same way as $   { \bf b)} $.
\end{dem}

\medskip
\begin{corollary}\label{firstcor}
  Let $ 1 \le a \le n-1$. Then $ g_a $ commutes with any symmetric polynomial
  in $ X_a $ and $ X_{a+1}$. 
\end{corollary}
\begin{dem}
From  
Lemma \ref{veryfirstusefullemma} and Lemma \ref{usefullemma} we have that 
    \begin{equation}
    g_a(X_a + X_{a+1}) \e(A)   = (X_a +X_{a+1}) g_a\e(A) 
    \end{equation}
    for all $ A \in \SetPar_n$ and so, in view
    of ${\bf a)} $ of Theorem \ref{teorem 6A}, we get that $ g_a $ commutes with $ X_a + X_{a+1} $.
    On the other hand, we also have that $  g_a $ commutes with $ X_a  X_{a+1} $ since
    \begin{equation}
\begin{aligned}      
      g_a (X_a X_{a+1} ) &=  g_a \left(q^{1-a} g_{a-1} g_{a-2} \cdots g_2 g_1^2 g_2 \cdots g_{a-2}  g_{a-1} \right)
      \left(q^{-a} g_{a}g_{a-1}   \cdots g_2 g_1^2 g_2 \cdots g_{a-1}  g_{a} \right) \\
&=   \left(q^{-a} g_a g_{a-1} g_{a-2} \cdots g_2 g_1^2 g_2 \cdots g_{a-2}  g_{a-1} g_a \right)
      \left(q^{1-a} g_{a-1}   \cdots g_2 g_1^2 g_2 \cdots g_{a-1}  \right)  g_{a} \\
&=  (X_a X_{a+1} )g_a .
 \end{aligned}      
    \end{equation}
Combining, we deduce that $ g_a $ commutes with all symmetric polynomials in $ X_a $ and $X_{a+1} $, as claimed.
\end{dem}

\medskip
Let once again $ \q $ be an indeterminate and set 
$ \OO = \Kdomain[\q, \q^{-1}]_{(\q -q) } $ that is 
$ \OO = \left\{ \frac{f(\q)}{g(\q)}  \in  \Kdomain[ \q, \q^{-1}] \,  \big| \, g(q) \neq 0 \right\} $.
Then $ \OO $ is a discrete valuation ring with maximal ideal $ (\q -q )\OO $,
satisfying $ \OO/(\q -q ) \OO = \Kdomain $, 
and with quotient field $ Q(\OO) = \Kdomain(\q) $. In other words, $ (\Kdomain(\q), \OO, \Kdomain ) $
is a modular system for $ \YY$.

\medskip
For 
the Yokonuma-Hecke algebra $ \YYOO$, defined over $ \OO $, the content function is given by
\begin{equation}
  \C_{\bs}(a) =  \left\{\begin{aligned} & {\q}^{\, c-r} && \mbox{if} && 1 \le a \le n \\
& \xi^{p} && \mbox{if} &&  n+1 \le a \le 2n
  \end{aligned}  \right.
\end{equation}
where $ a $ is positioned in the node labeled $\blambda[r,c,p] $ of $ \bs$. 
Note that $ \C_{\bs}(a) $ 
satisfies the
{\it separation condition} given by Mathas in Definition 2.8 of \cite{Mat-So}, since $ \q $ is an indeterminate. 
We may therefore
apply the results of section 3 of \cite{Mat-So} to $ \YYfractionfield $, and in particular we deduce that 
$ \YYfractionfield $ is semisimple.

\medskip
For $ a \in \{1, 2, \ldots, 2n \} $ we set $ {\mathcal C}(i) = \{ \C_{\T}(a) | \T \in \std_{d,n} \}$ and
define for $ \bt \in \std_{d, n } $ the element $ E_{\bt } \in  \YYfractionfield $ via
\begin{equation}\label{primitiveIdempotent}
  E_{\bt }  = \prod_{a=1}^{2n} \prod_{ \substack{c \in  {\mathcal C}(a) \\ c    \neq \C_{\bt}(a) }}
   \frac{ L_a - c}{ \C_{\bt}(a) -c}. 
\end{equation}

Then, by the results of section 3 in \cite{Mat-So}, the set $ \{ E_{\bt }| \bt \in \std_{d,n} \}$ is a complete
set of
primitive idempotents for $  \YYfractionfield $ and, furthermore, they are eigenvectors for
the action of the $L_a $'s on $  \YYfractionfield $ with associated eigenvalues $ \C_{\bt}(a) $, that
is 
\begin{equation}\label{eigLa}
L_a  E_{\bt } =  E_{\bt } L_a  = \C_{\bt}(a) E_{\bt } \mbox{ for } i= 1, 2,\ldots, 2n . 
\end{equation}  
Note that the formula 
in \eqref{primitiveIdempotent} defining $ E_{\bt} $ also makes sense
for $ \bt \in \tab_{d,n} \setminus \std_{d,n}$, 
but it can be shown that for such $ \bt $ either $ E_{\bt} = 0 $ or $ E_{\bt} = E_{\bs}$ for some $ \bs \in \std_{d,n}$.

\medskip 
In the notation of section \ref{idempotents in YY}, for $ \YYfractionfield $ we have
$ I = {\mathbb Z } $ and $ J = {\mathbb Z }/d  {\mathbb Z } $. 
For $ \bt \in \std_n$, we have 
$  \C_{\bs}(a) = \q^{i_a} $ for $ 1 \le a \le n $ and
$  \C_{\bs}(a) = \xi^{j_a} $ for $ a+1  \le a \le 2n $ where
$ i_a \in I $ and $ j_a \in J$, and so 
\begin{equation}\label{with this notation}
\begin{aligned}
   \ii_{\bt} &= (i_1, i_2, \ldots, i_n) \in \II \, \, \, \,  
 \jj_{\bt} = (  j_1           ,  j_2, \ldots, j_n) \in \JJn \\
  \kk_{\bt} & =\bigl( (i_1, j_a), ( i_2, j_2  ), \ldots, ( i_n, j_n    ) \bigr) \in \KKn. 
\end{aligned} 
\end{equation}
With this notation \eqref{primitiveIdempotent} is the identity 
\begin{equation}
E_{\bt} = e(\kk_{\bt}). 
\end{equation}

\medskip
There is an obvious factorization $ E_{\bt }= E_{\bt }^{(1)} E_{\bt }^{(2)}$ where
\begin{equation}
  E_{\bt }^{(1)}  = \prod_{a=1}^{n} \prod_{ \substack{c \in  {\mathcal C}(a) \\ c    \neq \C_{\bt}(a) }}
  \frac{ L_a - c}{ \C_{\bt}(a) -c} \, \, \mbox{  and  } \, \,
 E_{\bt }^{(2)}  = \prod_{a=n+1}^{2n} \prod_{ \substack{c \in  {\mathcal C}(a) \\ c    \neq \C_{\bt}(a) }}
  \frac{ L_a - c}{ \C_{\bt}(a) -c}.
\end{equation}

For $ \bs,\bt \in \std_{d,n} $ we
write $ \bs \underset{ \le n}{\sim} \bt $ if
$ \C_{\bs}(a) =\C_{\bt}(a)  $
for all $ 1 \le a \le n $, and $ \bs \underset{>n}{\sim} \bt $ if $ \C_{\bs}(a) =\C_{\bt}(a)  $
for all $ n+1 \le a \le 2n $.
Then $ \underset{ \le n}{\sim}$ and $\underset{ >n}{\sim} $ are equivalence relations
on $ \std_{d,n} $ and 
for $ \bt \in  \std_{d,n} $ we let $ [\bt]_{\le n} $ and $ [\bt]_{>n} $ denote
the equivalence classes represented by $ \bt$, with respect to $\underset{ \le n}{\sim} $ and 
$\underset{ >n}{\sim} $, respectively.
Let $  
\{ [\bt_1]_{\le n},  [\bt_2]_{\le n}  , \ldots,  [\bt_{N_1}]_{\le n}  \} = 
$ be the set of 
equivalence classes for $\underset{\le n}{\sim}$ and 
let $ \{ [\bs_1]_{>n},  [\bs_2]_{>n}  , \ldots,  [\bs_{N_2}]_{>n}  \} $ be the set of 
the equivalence classes for $\underset{>n}{\sim}$.

\begin{lemma}\label{E1andE2}
With the above notation we have 
  \begin{description}
  \item[a)]  $ \displaystyle E_{\bt }^{(1)}  = \sum_{ \bs  \in [\bt]_{\le n}    } E_{\bs}$. In particular
    $ X_a E_{\bt }^{(1)} =  E_{\bt }^{(1)} X_a  = \C_{\bt}(a) E_{\bt }^{(1)} $ for all $ 1\le a \le n $. 
  \item[b)]  $E^{(1)}_{\bt} = e(\ii_{\bt})$.
      \item[c)]  $ \displaystyle  E_{\bt }^{(2)}  = \sum_{ \bs    \in [\bt]_{>n}   } E_{\bs}. $
    In particular
    $ t_a E_{\bt }^{(1)} =  E_{\bt }^{(1)} t_a  = \pos_{\bt}(a) E_{\bt }^{(1)} $ for all $n+1 \le a \le 2n. $
  \item[d)] $E^{(2)}_{\bt} = e(\jj_{\bt}).$
  \item[e)] $ E_{\bt }^{(1)} E_{\bs }^{(2)}= \left\{ \begin{array}{ll} E_{\bu } & \mbox{if } \bu
    \in [\bt]_{ \le n} \cap [\bs]_{>n}\\
    0 &  \mbox{if } [\bt]_{ \le n} \cap [\bs]_{>n} = \emptyset. \end{array} \right. $
      \end{description}
\end{lemma}
\begin{dem}
To prove $   { \bf a)} $ we first observe that
  \begin{equation}
    E_{\bt}^{(1)}   E_{\bs}  = \left\{ \begin{array}{ll} E_{\bs} & \mbox{ if } \bs \underset{\le n}{\sim} \bt \\
      0 & \mbox{ otherwise } \end{array} \right.
\end{equation}    
  since, if $ \bs \underset{\le n}{\not\sim} \bt $ the product $  E_{\bt}^{(1)}   E_{\bs} $ contains
  the factor $ (L_a - \C_{\bs}(a) ) E_{\bs}=0$. From this we deduce that
  \begin{equation}
E_{\bt}^{(1)} = E_{\bt}^{(1)}1 = E_{\bt}^{(1)} ( \sum_{\bs \in \std_{d, n} } E_{\bs} ) =  
\sum_{ \bs  \underset{\le n}{\sim} \bt    } E_{\bs}
\end{equation}
  as claimed. The proof of $   { \bf b)} $ follows from $   { \bf a)} $ and the definitions.
  The proofs of $   { \bf c)} $ and $   { \bf d)} $ are similar and finally the proof of $   { \bf e)} $
  follows from $   { \bf a)} $.
\end{dem}

\medskip
It is useful to comment on the equivalence classes for $\underset{\le n}{\sim}$ and $\underset{>n}{\sim}$. The equivalence class $[ \bt ]_{>n}$ is given by permutation of the numbers
within the components of $ \bt$. More precisely, 
let $ \blambda =
( \lambda^{(1)} ,  \lambda^{(2)} , \ldots,  \lambda^{(d)} ) \in \Par_{d,n} $. 
Recall from section \ref{Partitions and tableaux} that we set
$ ||\blambda || = (| \lambda^{(1)}| , | \lambda^{(2)}| , \ldots, | \lambda^{(d)}| ) 
\in \Comp_{n}^d $. Suppose that $ \bt =
(| \bT^{(1)}| , | \bT^{(2)}| , \ldots, | \bT^{(d)}| ) 
\in \std_{d,n}(\blambda) $. Let 
$ || \bt ||  $ be the row standard $ ||\blambda ||$-tableau obtained from
$ \bt $ by placing, in increasing order, the numbers from $  \bT^{(1)} $ in the first row of $||\blambda ||$,
the numbers from $  \bT^{(2)} $ in the second row of $||\blambda ||$, and so on.
For example, for $ \bt $ as in \eqref{1.4} we have
\begin{equation}\label{isconstant}
  || \bt || = \raisebox{-.45\height}{\includegraphics[scale=0.7]{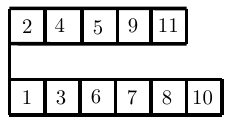}}\, \, 
\raisebox{-14\height}{.}
\end{equation}
Then the function $ \bt \mapsto || \bt ||$ defined this way is constant on
equivalence classes for $\underset{>n}{\sim}$
and it 
induces 
a bijection between  
$ \{ [\bs_1]_{>n},  [\bs_2]_{>n}  , \ldots,  [\bs_{N_2}]_{>n}  \} $ 
and $\Rstd_n^{\Comp, d} $. Note that 
\begin{equation}
\jj_{\bt} = (j_1, j_2, \ldots, j_n ) 
\end{equation}
where $ j_a $ is the row number for $ a $ in $ || \bt||$.

\medskip
The equivalence classes for $\underset{\le n}{\sim}$ are more difficult to describe.
For $ \sigma \in \Si_{d} $ and 
$ \bt= (\T^{(1)}, \T^{(2)}, \ldots, \T^{(d)}) 
\in \std_{d, n}$ we define 
\begin{equation}\label{419}
\sigma \cdot  \bt = (\T^{( \sigma^{-1}(1))}, \T^{( \sigma^{-1}(2))}, \ldots, \T^{(\sigma^{-1}(d))}) .
\end{equation}
Then clearly $ \sigma \cdot  \bt \in [ \bt]_{\le n} $. In general, however, there are more
elements in $ [ \bt]_{\le n} $, but we do not need a full description of 
$ [ \bt]_{\le n} $.

\medskip
Suppose that $ \T \in \Rstd_n^{\Comp, d} $ and
let $ A_j $ be the set of numbers that appear in the $j$'th row of $ \T$.
Then we define $ A_{\T} \in \SetPar_n $ as the set partition obtained from
$ \{ A_1, A_2, \ldots, A_d \} $ be removing the possible empty $ A_i$'s.
For example, if
\begin{equation}\label{421}
\T =   \raisebox{-.45\height}{\includegraphics[scale=0.7]{dib3.pdf}}
\end{equation}
we have $ A_{\T} = \{ \{ 2,3,4,6,7\} , \{ 8,9\}, \{1,5\}\} $. 
Clearly 
any set partition $ A \in \SetPar_n $ is of the form $ A= A_{\T } $ for some $ \T \in \Rstd_n^{\Comp, d} $
and then also of the form $ A_{ || \bt || } $ for some $ \bt \in \std_{d,n} $.
Combining, we get the following Lemma.
\begin{lemma}\label{lemma13}
Let $ A \in \SetPar_n $ and let $ \e( A) $ be as in Theorem \ref{teorem 6A}. 
Then 
\begin{equation}
  \e( A)   =
\sum_{  \substack{ \bt \in \std_{d,n} \\ A_{ || \bt || } = A}  }  E_{  \bt  }
\end{equation}
\end{lemma}
\begin{dem}
This follows from Lemma \ref{E1andE2} $   { \bf c)} $ and \eqref{overline e(A)first}
  together with the description of the classes for
$\underset{>n}{\sim}$ in the paragraphs before \eqref{419}.
\end{dem}

\medskip
We need an ordered set partition variation of $ A_{ \T} $. 
For $ \T \in \Rstd_n^{\Comp, d} $ we let once again
$ A_j $ be the set of numbers that appear in the $j$'th row of $ \T$ and 
define $ A_{\T}^{\ord} \in \OrdSetPar_n $ as the ordered set partition obtained from
$ ( A_1, A_2, \ldots, A_d ) $ by removing the possible empty $ A_i$'s.
For example, for $ \T $ as in \eqref{421} we have
$ A_{\T}^{\ord} = ( \{ 2,3,4,6,7\} , \{ 8,9\}, \{1,5\})  $. 
For $ A^{\ord} \in \OrdSetPar_n $ we then define 
\begin{equation}\label{firstwethendefine}
  \e( A^{\ord})   =
\sum_{  \substack{ \bt \in \std_{d,n} \\ A_{ || \bt || }^{\ord} = A^{\ord}}  }  E_{  \bt  }. 
\end{equation}
Note that, for $ A \in \SetPar_n $ we have that 
\begin{equation}
\e(A) = \sum_{  \substack{ B^{\ord} \in \OrdSetPar_n  \\ B = A }}    \e( B^{\ord}) . 
\end{equation}


\medskip
Recall that $ e_a = \frac{1}{d}\sum_{r=0}^{d-1} t_{a}^{r} t_{a+1}^{-r} $.
Using this we get from $   { \bf c)} $ of Lemma \ref{E1andE2} that 
\begin{equation}\label{613}
  e_a E_{\bt }^{(2)} = \left\{ \begin{array}{ll}  E_{\bt }^{(2)} & \mbox{ if } \pos_{\bt}(a) = \pos_{\bt}(a+1) \\
0 & \mbox{ if } \pos_{\bt}(a) \neq \pos_{\bt}(a+1). \end{array}
\right.    
\end{equation}
Via \eqref{613} we then obtain commutation formulas similar to the ones of Lemma
\ref{veryfirstusefullemma} and Lemma \ref{usefullemma}, as follows.

\begin{lemma}\label{veryfirstusefullemmaA}
  Suppose that $ \bt \in \std_{d,n} $ and that $a  $ and $ a+1 $ appear in the
  same component of $ \bt $. Then the following holds in 
  $ \YYfractionfield $. 
  \begin{description}
  \item[a)]     $  g_a^2 E_{\bt} =  (\q + (\q-1) g_a) E_{\bt}. $
\item[b)]   $     g_a X_{a+1} E_{\bt} = X_{a}( g_a  + (\q-1) )                      E_{\bt}.$
  \item[c)]  $     g_a X_{a} E_{\bt} = X_{a+1}( g_a  + (1-\q) )                     E_{\bt}. $
  \end{description}
\end{lemma}
\begin{dem}
  From \eqref{613} we have
\begin{equation} e_a E_{\bt} = e_a E_{\bt}^{(1)} E_{\bt}^{(2)} =
  e_a E_{\bt}^{(2)} E_{\bt}^{(1)} = E_{\bt}^{(2)} E_{\bt}^{(1)} = E_{\bt}
  \end{equation}
and so we can argue exactly as we did in Lemma \ref{veryfirstusefullemma}. 
 \end{dem}

\begin{lemma}\label{usefullemmaA}
  Suppose that $ \bt \in \std_{d,n} $ and that $a  $ and $ a+1 $ appear in 
  different components of $ \bt $. Then the following formulas hold in 
  $ \YYfractionfield $. 
  \begin{description}
  \item[a)]     $  g_a^2 E_{\bt}= \q E_{\bt} .$
  \item[b)]  $     g_a X_{a+1} E_{\bt} = X_a g_a E_{\bt}. $
\item[c)]  $ g_a X_a E_{\bt} =     X_{a+1} g_a E_{\bt} .  $
  \end{description}
\end{lemma}
\begin{dem}
  From \eqref{613} we have this time
\begin{equation} e_a E_{\bt} = e_a E_{\bt}^{(1)} E_{\bt}^{(2)} =
  e_a E_{\bt}^{(2)} E_{\bt}^{(1)} = 0
  \end{equation}
and so the same arguments that were applied for Lemma \ref{usefullemma} will show the Lemma. 
\end{dem}

\medskip
Using the previous two Lemmas we  arrive at the following commutation relations, involving
the $ E_{\bt}$'s.

\begin{theorem}\label{usefullemmaA}
  Let $ \bt \in \std_{d,n} $ and $ 1 \le a \le n-1$. Define $\bs= \sigma_a \bt$. 
  Then the following commutation formulas hold in $ \YYfractionfield $. 
  \begin{description}
  \item[a)] Suppose that $a  $ and $ a+1 $ appear in 
    the same component of $ \bt $. Then we have 
    \begin{equation}\label{619}
      \left(g_a +(1-\q)\Bigl(1- \C_{\bt}(a)/\C_{\bt}(a+1)\Bigr)^{-1} \right) E_{\bt} =
      E_{\bs} \left(g_a +(1-\q) \Bigl(1- \C_{\bs}(a)/\C_{\bs}(a+1) \Bigr)^{-1} \right) .
    \end{equation}
  \item[b)] Suppose that $a  $ and $ a+1 $ appear in 
    different components of $ \bt $. Then we have 
    \begin{equation}\label{620}
     g_a E_{\bt}  =  E_{\bs} g_a.
    \end{equation} 
  \end{description}
\end{theorem}
\begin{dem}
We first prove $   { \bf a)} $. 
If $ a $ and $ a+1 $ are in the same component and in the
same row or column of $ \bt $, then the left hand side of
\eqref{619} is zero, but also the right hand side of \eqref{619} is zero, since $ E_{\bs} =0 $. 
Hence we may assume that $ \bs \in \std_{d,n}$ and 
proceed to show that   
    \begin{equation}\label{inorder}
      \left(g_a +\frac{(\q-1)\C_{\bs}(a)}{\C_{\bt}(a) - \C_{\bs}(a)} \right) E_{\bt} =
      E_{\bs} \left(g_a - \frac{(\q-1)\C_{\bt}(a)}{\C_{\bt}(a) - \C_{\bs}(a)} \right) 
    \end{equation}
from which $   { \bf a)} $ follows. 
Set $ E^{(1)} = E_{\bs }^{(1)} + E_{\bt}^{(1)} $. Then, by
Lemma \ref{E1andE2} we have
for $ 1 \le a \le n $ that
\begin{equation}
  E_{\bt}^{(1)} = \frac{X_a -\C_{\bs}(a)}{\C_{\bt}(a) -\C_{\bs}(a)}E^{(1)} \, \, \, \,  \mbox{and} \, \, \, 
  E_{\bs}^{(1)} = \frac{X_a -\C_{\bt}(a)}{\C_{\bs}(a) -\C_{\bt}(a)}E^{(1)} .
\end{equation}  
Furthermore, noting that $ E^{(1)} $ is a symmetric polynomial in $ X_a $ and $X_{a+1}$, we get from
Corollary \ref{firstcor} that $ g_a $ commutes with $ E^{(1)} $. Hence
\begin{equation}\label{622}
\begin{aligned}      
   g_a E_{\bt}^{(1)} E_{\bt}^{(2)}  &=  \frac{g_a(X_a -\C_{\bs}(a))}{\C_{\bt}(a) -\C_{\bs}(a)}E^{(1)} E_{\bt}^{(2)} =  \frac{X_{a+1}(g_a +1-\q) -g_a \C_{\bs}(a)}{\C_{\bt}(a) -\C_{\bs}(a)}E^{(1)} E_{\bt}^{(2)}  \\
   &=  \frac{X_{a+1}(E^{(1)}_{\bs} +E^{(1)}_{\bt} )
     (g_a +1-\q) - \C_{\bs}(a) (E^{(1)}_{\bs} +E^{(1)}_{\bt} )g_a}{\C_{\bt}(a) -\C_{\bs}(a)} E_{\bt}^{(2)}
\\
   &=  \frac{(\C_{\bs}(a+1) E^{(1)}_{\bs} +\C_{\bt}(a+1)E^{(1)}_{\bt} )
  (g_a +1-\q) - \C_{\bs}(a) (E^{(1)}_{\bs} +E^{(1)}_{\bt} )g_a}{\C_{\bt}(a) -\C_{\bs}(a)} E_{\bt}^{(2)}
\\
   &=  \frac{(\C_{\bt}(a) E^{(1)}_{\bs} +\C_{\bs}(a)E^{(1)}_{\bt} )
  (g_a +1-\q) -\C_{\bs}(a) (E^{(1)}_{\bs} +E^{(1)}_{\bt} )g_a }{\C_{\bt}(a) -\C_{\bs}(a)} E_{\bt}^{(2)}
\\
&=  \frac{(\C_{\bt}(a)  -\C_{\bs}(a) ) E^{(1)}_{\bs}g_a+(1-\q) (\C_{\bt}(a) E_{\bs}^{(1)}  +\C_{\bs}(a) E_{\bt}^{(1) })}{\C_{\bt}(a) -\C_{\bs}(a)} E_{\bt}^{(2)}
\\
   &= E^{(1)}_{\bs}g_aE_{\bt}^{(2)}  -(\q-1) \frac{  \C_{\bt}(a) E_{\bs}^{(1)}  +\C_{\bs}(a) E_{\bt}^{(1) }}{\C_{\bt}(a) -\C_{\bs}(a)} E_{\bt}^{(2)} 
 \end{aligned} 
\end{equation}
where we used $   { \bf a)} $ of 
Lemma \ref{veryfirstusefullemmaA} for the second equality: note that
$ E^{(1)} $ and $ E^{(2)}_{\bt} $ commute. Comparing the first and the last expression in \eqref{622}
we get
\begin{equation}\label{623}
  \left(g_a +(\q-1) \frac{\C_{\bs}(a)}{\C_{\bt}(a) - \C_{\bs}(a)}\right) E_{\bt}^{(1)} E_{\bt}^{(2)}
= E_{\bs}^{(1)}\left(g_a -(\q-1) \frac{\C_{\bt}(a)}{\C_{\bt}(a) - \C_{\bs}(a)}\right)  E_{\bt}^{(2)} .
\end{equation}
Finally, 
since $ a $ and $ a+1 $ appear in the same component of $ \bt$ we have 
$ E_{\bt}^{(2)}  = E^{(2)}_{\bs}  $ and $ g_a  E^{(2)}_{\bs} = E^{(2)}_{\bs} \! g_a $ and
so the last expression of \eqref{623} is equal to 
$ E_{\bs}^{(1)} E_{\bs}^{(2)} \left(g_a -(\q-1) \frac{\C_{\bs}(a)}{\C_{\bt}(a) - \C_{\bs}(a)}\right)   $, which
shows \eqref{inorder} and therefore also 
$   { \bf a)} $. 

\medskip
To show $   { \bf b)} $, we first observe that $   { \bf b)} $ and $   { \bf c)} $ of
Lemma \ref{usefullemmaA} imply that
$ g_a  E_{\bt}^{(1)} E_{\bt}^{(2)}  =   E_{\bs}^{(1)} g_a E_{\bt}^{(2)}   $ whereas 
relation \eqref{eq three} implies that $ g_a E_{\bt}^{(2)} =  E_{\bs}^{(2)} g_a $. Combining, we 
have also proved $   { \bf b)} $. 
\end{dem}

\medskip
\medskip
Recall the cellular basis $ \{ c_{\bs \bt} | \bs, \bt \in \std(\blambda),  \blambda \in \Par_{d,n} \}$ 
for $ \YY$ that was constructed in \cite{ERH}. Following the general theory laid out in
\cite{Mat-So} we now define {\it the seminormal basis} $ \{ f_{\bs \bt} | \bs, \bt \in \std(\blambda),  \blambda \in \Par_{d,n} \}$ 
for $ \YYfractionfield $ as follows
\begin{equation}
  f_{\bs \bt} = E_{\bs} c_{\bs \bt} E_{\bt} \mbox{  where  } \bs, \bt \in \std(\blambda),
\blambda \in \Par_{d,n}.
 \end{equation}
The following Theorem is an analogue of {\it Young's seminormal form} in
the setting.
\begin{theorem}\label{YSF}
Suppose that $ \bt $ and $ \bu $ are standard multitableaux of the same shape
and that $ \bs = \sigma_a \bt $ for some $ 1 \le a < n$. 
\begin{description}
\item[a)] 
If $ \bs $ is standard the following formulas hold in $ \YYfractionfield$
\begin{equation}\label{is bs standard}
  g_a f_{\bt \bu} = \left\{ \begin{array}{ll} \displaystyle \frac{\q-1}{1 - \C_{\bt}(a) \C_{\bt}(a+1)^{-1} } f_{\bt \bu}
    + f_{\bs \bu} & \mbox{if } \bs \lhd \bt,  \pos_{\bt}(a) =\pos_{\bt}(a+1) \\[10pt]  \displaystyle
 \frac{\q-1}{1 - \C_{\bt}(a) \C_{\bt}(a+1)^{-1} } f_{\bt \bu}
 +  \frac{ (\q \C_{\bt}(a)- \C_{\bt}(a+1)) ( \C_{\bt}(a)- \q \C_{\bt}(a+1))  }{( \C_{\bt}(a)-  \C_{\bt}(a+1))^2} f_{\bs \bu}
 & \mbox{if } \bs \rhd \bt,  \pos_{\bt}(a) =\pos_{\bt}(a+1) \\[10pt] 
f_{\bs \bu}
& \mbox{if }   \pos_{\bt}(a) <\pos_{\bt}(a+1)
\\[10pt] 
\q f_{\bs \bu}
 & \mbox{if }   \pos_{\bt}(a) >\pos_{\bt}(a+1).
\end{array}
\right.
\end{equation}
\item[b)] 
If $ \bs $ is non-standard the following formulas hold in $ \YYfractionfield$
\begin{equation}\label{is bs non standard}
  g_a f_{\bt \bu} = \left\{ \begin{array}{ll} \q f_{\bt \bu}  & \mbox{if } a \mbox{ and } a+1 \mbox{ are in the same row of } \bt \\[10pt]
    -f_{\bt \bu}& \mbox{if }
a \mbox{ and } a+1 \mbox{ are in the same column of } \bt.
\end{array}
\right.
\end{equation}
\end{description}
\end{theorem}
\begin{dem}
The proof is an adaption of the proof of Theorem 3.36 in \cite{mathasbook}, which corresponds 
to the Hecke algebra case.
If $ \pos_{\bt}(a) = \pos_{\bt}(a+1) $ we get 
by $   { \bf a)} $ of Theorem \ref{usefullemmaA} that
for any $ \bnu \in \std_{d,n} $ we have
\begin{equation}
\begin{aligned}
  E_{\bnu} g_a E_{\bt} &=
  E_{\bnu}
  \left(g_a + \frac{(\q-1) \C_{\bs}(a)}{\C_{\bt}(a) - \C_{\bs}(a)}
- \frac{(\q-1) \C_{\bs}(a)}{\C_{\bt}(a) - \C_{\bs}(a)}
\right) E_{\bt}
 \\ & =
  E_{\bnu} E_{\bs}
  \left(g_a - \frac{(\q-1) \C_{\bs}(a)}{\C_{\bt}(a) - \C_{\bs}(a)}
\right)   - \frac{(\q-1) \C_{\bs}(a)}{\C_{\bt}(a) - \C_{\bs}(a)} E_{\bnu} E_{\bt}
\end{aligned}
\end{equation}  
which is non-zero only if $ \bnu = \bs $ or $ \bnu = \bt $ and hence
$    g_a f_{\bt \bu} = c_1 f_{\bt \bu} + c_2 f_{\bs \bu} $ for some coefficients $c_1, c_2 \in \Kdomain(\q)$.
Similarly, if $ \pos_{\bt}(a) \neq \pos_{\bt}(a+1) $ we get 
$   { \bf b)} $ of Theorem \ref{usefullemmaA} that
\begin{equation}
  E_{\bnu} g_a E_{\bt} = E_{\bnu}  E_{\bs} g_a 
\end{equation}
which is non-zero only if $ \bnu = \bs $ and so
$    g_a f_{\bt \bu} =  c f_{\bs \bu} $ for some coefficient $c \in \Kdomain(\q)$. 
Given this, we proceed as in the proof of Theorem 3.36 of \cite{mathasbook} to arrive 
at \eqref{is bs standard} and
\eqref{is bs non standard}.
\end{dem}

\medskip
Let 
$ \YYfractionfielda = 
\spa_{\Kdomain(\q)} \{ f_{\bs \bt} |
\bs, \bt \in \std_{d,n}, shape(\bs) = shape(\bt),  \pos_{\bs}(a) = \pos_{\bs}(a+1) \}  $.
Then by 
\eqref{eigLa} the eigenvalues of $ 1-X_a X_{a+1}^{-1} $ acting on $\YYfractionfielda $
are non-zero and so $ 1-X_a X_{a+1}^{-1} $ is invertible on $\YYfractionfielda  $. 
Using this we get the following Corollaries of Theorem \ref{YSF}.

\begin{corollary}\label{YSFcorA}
Suppose that $ \bt $ and $ \bu $ are standard multitableaux of the same shape
and that $ \bs = \sigma_a \bt $ for some $ 1 \le a < n$. 
\begin{description}
\item[a)] 
If $ \bs $ is standard the following formulas hold in $ \YYfractionfield$
\begin{equation}\label{is bs standard}
  (g_a+1) f_{\bt \bu} = \left\{ \begin{array}{ll}
    \displaystyle \frac{X_a- \q X_{a+1}}{X_a -  X_{a+1} } f_{\bt \bu}
    + f_{\bs \bu} & \mbox{if } \bs \lhd \bt,  \pos_{\bt}(a) =\pos_{\bt}(a+1) \\[10pt]  \displaystyle
 \frac{X_a - \q X_{a+1}}{X_a-  X_{a+1 }} f_{\bt \bu}
 +  \frac{ (\q X_a- X_{a+1}) ( X_a- \q X_{a+1})  }{( X_a-  X_{a+1})^2} f_{\bs \bu}
 & \mbox{if } \bs \rhd \bt,  \pos_{\bt}(a) =\pos_{\bt}(a+1) \\[10pt] 
f_{\bt \bu} +f_{\bs \bu}
& \mbox{if }   \pos_{\bt}(a) <\pos_{\bt}(a+1)
\\[10pt] 
f_{\bt \bu} +\q f_{\bs \bu}
 & \mbox{if }   \pos_{\bt}(a) >\pos_{\bt}(a+1).
\end{array}
\right.
\end{equation}
\item[b)] 
If $ \bs $ is non-standard the following formulas hold in $ \YYfractionfield$
\begin{equation}\label{is bs non standard}
  g_a f_{\bt \bu} = \left\{ \begin{array}{ll} (\q +1)f_{\bt \bu}  & \mbox{if } a \mbox{ and } a+1 \mbox{ are in the same row of } \bt \\[10pt]
    0& \mbox{if }
a \mbox{ and } a+1 \mbox{ are in the same column of } \bt.
\end{array}
\right.
\end{equation}
\end{description}
\end{corollary}

\begin{corollary}\label{YSFcor}
Let $ \bt $ and $ \bu $ be standard multitableaux of the same shape
such that $ \bs = \sigma_a \bt $ for some $ 1 \le a < n$. 
\begin{description}[leftmargin=!]
\item[a)] 
Suppose that $ \bs $ is standard 
and that 
$ a $ and $ a+1 $ appear in the same component of $ \bt$. 
Then 
\begin{equation}\label{is bs standard}
  \left(g_a+ \frac{1-\q}{1 - X_a X_{a+1}^{-1} } \right)
  f_{\bt \bu}  = \left\{ \begin{array}{ll} \displaystyle 
     f_{\bs \bu} & \mbox{if } \bs \lhd \bt \\[10pt]  \displaystyle
   \frac{ (\q X_a- X_{a+1}) ( X_a- \q X_{a+1} ) }{( X_a-  X_{a+1})^2} f_{\bs \bu}
 & \mbox{if } \bs \rhd \bt. 
\end{array}
  \right.
\end{equation}

\item[b)] 
  Suppose that $ \bs $ is non-standard, that is $ a $ and $ a+1 $ appear in the same row or column
  of $ \bt$ (and hence also in the same 
  component of 
  $ \bt $). Then 
\begin{equation}\label{is bs non standard}
  \left(g_a+ \frac{1-\q}{1 - X_a X_{a+1}^{-1} } \right)
  f_{\bt \bu} = 0.
\end{equation}
\item[c)]
Suppose that $ a $ and $ a+1 $ appear in different components of $ \bs $. Then 
\begin{equation}  
g_a f_{\bt \bu} = \left\{ \begin{array}{ll} \displaystyle 
      f_{\bs \bu} & \mbox{if } \pos_{\bs}(a) < \pos_{\bs}(a+1) \\[10pt]  \displaystyle
\q f_{\bs \bu}   
 & \mbox{if } \pos_{\bs}(a) > \pos_{\bs}(a+1).
\end{array}
\right.
\end{equation}  
\end{description}
\end{corollary}

Our main interests, however, do not lie in the representation theory of $ \YYfractionfield $, but rather in
the representation theory of 
$ \YYresiduefield$, for which the content function given in \eqref{contentfunction}, 
does not satisfy the separation condition of \cite{Mat-So}. Therefore, 
idempotents $ E_{\bt}  $ satisfying \eqref{eigLa} do not exist for $\YYresiduefield $.

\medskip
We can however still apply
the results of section 4 in \cite{Mat-So}.
Let therefore $ \sim_{q, t}$ be the equivalence relation
on $ \std_{d, n} $ given by 
\begin{equation}\label{equivalence qt}
\bs \sim_{q,t} \bt \, \mbox{ if and only if } c_{\bs}(a) = c_{\bt}(a) \, \, \mbox{for all } a=1,2,\ldots, 2n.
\end{equation}
Let $ \TT $ be an equivalence class for $ \sim_{q,t} $, and let $ \bt \in \TT$.
Then there are residue sequences $ \ii =(i_1, i_2, \ldots, i_n) \in \II $ and $ \jj =
(j_1, j_2, \ldots, j_n ) \in \JJn $ such that
\begin{equation}
  (c_{\bt}(1), c_{\bt}(2), \ldots, c_{\bt}(n)) = (q^{i_1}, q^{i_2}, \ldots, q^{i_n} ) \, \, \mbox{and} \, \,
(c_{\bt}(n+1), c_{\bt}(n+2), \ldots, c_{\bt}(2n)) = (\xi^{j_1}, \xi^{j_2}, \ldots, \xi^{j_n} ).
  \end{equation}  
Note that $ \ii $ and $ \jj $ are independent of the choice of $ \bt \in \TT $ and so we
define $ \res(\TT) \in \KKn $ via
\begin{equation}
  \res(\TT) = ((i_1, j_1), (i_2, j_2), \ldots, (i_n, j_n))
\end{equation}
The function $\res: \TT \mapsto \KKn$ is injective. 

\medskip
For an equivalence class $ \TT $ for $ \sim_{q,t} $ we next introduce the element $ \hat{E}_{ \TT}  $ via 
\begin{equation}\label{classidempotent}
\hat{E}_{ \TT} = \sum_{ \bt\in \TT} E_{\bt} .
\end{equation}
Then by definition $ \hat{E}_{ \TT} \in \YYfractionfield $, but by the general results in section 4 of \cite{Mat-So},
in fact $ \hat{E}_{ \TT} \in \YYOO $ and so by reduction modulo
$ (\q -q ) \OO$ we obtain 
an element $ E_{ \TT} $ of $ \YYresiduefield$.

\medskip
Our main reason for considering $ E_{ \TT} $ is the following Theorem. In the case of the 
usual Hecke algebra, that is the case $ d= 1 $, it was
proved by Hu-Mathas in \cite{hu-mathas}. 
\begin{theorem}\label{teo17}
  Suppose that $ \TT $ is an equivalence class for $ \sim_{q,t} $ and that $ \res(\TT) = \kk$.
  Then the following identity holds in $\YYresiduefield$
  \begin{equation}
E_{\TT} = e(\kk) 
  \end{equation}    
where $ e(\kk) \in \YYresiduefield $ is as in \eqref{eigenspace}. 
\end{theorem}
\begin{dem}
Suppose that $ \kk = ((i_1, j_1), (i_2, j_2), \ldots, (i_n, j_n)) $ and 
that $ \TT =\{\bt_1, \bt_2, \ldots, \bt_N \} $.
For $ m \in \Z $ we define $ \logq(\q^{ \,m }) =m$ 
  and set 
  $l_x = \logq( \hat{c}_{\bt_x}(a)) $ for $ x =1, 2, \ldots, N$. Then $  q^{ l_x } = q^a$ for all $ x$
and for $ 1 \le a \le n $ we have 
\begin{equation}
  \begin{aligned}      
    (X_a - q^{ i_a } )^{N}E_{\TT}  &= (X_a-\q^{l_1})  (X_a-\q^{l_2})  \cdots (X_a-\q^{l_{N}})  \hat{E}_{\TT} \, \,  \mbox{\rm  mod } (\q -q )\OO \\
    &= \sum_{i=1}^N(X_a-\q^{l_1})  (X_a-\q^{l_2})  \cdots (X_a-\q^{l_{N}})  \hat{E}_{\bt_i} \, \,  \mbox{\rm  mod } (\q -q )\OO  =0.
  \end{aligned}      
  \end{equation}
For $n+1  \le a \le 2n $ we get, arguing the same way, that
$ (X_a - \xi^{ j_a } )E_{\TT} = 0 $ and so we have $ E_{\TT} \in e(\kk) \YYresiduefield$. 
But the $ E_{\TT}$'s are a complete set of idempotents for $ \YYresiduefield$ and so  
$ E_{\TT} = e(\kk) $, as claimed.
\end{dem}

\medskip
We have a couple of Corollaries of the Theorem. 

\begin{corollary}
  For $ \kk \in \KK $ we have $ e(\kk) =0 $ unless $ \kk = \res(\TT) $ for some equivalence class $ \TT $ for
  $ \sim_{q,t} $. 
\end{corollary}

\begin{corollary}[see also \cite{WeidengCuiJinkuiWan}]\label{eigenvalues}
  Let $ \kk =  ((i_1, j_1), (i_2, j_2), \ldots, (i_n, j_n)) \in \KK $.
  Then the minimal polynomial for $ X_a $ acting in $ e(\kk) \YYresiduefield  $ is
  $ (x-q^{i_a})^{N_a} $ for a natural number $ N_a $. Hence,
the minimal polynomial for $ X_a $ acting in $ \YYresiduefield  $ is  
  of the form $ \prod_i (x-q^i)^{N_{a,i}}$ for $ i \in I $ and natural numbers 
  $ N_{a,i} $ and therefore the eigenvalues for $ X_a $ acting in $ \YYresiduefield $ are of the form $  q^i$
  for $ i \in I$. 
\end{corollary}

Suppose now that $ f(L_1, L_2, \ldots, L_{2n}) \in
\Kdomain[\q, \q^{-1}][L_1^{\pm 1}, L_2^{\pm 1} , \ldots, L_{2n}^{\pm 1} ] $. Then
we have $ f(L_1, L_2, \ldots, L_{2n}) \in \YYOO $ and therefore
$  f(L_1, L_2, \ldots, L_{2n}) \hat{E}_{\TT} \in \YYOO$. 
Viewing $ \YYOO$ inside $ \YYfractionfield $ we get
\begin{equation}\label{631}
f(L_1, L_2, \ldots, L_{2n}) \hat{E}_{\TT} = \sum_{\bt \in \TT} f(c_{\bt}(1), c_{\bt}(2), \ldots, c_{\bt}(2n)) {E}_{\bt}.
\end{equation}
In particular, the right hand side of \eqref{631} belongs to $ \YYOO$ and may therefore
also be viewed as an element of $ \YYresiduefield$
by reduction modulo
$ (\q -q ) \OO $.
\medskip

Suppose that $ \kk =  ((i_1, j_1), (i_2, j_2), \ldots, (i_n, j_n)) \in \KKn $. 
Then, by the last Corollary, the only eigenvalue for the action of $X_a $ 
via left multiplication on $ e(\kk) \YY$ 
is $ q^{i_a} $, and similarly the only eigenvalue for the action of $X_{a+1} $ is $ q^{i_{a+1}} $
and therefore, since $ X_a $ and $ X_{a+1} $ commute,  
the only eigenvalue for the action of $ X_a X_{a+1}^{-1} $ on
$ e(\kk) \YY$ 
is $  q^{i_a -i_{a+1}} $. Hence, if $ i_a \neq i_{a+1 }$ the only eigenvalue for the action of $ 1-X_a X_{a+1}^{-1} $
on $ e(\kk) \YY$ is $ 1-  q^{i_a -i_{a+1}} \neq 0 $, and therefore $ 1-X_a X_{a+1}^{-1}$ is invertible on 
$ e(\kk) \YY$. In more detail, on Jordan blocks for $ 1-X_a X_{a+1}^{-1} $ of length $ N $ we
have
\begin{equation}\label{inmoredetail}
  \begin{aligned}
    (1-X_a X_{a+1}^{-1})^{-1} & = (X_a X_{a+1}^{-1} - q^{i_a-i_{a+1}})^{N-1}( 1-q^{i_a-i_{a+1}})^{-N} \\ &+
(X_a X_{a+1}^{-1} - q^{i_a-i_{a+1}})^{N-2}( 1-q^{i_a-i_{a+1}})^{-N+1}
    +\ldots  \\ &+
   (X_a X_{a+1}^{-1} - q^{i_a-i_{a+1}})^{1}( 1-q^{i_a-i_{a+1}})^{-2} \\ &+ ( 1-q^{i_a-i_{a+1}})^{-1} 
  \end{aligned}
\end{equation}
acting on $ e(\kk) \YY$. Note that \eqref{inmoredetail} makes sense since 
$q^{i_a-i_{a+1}}-1 \neq 0 $. Choosing $ N $ to be larger than the length of all the Jordan blocks, we get that 
\eqref{inmoredetail} holds on all $ e(\kk) \YY$. 

\medskip
For $ i \in I $ we define $ \hat{i} \in \Z$ by the conditions $ 0 \le \hat{i} < {\rm char }_q (\Kdomain) $
and $ \hat{i}\, \, {\rm mod} \, \, {\rm char }_q (\Kdomain) = i $. Then combining 
\eqref{631} and \eqref{inmoredetail} we get the following lift of $ (1-X_a X_{a+1}^{-1})^{-1} {E}_{\TT}  $ to
$ \YYOO $. 

\begin{equation}\label{inverselift}
  \begin{aligned}
    lift
    & = \sum_{\bt \in \TT} \left(\C_{\bt}(a)/\C_{\bt}(a+1) -
    \q^{\, \hat{i}_a-\hat{i}_{a+1}}\right)^{N-1}(1- \q^{\, \hat{i}_a-\hat{i}_{a+1}})^{-N} E_{\bt} \\ &+
\sum_{\bt \in \TT} \left(\C_{\bt}(a)/\C_{\bt}(a+1) - \q^{\, \hat{i}_a-\hat{i}_{a+1}}\right)^{N-2}(1- \q^{\, \hat{i}_a-\hat{i}_{a+1}})^{-N+1} E_{\bt}
    +\ldots  \\ &+
    \sum_{\bt \in \TT}   \left(\C_{\bt}(a)/\C_{\bt}(a+1) - \q^{\, \hat{i}_a-\hat{i}_{a+1}}\right)^{1}(1- \q^{\, \hat{i}_a-\hat{i}_{a+1}})^{-2} E_{\bt} \\ &+
    \sum_{\bt \in \TT} (1- \q^{\, \hat{i}_a-\hat{i}_{a+1}})^{-1} E_{\bt}. 
  \end{aligned}
\end{equation}
Using the geometric sum, \eqref{inverselift} becomes 
\begin{equation}\label{inverselift2}
\sum_{\bt \in \TT}  \left(\left(\frac{\C_{\bt}(a)/\C_{\bt}(a+1)- \q^{\, \hat{i}_a-\hat{i}_{a+1}}   }
      {1- \q^{\hat{i}_a -\hat{i}_{a+1}} } \right)^N -1 \right) \bigg/ \bigg(\C_{\bt}(a)/\C_{\bt}(a+1)-1 \bigg) \hat{E}_{\bt}.
\end{equation}
On the other hand, since $ \frac{\C_{\bt}(a)/\C_{\bt}(a+1)- \q^{\, \hat{i}_a-\hat{i}_{a+1}}   }
{1- \q^{\hat{i}_a -\hat{i}_{a+1}} }  $ belongs to $ (q-\q ) \OO$, there exists a sufficiently large $ N $
such that for all $ \bt \in \TT $ 
\begin{equation}\label{inverselift3}
\left(\frac{\C_{\bt}(a)/\C_{\bt}(a+1)- \q^{\, \hat{i}_a-\hat{i}_{a+1}}   }
{1- \q^{\hat{i}_a -\hat{i}_{a+1}} } \right)^N E_{\bt}  = 0 \mbox{ mod }  (q-\q ) \OO .
\end{equation}
Combining \eqref{inverselift}, \eqref{inverselift2} and \eqref{inverselift3}
we then deduce that also 
\begin{equation}\label{636}
  (1-X_a X_{a+1}^{-1})^{-1} \hat{E}_{\TT}  =
\sum_{\bt \in \TT} \big(1 -\C_{\bt}(a)/\C_{\bt}(a+1) \big)^{-1} \hat{E}_{\bt}
\end{equation}
is a lift of $ (1-X_a X_{a+1}^{-1})^{-1} {E}_{\TT}  $ to
$ \YYOO $. Note that, apriori, it is not obvious that the right hand side of 
\eqref{636} belongs to $ \YYOO$.

\medskip
Unfortunately, if $ i_a = i_{a+1} $ the linear map 
$ (1-X_a X_{a+1}^{-1})^{-1} $ is not defined on $ e(\kk) \YYresiduefield $.
To overcome this obstacle, in order to extend $ (1-X_a X_{a+1}^{-1})^{-1} $ to all $ \YYresiduefield $
we introduce 
$ \F_a \in   \YYresiduefield $ as follows
\begin{equation}\label{intertwinedef}
  \F_a =
   g_a + (1-q) \sum_{\substack{ \kk \in \KKn \\  i_a \neq i_{a+1} \\ j_a = j_{a+1} }}
  (1-X_a X_{a+1}^{-1})^{-1} e(\kk) +
  \sum_{\substack{ \ii \in \KKn \\  i_a = i_{a+1} \\ j_a = j_{a+1} }} e(\kk) 
\end{equation}
where $ \kk  =  ((i_1, j_1), (i_2, j_2), \ldots, (i_n, j_n)) $.
These $ \F_a$'s are the {\it intertwining elements} considered in \cite{Ro}, and originally for $ d=1 $ in
\cite{brundan-klesc}.

\medskip
We now use the results developed so far in this section to prove the following Theorem
involving the $\F_a$'s.
\begin{theorem}\label{commutation intertwiner}
Let $ \TT $ be an equivalence class for $ \sim_{q,t} $ and
set
$ \res(\TT) = \kk  =  ((i_1, j_1), (i_2, j_2), \ldots, (i_n, j_n)) \in \KKn  $. 
Suppose that 
$ 1 \le a < n  $.     
Then we have 
\begin{align}
\label{intertwine one} \F_a  e(\kk)  & = e(\sigma_a \kk) \F_a            \\
\label{intertwine two} \F_a  X_b  & = X_b \F_a  \, \, \, \, \, \, \, \, \, \, \, \,     \mbox{if } b  \neq  a, a+1        \\
 \label{intertwine three} \F_a  \F_b  & = \F_b \F_a \, \, \, \, \, \, \, \, \, \, \, \, \,  \mbox{if } | a-b | > 1           \\
 \label{intertwine four} \F_a X_{a+1} e(\kk)    &   =  \left\{\begin{aligned} & X_a \F_a e(\kk) & &  \mbox{if }
(i_a \neq i_{a+1}, \,  j_a = j_{a+1}) \mbox{ or }  j_a \neq j_{a+1} \\
&  (X_a \F_a +q X_{a+1} -X_a)  e(\kk) & & \mbox{if } i_a = i_{a+1}, \, j_a = j_{a+1} \\
 \end{aligned} \right.    \\
\label{intertwine five}  X_{a+1} \F_a e(\kk)    &   =  \left\{\begin{aligned} & X_a \F_a e(\kk) & &  \mbox{if }
(i_a \neq i_{a+1}, \, j_a = j_{a+1}) \mbox{ or } j_a \neq j_{a+1} \\
&  ( \F_a X_a +q X_{a+1} -X_a)  e(\kk) & & \mbox{if } i_a = i_{a+1}, \, j_a = j_{a+1}  \end{aligned} \right.    \\
\label{intertwine six}  \F_a^2 e(\kk)    &   =  \left\{\begin{aligned} &
\frac{(q X_{a} -X_{a+1}) (X_{a} -qX_{a+1})}{(X_a-X_{a+1})^2} e(\kk)  & &  \mbox{if }
i_a \neq i_{a+1}, \, j_a = j_{a+1} \\
&  (q+1) \F_a e(\kk) & & \mbox{if } i_a = i_{a+1}, \,  j_a = j_{a+1} \\
&  q \F_a e(\kk) & & \mbox{if } j_a \neq j_{a+1}
\end{aligned} \right.    \\
\label{intertwine seven} \F_a \F_{a+1} \F_{a}   e(\kk)    &   =  \left\{\begin{aligned}
& (\F_{a+1} \F_{a} \F_{a+1}  +q \F_a -q \F_{a+1} ) e(\kk) & &  \mbox{if } i_a = i_{a+1} = i_{a+2}, \,
j_a = j_{a+1} = j_{a+2}
\\
& (\F_{a+1} \F_{a} \F_{a+1}+Z_a) e(\kk) & &  \mbox{if }  i_a = i_{a+2} \neq i_{a+1}, \,
j_a = j_{a+1} = j_{a+2}       \\
& \F_{a+1} \F_{a} \F_{a+1} e(\kk) & &  \mbox{otherwise}  \\
\end{aligned} \right.    
\end{align}
where $ Z_a = (1-q)^2\frac{ (X_aX_{a+2} -X_{a+1}^2) (X_a X_{a+1} -q X_{a+1} X_{a+2})}
{(X_a-X_{a+1})^2(X_{a+1}-X_{a+2})^2}$.
  \end{theorem}  
\begin{dem}
We use the concrete description of $ e(\kk) $ given in \eqref{classidempotent}
and Theorem \ref{teo17}.

\medskip  
We first show \eqref{intertwine one}. If $i_a \neq i_{a+1} $ and $j_a = j_{a+1} $ 
it follows directly
from $ {\bf a)} $ of Theorem \ref{usefullemmaA} and
\eqref{636}. If $i_a = i_{a+1} $ and $j_a = j_{a+1} $ it is the statement that $ g_a e(\kk) = e(\kk) g_a $
or equivalently $g_a E_{\TT} =      E_{\TT} g_a$.
Under these assumptions we have for any $ \bt \in \TT $ that 
$ a$ and $a+1$ neither appear in the same row of $ \bt$, nor in the same column
  of $ \bt$ and so 
$ \TT $ is the disjoint union of $ \TT^1 $ and $ \TT^2 $ where 
\begin{equation}\label{disjoint}
 \TT^1 = \{ \bt \in \TT \, | \, \sigma_a \bt \lhd \bt  \} \, \, \,  \mbox{and} \, \, \,  \TT^2 = \{ \bt \in \TT \, | \, \sigma_a \bt \rhd \bt  \}.
\end{equation}
Suppose that $ \bt_1 \in \TT^1 $, and therefore $ \bt_2 = \sigma_a \bt_1 \in \TT^2$. Applying $ {\bf a)} $ of Theorem
\ref{usefullemmaA} on each of the pairs $ (\bt_1, \bt_2) $ and $ (\bt_2, \bt_1) $ we get 
\begin{equation}\label{summing}
  \begin{aligned}
    g_a E_{\bt_1} &=  E_{\bt_2}g_a + \frac{1-\q}{\C_{\bt_1}(a) - \C_{\bt_2}(a)}\bigg( \C_{\bt_1}(a)E_{\bt_2} +\C_{\bt_2}(a)E_{\bt_1}\bigg) \\
    g_a E_{\bt_2} &=  E_{\bt_1}g_a + \frac{1-\q}{\C_{\bt_2}(a) - \C_{\bt_1}(a)}\bigg( \C_{\bt_2}(a)E_{\bt_1} +\C_{\bt_1}(a)E_{\bt_2}\bigg).
\end{aligned}    
\end{equation}
Summing the two equations in \eqref{summing} we get $ g_a (E_{\bt_1} + E_{\bt_2} ) =  (E_{\bt_1} + E_{\bt_2} ) g_a $ and via 
\eqref{disjoint}
and \eqref{classidempotent}
we then arrive at $g_a \hat{E}_{\TT} =      \hat{E}_{\TT} g_a$ 
and hence $g_a e(\kk) =      e(\kk) g_a$, as claimed.
Finally, if $ j_a \neq j_{a+1} $ we have that \eqref{intertwine one}
is the statement $g_a E_{\TT} = E_{\SSS} g_a $
which follows immediately from $ {\bf b)} $ of Theorem of \ref{usefullemmaA}.

\medskip
For the proofs of the remaining identities
\eqref{intertwine two}-\eqref{intertwine seven} 
we work in $ \YYfractionfield $ and check that the left hand side of each identity acts
the same way as the right hand side
on the seminormal basis $ \{ f_{\bt \bu}\} $.

\medskip
In view of this, \eqref{intertwine two} is an immediate consequence of 
Corollary \ref{YSFcorA} and Corollary \ref{YSFcor} since $ \F_a f_{\bt \bu} $ is given by either
$ {\bf a)} $, $ {\bf b)} $ of Corollary \ref{YSFcorA} or by 
$ {\bf a)} $, $ {\bf b)} $ or $ {\bf c)} $ of Corollary \ref{YSFcor}.
In each case $\F_a f_{\bt \bu} $, is just a linear combination of
$   f_{\bt \bu} $ and $   f_{\bs \bu} $ where $ \bs = \sigma_a \bt$, 
and hence $ X_b $ acts the same way on $ f_{\bt \bu} $ and $\F_a f_{\bt \bu} $, which proves 
\eqref{intertwine two}. 

\medskip
The other identities 
\eqref{intertwine three}-\eqref{intertwine seven} are 
proved essentially the same way as \eqref{intertwine two}. We explain 
the second identity of \eqref{intertwine seven}, which is the most complicated of
the identities 
\eqref{intertwine three}-\eqref{intertwine seven}, and leave the rest to the reader.

\medskip
By the definition of $ \F_a $ in \eqref{intertwinedef},
the second identity of \eqref{intertwine seven} corresponds to
\begin{equation}\label{correspondstoproving}
 \F_a (g_{a+1} +1) \F_a E_{\bt}   = \Big(\F_{a+1} (g_{a} +1) \F_{a+1} +Z_a \Big) E_{\bt}.
\end{equation}
where $ \bt \in \TT $ satisfies that 
$ a, a+1 $ and $ a+2$ appear in the same component of $ \bt $
and has residues $ i_a = i_{a+1} \neq i_{a+1}$.
Let us assume that $ a$ appear above $a+1$ in $ \bt $ and $ a+1 $ appears
above $ a+2 $ in $ \bt$ and let $ \bt_d = \sigma_a \sigma_{a+1} \sigma_a \bt $. Then by
Corollary \ref{YSFcorA} and Corollary \ref{YSFcor}
we get for the left and right hand sides of 
\eqref{correspondstoproving} 
\begin{equation}\label{correspondstoprovingleft}
\begin{aligned}
  \F_a (g_{a+1} +1) \F_a E_{\bt} f_{\bt \bu} &=   \dfrac{X_a-\q X_{a+2}}{X_a-X_{a+2}} \times
\dfrac{(\q X_{a+1}-X_{a})(X_{a+1}-\q X_{a})}{(X_a-X_{a+1})^2 }f_{\bt \bu} + f_{\bt_d \bu} 
     \\[8pt]
    \F_{a+1} (g_{a} +1) \F_{a+1} E_{\bt} f_{\bt \bu} &=
\dfrac{X_a-\q X_{a+2}}{X_a-X_{a+2}} \times
\dfrac{(\q X_{a+2}-X_{a+1})(X_{a+2}-\q X_{a+1})}{(X_{a+1}-X_{a+2})^2 }f_{\bt \bu} + f_{\bt_d \bu} 
\end{aligned}  
\end{equation}
and using the identity 
\begin{equation}
  \dfrac{X_a-\q X_{a+2}}{X_a-X_{a+2}} \left(
  \dfrac{(\q X_{a+1}-X_{a})(X_{a+1}-\q X_{a})}{(X_a-X_{a+1})^2 } -
  \dfrac{(\q X_{a+2}-X_{a+1})(X_{a+2}-\q X_{a+1})}{(X_{a+1}-X_{a+2})^2 } \right) =Z_a 
  \end{equation}
we arrive at the claim. If $ a, a+1 $ and $ a+2 $ are positioned differently in $ \bt $ then
the coefficient of each $ f_{\bt_d \bu} $ in \eqref{correspondstoprovingleft} will be 
another scalar $ \alpha_{\bt_d \bu} \in \Kdomain(\q)$ but the coefficients of
the $f_{ \bt \bu} $'s remain the same and so the previous argument carries over. This finishes the
proof of the Theorem.
\end{dem}

\medskip
We need to generalize the previous results slightly. 
Let us introduce two new equivalence relations $ \sim_{q, || \cdot| |} $
and $ \sim_{q, \ord, || \cdot| |} $ on $ \std_{d, n} $ via
\begin{equation}
\begin{aligned}  
  \bs \sim_{q, || \cdot| |} \bt \, & \mbox{ if and only if } c_{\bs}(a) = c_{\bt}(a) \, \, \mbox{for } a=1,2,\ldots, n
  \mbox{ and } A_{ || \bs || } = A_{ || \bt || }
\\
 \bs \sim_{q,\ord,  || \cdot| |} \bt \, & \mbox{ if and only if } c_{\bs}(a) = c_{\bt}(a) \, \, \mbox{for } a=1,2,\ldots, n
  \mbox{ and } A_{ || \bs || }^{\ord} = A_{ || \bt || }^{\ord}.
\end{aligned}  
\end{equation}

The equivalence relation $  \sim_{q,t} $, defined in \eqref{equivalence qt},
may be reformulated in the following way 
\begin{equation}
  \bs \sim_{q,t} \bt \, \mbox{ if and only if } c_{\bs}(a) = c_{\bt}(a) \, \, \mbox{for } a=1,2,\ldots, n
  \mbox{ and }  || \bs ||  =  || \bt || 
\end{equation}  
and so we get that the equivalence classes for $\sim_{q, \ord, || \cdot| |} $ are disjoint unions of
equivalence classes for $  \sim_{q, t}     $. Moreover, the equivalence classes for
$\sim_{q,  || \cdot| |} $ are disjoint unions of
equivalence classes for $  \sim_{q, \ord, || \cdot| |     }     $, and therefore also disjoint unions
of equivalence classes for $  \sim_{q,t} $.

\medskip
Let $ \ii \in \II $. For $ A \in \SetPar_n $ and $ A^{\ord} \in \OrdSetPar_n $ we 
let $[ \ii, A] =  [ \ii, A]_{q, || \cdot| |} $ and $[ \ii, A^{\ord}] =  [ \ii, A]_{q, \ord, || \cdot| |} $ be the associated
equivalence classes for $ \sim_{q, || \cdot| |} $ and $ \sim_{q, \ord || \cdot| |} $, that is 
\begin{equation}
\begin{aligned}  
[ \ii, A] = \{ \bs \in \std_{d, n} \, |\, c_{\bs}(a) &= i_a \, \, \mbox{for } a=1,2,\ldots, n   \mbox{ and }
A_{ || \bs || } = A \} \\
[ \ii, A^{\ord}] = \{ \bs \in \std_{d, n} \, |\, c_{\bs}(a) & = i_a \, \, \mbox{for } a=1,2,\ldots, n   \mbox{ and }
A^{\ord}_{ || \bs || } = A^{ \ord} \} .
\end{aligned}  
\end{equation}
We then define
\begin{equation}\label{bytheaboveremarks}
\e( \ii, A) =  \sum_{ \bt \in [\ii, A]} E_{\bt} \mbox{      and      }  \e( \ii, A^{\ord}) =  \sum_{ \bt \in [\ii, A^{\ord}]} E_{\bt}. 
\end{equation}
By the above remarks on equivalence classes,
$ \e( \ii, A) $ and $ \e( \ii, A^{\ord}) $ are sums of $ E_{\TT} $'s 
and then, in view of \eqref{classidempotent}, 
also $ \e( \ii, A) $ and $ \e( \ii, A^{\ord}) $
are elements of $ \YYfractionfield $ and may therefore be considered as 
elements of $ \YYresiduefield $ via reduction modulo
$ (\q -q ) \OO $.

\medskip
Arguing as above, and using Lemma \ref{lemma13} and \eqref{firstwethendefine} we see 
that the sums for $ \e(A) $ and $ \e(A^{\ord} )$ also involve disjoint equivalence classes for
$  \sim_{q, t}     $ and so we have 
\begin{equation}
\e(A) \in \YYresiduefield \mbox{ for } A \in \SetPar_n \, \, \mbox{and} \, \, \, \, 
\e(A^{\ord}) \in \YYresiduefield \mbox{ for } A^{\ord} \in \OrdSetPar_n.
\end{equation}

We then get formulas
\begin{equation}\label{bytheaboveremarksA}
  \e( \ii, A) =  e(\ii) \e(A) \mbox{      and      }
  \e( \ii, A^{\ord}) =  e(\ii) \e(A^{\ord}) . 
\end{equation}

\medskip
Suppose that $  \TT $
is an equivalence class for $  \sim_{q, t}     $ and that $ \TT \subseteq [ \ii, A^{\ord}] $
for some $ \ii \in \II $ and $ A^{\ord} \in \OrdSetPar_n $.
Then for $ 1 \le a < n $ and $ \bt \in \TT $ the numbers 
$ a $ and $ a+1 $ appear in the same component of $ \bt $ if and only they appear in the same block
of $ A^{\ord} $ that is $ a \sim_A (a+1)$. The same observation also holds true for
$ \TT \subseteq [\ii, A] $
for some $ \ii \in \II $ and $ A \in \SetPar_n $.

\medskip
There is a left $ \Si_n$-action on pairs $ (\ii, A ) \in \II \times \SetPar_n $ given by
$ \sigma  (\ii, A) = (\sigma \ii, \sigma A)$ and similarly 
there is a left $ \Si_n$-action on pairs $ (\ii, A^{\ord} ) \in \II \times \OrdSetPar_n $ 
via 
$ \sigma  (\ii, A^{\ord}) = (\sigma \ii, \sigma A^{\ord})$ where $ (\sigma, A^{\ord} ) \mapsto \sigma A^{\ord} $
is the natural left $ \Si_n$-action on $ \OrdSetPar_n$.

\medskip
Using these observations and notations we can now reformulate
the intertwining elements $ \F_a $ from \eqref{intertwinedef} as follows

\begin{equation}\label{intertwineE}
\begin{aligned}
  \F_a &=
   g_a + (1-q) \sum_{\substack{\ii \in \II \\ A \in \SetPar_n \\  i_a \neq i_{a+1} \\   a \sim_A (a+1) }}
  (1-X_a X_{a+1}^{-1})^{-1} \e(\ii, A) +
  \sum_{\substack{  \ii \in \II \\ A \in \SetPar_n \\  i_a = i_{a+1} \\   a \sim_A (a+1)  }} \e(\ii,A) \\
&=
 g_a + (1-q) \sum_{\substack{\ii \in \II \\ A^{\ord} \in \OrdSetPar_n \\  i_a \neq i_{a+1} \\   a \sim_A (a+1) }}
  (1-X_a X_{a+1}^{-1})^{-1} \e(\ii, A^{\ord}) +
  \sum_{\substack{  \ii \in \II \\ A^{\ord} \in \OrdSetPar_n \\  i_a = i_{a+1} \\   a \sim_A (a+1)  }} \e(\ii,A^{\ord}).
\end{aligned}  
\end{equation}

\medskip
We get the following analogue of Theorem \ref{commutation intertwiner}.
\begin{theorem}\label{commutation intertwiner setpar}
Suppose that $ \ii \in \II  $ and $ A \in \SetPar_n $.
Then for 
$ 1 \le a < n  $
we have that
\begin{align}
\label{intertwine one setpar} \F_a  \e(\ii, A)  & = \e(\sigma_a ( \ii, A)) \F_a            \\
\label{intertwine two setpar} \F_a  X_b  & = X_b \F_a  \, \, \, \, \, \, \, \, \, \, \, \,     \mbox{if } b  \neq  a, a+1        \\
 \label{intertwine three setpar} \F_a  \F_b  & = \F_b \F_a \, \, \, \, \, \, \, \, \, \, \, \, \,  \mbox{if } | a-b | > 1           \\
 \label{intertwine four setpar} \F_a X_{a+1} e(\ii,A)    &   =  \left\{\begin{aligned} & X_a \F_a \e(\ii,A) & &  \mbox{if }
(i_a \neq i_{a+1}, \,   a \sim_A (a+1)) \mbox{ or }  a \not\sim_A (a+1) \\
&  (X_a \F_a +q X_{a+1} -X_a)  \e(\ii,A) & & \mbox{if } i_a = i_{a+1}, \, a \sim_A (a+1) \\
 \end{aligned} \right.    \\
\label{intertwine five setpar}  X_{a+1} \F_a \e(\ii,A)    &   =  \left\{\begin{aligned} & X_a \F_a \e(\ii,A) & &  \mbox{if }
(i_a \neq i_{a+1}, \, a \sim_A (a+1)) \mbox{ or } a \not\sim_A (a+1) \\
&  ( \F_a X_a +q X_{a+1} -X_a)  \e(\ii,A) & & \mbox{if } i_a = i_{a+1}, \, a \sim_A (a+1)  \end{aligned} \right.    \\
\label{intertwine six setpar}  \F_a^2 \e(\ii,A)    &   =  \left\{\begin{aligned} &
\frac{(q X_{a} -X_{a+1}) (X_{a} -qX_{a+1})}{(X_a-X_{a+1})^2} \e(\ii,A)  & &  \mbox{if }
i_a \neq i_{a+1}, \, a \sim_A (a+1) \\
&  (q+1) \F_a \e(\ii,A) & & \mbox{if } i_a = i_{a+1}, \,  a \sim_A (a+1)  \\
&  q \F_a \e(\ii,A) & & \mbox{if } a \not\sim_A (a+1)
\end{aligned} \right.    \\
\label{intertwine seven setpar} \F_a \F_{a+1} \F_{a}   \e(\ii,A)    &   =  \left\{\begin{aligned}
& (\F_{a+1} \F_{a} \F_{a+1}  +q \F_a -q \F_{a+1} ) \e(\ii,A) & &  \mbox{if } i_a = i_{a+1} = i_{a+2}, \,
a \sim_A (a+1) \sim_A (a+2)
\\
& (\F_{a+1} \F_{a} \F_{a+1}+Z_a) \e(\ii,A) & &  \mbox{if }  i_a = i_{a+2} \neq i_{a+2}    , \,
a \sim_A (a+1) \sim_A (a+2)       \\
& \F_{a+1} \F_{a} \F_{a+1} \e(\ii,A) & &  \mbox{otherwise}  \\
\end{aligned} \right.    
\end{align}
where $ Z_a = (1-q)^2\frac{ (X_aX_{a+2} -X_{a+1}^2) (X_a X_{a+1} -q X_{a+1} X_{a+2})}
{(X_a-X_{a+1})^2(X_{a+1}-X_{a+2})^2}$.

\medskip
There are similar formulas with $ A \in \SetPar_n $ replaced by $ A^{\ord} \in \OrdSetPar_n$.
\end{theorem}





\section{The bt-algebras $ \E  $ and $ \Eord$. }\label{bt-algebra}
In this section we introduce two of the main objects of the present paper, namely
the {\it bt-algebras} $\E$ and $\Eord$.

\medskip
The original bt-algebra $ \E$ was defined by Aicardi and Juyumaya in \cite{AicardiJuyumaya1}.
There are several slightly different presentations of it, but we shall use the following one. 

\begin{definition}\label{btdefinition}
  \normalfont
Let $ \Kdomain $ be a field and let $ q \in \Kdomain^{\times}$. Then the bt-algebra 
$ \E = \Efield $ is defined as the $ \Kdomain$-algebra generated by
  $ \{g_1, g_2,\ldots, g_{n-1} \} $ and $ \{ e_1, e_2, \ldots, e_{n-1}\} $ subject to the following relations
\begin{align}
  \label{eq Eone} g_a e_a    & = e_a g_a    &   &   \mbox{for all }  a   \\
\label{eq Etwo} g_a g_b g_a   & = g_b g_a g_b   &   &   \mbox{for }  | a-b | = 1   \\
\label{eq Ethree} g_a g_b    & = g_b g_a       &   &  \mbox{for }  | a-b | > 1   \\
\label{eq Efour} e_a g_b g_a   & = g_b g_a e_b   &   &     \mbox{for }  | a-b | = 1   \\
\label{eq Efive} e_a e_b g_b   & = e_a g_b e_a = g_be_a e_b   &   &     \mbox{for }   | a-b | = 1   \\
\label{eq Eseven} e_a e_b   & = e_b e_a     &   &  \mbox{for all }  a, b  \\
\label{eq Eeight} g_a e_b   & = e_b g_a     &   &     \mbox{for }   | a-b | > 1   \\
\label{eq Enine} e_a^2   & = e_a     &   & \mbox{for all }  a \\
\label{eq Ten} g_a^2  & = q +  (q-1) e_a g_a   &  &  \mbox{for all }  a.        
        \end{align}
\end{definition}

\medskip
Note that the definition of $ \Efield $ does not require any restrictions on the characteristic of 
$ \Kdomain$, contrary to the definition of $ \YYresiduefield $.

\medskip
Suppose that $ \Kdomain $ is a field for which $ \YYresiduefield $ is defined. 
Then there is a homomorphism $ \iota: \Efield \rightarrow \YYresiduefield $ given by
$ \iota(g_a) = g_a $ and $  \iota(e_a) = e_a $.
In Theorem 14 of \cite{ERH} it was shown
that $ \iota $ is injective if $ d \ge n$ and $ \Kdomain $ contains
a primitive $ d$'th root of unity, but in general $ \iota $ is not injective.

\medskip
Note that $ g_a \in \E$ is invertible, with inverse $ g_a^{-1} = 
q^{-1} g_a +(q^{-1}-1) e_a $.
For $ 1 \le a < b \le n $ we can therefore
mimic \eqref{e_{ab}} to introduce elements $ e_{ab}   \in \E $, denoted the same way as
$ e_{ab} \in  \YY$, and similarly we introduce $ e_{ba} =  e_{ab} $. 
Mimicking \eqref{e_A}, we can then also introduce elements $ e(A) \in \E$, for $A \in \SetPar_n $.

\medskip
With these definitions, if $ \Kdomain $ a field such that $ \YYresiduefield $ is defined, we have
\begin{equation}\label{mimicking}
\iota(e_a) = e_{a}, \, \, \, \iota(e(A)) = e(A)
\end{equation}
for $ 1 \le a \le n $ and $ A \in \SetPar_n$. 

\medskip
For $ w= \sigma_{i_1} \sigma_{i_2} \cdots \sigma_{i_N} \in \Si_n  $ a reduced expression
for $ w $
we define $ g_{w} = g_{i_1} g_{i_2} \cdots g_{i_N} $. 
By \eqref{eq Etwo}, this is independent of the choice of reduced expression. 
In \cite{RH1} it was shown that the set $ \{ g_w e(A) | w \in \Si_n, A\in \SetPar_n \} $ is
a basis for $ \E$. In particular, $ \dim \E = b_n n !$. 
  
\medskip
The $ e(A)$'s are idempotents in $ \E$, but they are not orthogonal idempotents. 
To remedy this, following  \cite{ERH}
we introduce orthogonal idempotents $ \{\e(A) \, | \, A \in \SetPar_n \} $ for $ \E$ as follows. 
Let $ (\SetPar, \subseteq) $ be the poset structure introduced in  
the paragraph following \eqref{overline e(A)first}. 
We then define the $ \e(B) $'s as the unique solutions 
of the equations
\begin{equation}\label{uniquesolutions}
 e(A) =  \sum_{\substack{  B \in \SetPar_n \\  A \subseteq B}} \e(B).
\end{equation}
They form a complete set of orthogonal idempotents for $ \E$, that is
\begin{equation}\label{complete}
\e(A) \e(B) = \delta_{AB} \e(A) \, \, \, \, \,  \mbox{and} \, 
\sum_{ A \in \SetPar_n} \e(A) = 1. 
\end{equation}

\medskip
To express the $ \e(A) $'s in terms of the $ e(A)$'s, we use M\"obius inversion as in section 
\ref{idempotents in YY}. 
We get 
\begin{equation}\label{mobiusinversionA}
 \e(A) =  \sum_{\substack{  B \in \SetPar_n \\  A \subseteq B}} \mu(A,B) e(B).
\end{equation}
Since $ \mu(A,B) \in \Z $ 
it follows from \eqref{mobiusinversionA} that $ \e(A) \in \E$, as claimed.
Moreover, if $ \Kdomain $ a field such that $ \YYresiduefield $ is defined we get
via
\eqref{firstMobius inversion} that
\begin{equation}
\iota(\e(A)) = \e(A). 
\end{equation}

\medskip
For general fields $ \Kdomain$, we have the following Lemmas.
\begin{lemma}\label{afirstlemmaAA}
  For $ 1\le a \le n-1 $ and $ A \in \SetPar_n $ the following holds in $ \E$. 
  \begin{description}
\item[a)] $g_a e(A) = e(\sigma_a  A) g_a$
\item[b)] $ g_a \e(A) = \e( \sigma_a  A) g_a$.
\end{description}
\end{lemma}
\begin{dem}
  $   { \bf a)} $ is Corollary 1 of \cite{RH1}, and $   { \bf b)} $ follows from
  it via \eqref{mobiusinversionA}.
\end{dem}

\begin{lemma}\label{lemma28}
Suppose that $ A \in \SetPar_n $ and that $a \sim_A (a+1) $. Then the following holds in 
  $ \E $
  \begin{equation}  g_a^2 \e(A) =  (q + (q-1) g_a) \e(A).
 \end{equation}    
\end{lemma}
\begin{dem}
The proof is the same as the proof of part $   { \bf a)} $ of Lemma \ref{veryfirstusefullemma}.
\end{dem}

\begin{lemma}\label{lemma29}
Suppose that $ A \in \SetPar_n $ and that $a \not\sim_A (a+1) $. Then the following identity holds in 
  $ \E $ 
  \begin{equation}  g_a^2 \e(A) =  q \e(A). 
 \end{equation}    
\end{lemma}
\begin{dem}
This is proved the same way as $   { \bf a)} $ of Lemma \ref{usefullemma}.
\end{dem}

\medskip

Using the above Lemmas we get the following Theorem that we need later on.

\begin{theorem}\label{ortoBT}
Let $ \Kdomain $ be a field and let $ q \in \Kdomain^{\times}$. Then $ \E = \Efield $ is 
isomorphic to the $ \Kdomain$-algebra generated by
  $ \{g_1, g_2,\ldots, g_{n-1} \} $ and $ \{ \e(A) | A \in \SetPar_n \} $,  subject to the following relations
\begin{align}
\label{eqbt-1} \sum_{ A \in \SetPar_n} \e(A)    & = 1   &   &      \\
\label{eqbt0} \e(A) \e(B)   & = \delta_{A B}  \e(A)   &   &   \mbox{for all }  A, B \in \SetPar_n   \\
\label{eqbt1} g_a \e(A)    & = \e(\sigma_a A) g_a    &   &     \mbox{for all }  1 \le a < n  \mbox{ and } A \in \SetPar_n    \\
\label{eqbt8} g_a^2  \e(A) & = q \e(A)     &   & \mbox{if  }  a \not\sim_A (a+1) \\
\label{eqbt9}  g_a^2  \e(A) & = (q+(q-1)g_a) \e(A)     &  & \mbox{if  }  a \sim_A (a+1)  .     \\
\label{eqbt2} g_a g_b g_a   & = g_b g_a g_b   &   &   \mbox{for }  | a-b | = 1   \\
\label{eqbt3} g_a g_b    & = g_b g_a       &   &  \mbox{for }  | a-b | > 1. 
\end{align}
\end{theorem}
\begin{dem}
  Letting $ \Eprime $ be the algebra defined by \eqref{eqbt0}--\eqref{eqbt9} we must show that
  $ \Eprime \cong \E $. Let $ F: \Eprime \rightarrow \E $ be the homomorphism given by
\begin{equation}
 F(g_a ) = g_a \mbox{ and   }  F(\e(A) ) =    \sum_{\substack{  B \in \SetPar_n \\  A \subseteq B}} \mu(A,B) e(B). 
\end{equation}
To show that $ F $ exists we must verify that the defining relations for $ \Eprime $ hold when
$ g_a $ and $ \e(A)  $ are replaced by $ F(g_a) $ and $F( \e(A))  $. For this we first observe that \eqref{eqbt-1}
and    
\eqref{eqbt0} follow from \eqref{complete}, whereas \eqref{eqbt1} follows from $   { \bf b)} $ of 
Lemma \ref{afirstlemmaAA}. We get \eqref{eqbt8} and \eqref{eqbt9} from
Lemma \ref{lemma28} and Lemma \ref{lemma29} and finally \eqref{eqbt2} and \eqref{eqbt2} follow 
directly from \eqref{eq Etwo} and \eqref{eq Ethree}.

\medskip
For $ 1 \le a < n $ we set as before
$ A_a=\{ \{1\}, \{2\}, \ldots, \{a,a+1\}, \ldots, \{n-1\},  \{n\} \} \in \SetPar_n$. 
Then $ e_a =  e(A_a) \in \E $. 
We then define a homomorphism $ G: \E \rightarrow \Eprime $ as follows
\begin{equation}\label{Gga}
 G(g_a ) = g_a \mbox{ and   }  G(e_a)  =    \sum_{\substack{  B \in \SetPar_n \\  A_a \subseteq B}}  \e(B). 
\end{equation}
To show that $ G $ exists we must verify that the defining relations for $ \E $ hold when
$ g_a $ and $ e_a  $ are replaced by $ G(g_a) $ and $G( e_a)  $. 
Here \eqref{eq Etwo}, \eqref{eq Ethree} and \eqref{eq Eseven} are immediate consequences of
the definitions. Relation \eqref{eq Eone} follows from \eqref{Gga} and 
$ \sigma_a A_a = A_a $ and \eqref{eq Efour} follows from \eqref{Gga} and the fact that
for any $ B \in \SetPar_n $ we have
$ A_a \subseteq B $ if and only if $ \sigma_b \sigma_a A_b \subseteq B $ whenever
$ | a-b | = 1 $, and the relations \eqref{eq Efive} and \eqref{eq Eeight} are proved the same way.
Relation \eqref{eq Enine} follows from the fact that a sum of orthogonal idempotents is an idempotent. 
Finally, relation \eqref{eq Ten} is proved as follows
\begin{equation}
\begin{aligned}
  g_a^2 &=  \sum_{ A \in \SetPar_n}  g_a^2 \e(A) =
  \sum_{ \substack{  A \in \SetPar_n \\  a \sim_A (a+1) }}  g_a^2 \e(A)  +
  \sum_{ \substack{  A \in \SetPar_n \\  a \not\sim_A (a+1) }}  g_a^2 \e(A)   \\
& =\sum_{ \substack{  A \in \SetPar_n \\  a \sim_A (a+1) }}  (q+(q-1)g_a \e(A)  +
  \sum_{ \substack{  A \in \SetPar_n \\  a \not\sim_A (a+1) }}  q \e(A)   \\
& =   \sum_{ \substack{  A \in \SetPar_n }}  q \e(A) +
  (q-1)g_a  \sum_{ \substack{  A \in \SetPar_n \\  a \sim_A (a+1) }}   \e(A)   
 \\
  & = q +
  (q-1)g_a  G(e_a) .   
    \end{aligned}
\end{equation}

\medskip \medskip
Clearly $ F \circ G = Id_{ \E} $ and so $ G $ is injective. On the other hand, it follows from
the relations \eqref{eqbt-1}--\eqref{eqbt3} that $\{  g_w \e(A) | w \in \Si_n, A \in \SetPar_n \}  $
is a generating set for $ \Eprime $ and so $ \dim \Eprime \le b_n n! $, which shows that
$ G $ is indeed an isomorphism, as claimed. 
\end{dem}



\medskip \medskip
There is a diagram calculus associated with $\E $
which plays an important role
in applications of $ \E $ to knot theory, see for example \cite{AicardiJuyumaya1},
\cite{AicardiJuyumaya}, \cite{RH1}, \cite{RH2}. 
That diagram calculus is derived from
the relations \eqref{eq Eone}--\eqref{eq Ten}, but the relations 
\eqref{eqbt-1}--\eqref{eqbt3} can be used to derive another diagram calculus for $ \E$, that we shall now explain.

\medskip
A {\it bt-diagram} $ D $ for $ \E$, with respect to the relations \eqref{eqbt-1}--\eqref{eqbt3},
is a braid $w$ 
on $ n$ strands joining the northern points $ \{1, 2, \ldots, n \} $
with the southern points $ \{1^{\prime}, 2^{\prime}, \ldots, n^{\prime} \} $
such that each pairs of points in $ \{1^{\prime}, 2^{\prime}, \ldots, n^{\prime} \} $
may be connected by a \lq tie\rq\ {\!.} The set of these ties define a 
set partition $ A_D \in \SetPar_n$.
Here are two examples of bt-diagrams for
$ {\mathcal E}_5(q) $ with associated set partitions $A_D= \{ \{1,3\},  \{2,4,5\} \} $ and
$A_E =\{ \{1,2,3\},  \{4,5\} \} $, respectively.
\begin{equation}\label{btdiagram}
D=  \raisebox{-.45\height}{\includegraphics[scale=0.7]{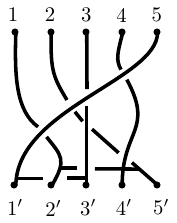}} \raisebox{-20\height}{,} \, \, \, \, 
E=  \raisebox{-.45\height}{\includegraphics[scale=0.7]{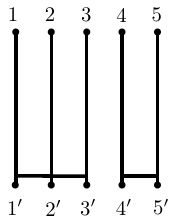}}
\raisebox{-20\height}{.}
\end{equation}
Let $ D $ be a bt-diagram with associated set partition $ A_D $.
Then there is a natural \lq upper\rq\ set partition $ A_D^{up} $ associated with $ D$, obtained from
$ A_D $ by pulling all the ties up to $ \{1, 2, \ldots, n \} $. 
For example, for $ D$ and $ E $ in \eqref{btdiagram} we have
$ A_D^{up} = \{ \{1,5\},  \{2,3,4\} \} $ and $ A_D^{up} = B  =\{ \{1,2,3\},  \{4,5\} \}$.

\medskip
Let $  {\mathcal E}^{diag}_n(q) $ be the free $ \Kdomain $-vector space on
the set of all bt-diagrams for $ \E$. 
There is a multiplication in $  {\mathcal E}^{diag}_n(q) $ given 
by vertical 
concatenation of bt-diagrams, that is for bt-diagrams $ D$ and $ D_1 $ 
the product $ D D_1  $ is the join of the diagrams $ D $ and $ D_1 $ with $ D_1 $ on top
of $ D$. The product $ D D_1 $ is defined to be zero if the lower set partition for $ D $ does not coincide
with the set partition for $ D_1$. We further impose relations in 
$ {\mathcal E}^{diag}_n(q) $
corresponding to \eqref{eqbt8}--\eqref{eqbt3}. We illustrate them as follows. 
\begin{equation}\label{illustrate1}
  \raisebox{-.45\height}{\includegraphics[scale=0.7]{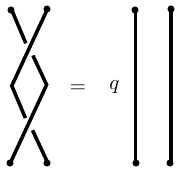}} \raisebox{-20\height}{,} \, \,  \, \,  \, \,  \, \, 
  \raisebox{-.43\height}{\includegraphics[scale=0.7]{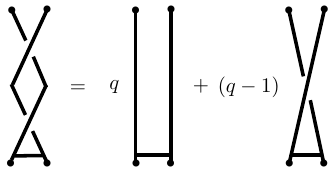}}
  \end{equation}  
\begin{equation}\label{529}
  \raisebox{-.45\height}{\includegraphics[scale=0.7]{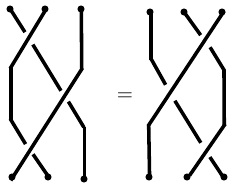}} \raisebox{-20\height}{,} \, \,  \, \,  \, \,  \, \, 
   \raisebox{-.43\height}{\includegraphics[scale=0.7]{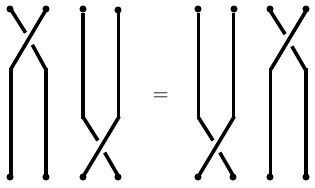}}
  \raisebox{-20\height}{.}
\end{equation}  
Finally, the one-element $ 1 $ in $ {\mathcal E}^{diag}_n(q) $ is the sum
of all straight vertical line diagrams with all possible associated set partitions, for example
\begin{equation}\label{illustrate3}
1=  \raisebox{-.45\height}{\includegraphics[scale=0.7]{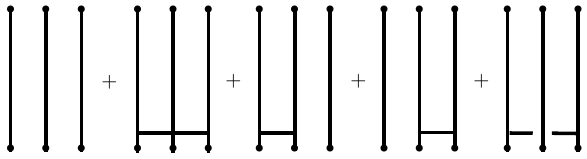}} \in  {\mathcal E}^{diag}_3(q) .
\end{equation}

\medskip
With these definition, one obtains an algebra isomorphism $ \E \cong  {\mathcal E}^{diag}_n(q) $.

\medskip
Partly inspired by Theorem \ref{ortoBT} we now introduce the {\it ordered bt-algebra} $ \Eord$.

\begin{definition}\label{orderedBT}
  \normalfont
Let $ \Kdomain $ be a field and let $ q \in \Kdomain^{\times}$. Then the ordered bt-algebra
$ \Eord = \Eordfield $ is defined as the $ \Kdomain$-algebra generated by
$ \{g_1, g_2,\ldots, g_{n-1} \} $ and
$ \{ \e(A^{\ord}) | A^{\ord} \in \OrdSetPar_n \} $
subject to the following relations
\begin{align}
\label{ordeqbt-1} \sum_{ A^{\ord} \in \OrdSetPar_n} \e(A^{\ord})    & = 1   &   &      \\
\label{ordeqbt0} \e(A^{\ord}) \e(B^{\ord})   & = \delta_{A^{\ord} B^{\ord}}  \e(A^{\ord})   &   &   \mbox{for all }  A^{\ord}, B^{\ord} \in \OrdSetPar_n   \\
\label{ordeqbt1} g_a \e(A^{\ord})    & = \e(\sigma_a A^{\ord}) g_a    &   &     \mbox{for all }  1 \le a < n  \mbox{ and } A^{\ord} \in \OrdSetPar_n    \\
\label{ordeqbt8} g_a^2  \e(A^{\ord}) & = q \e(A^{\ord})     &   & \mbox{if  }  a \not\sim_{A} (a+1) \\
\label{ordeqbt9}  g_a^2  \e(A^{\ord}) & = (q+(q-1)g_a) \e(A^{\ord})     &  & \mbox{if  }  a \sim_{A}(a+1)  .     \\
\label{ordeqbt2} g_a g_b g_a   & = g_b g_a g_b   &   &   \mbox{for }  | a-b | = 1   \\
\label{ordeqbt3} g_a g_b    & = g_b g_a       &   &  \mbox{for }  | a-b | > 1. 
        \end{align}
\end{definition}

\medskip
\begin{theorem}\label{orderedBTinclusion}
 Suppose that $ \Kdomain $ is a field for which $ \YY = \YYresiduefield  $ exists. 
Suppose furthermore that $ d \ge n $ and that $ \Kdomain $ contains a $ d$'th root of unity. Then
there is a $ \Kdomain $-algebra inclusion $ \epsilon: \Eord \rightarrow \YYresiduefield $ given by 
\begin{equation}
 \epsilon(g_a ) = g_a \mbox{   and    } \, \, \epsilon(\e(A^{\ord})) = \e( A^{\ord}) 
\end{equation}
and a $ \Kdomain $-algebra inclusion
$ \epsilon_1: \E \rightarrow \Eord$ given by 
\begin{equation}
  \epsilon_1(g_a ) = g_a \mbox{   and    } \,\, \epsilon_1(\e(A) ) =
  \sum_{\substack{  B^{\ord} \in \OrdSetPar_n \\  B=A } } \e(B).
\end{equation}
\end{theorem}
\begin{dem}
  Using Lemma \ref{veryfirstusefullemmaA} and \ref{usefullemmaA} together with 
  the definitions one checks that $ \epsilon $ and $\epsilon_1 $ are well-defined.
  On the other hand, we have that $ \epsilon \circ \epsilon_1 = \iota$
  where $ \iota: \E \rightarrow \YYresiduefield $ is as above, and since $ \iota $ is injective, we conclude that
  also $  \epsilon_1 $ is injective, as claimed. Finally, to show the injectivity of $ \epsilon $ we may argue the same way 
  as in \cite{ERH}, where the injectivity of $ \iota$ is shown.
  
\end{dem}  

\begin{remark}
  \normalfont
  It follows from the proof of the Theorem that
  \begin{equation}
\dim \Eord = n! f_n 
  \end{equation}    
where $ f_n $ is the Fubini number introduced in the paragraph after \eqref{ordered set partitions}.
\end{remark}  
  
There is also a diagram calculus associated with $ \Eord$, given by {\it ordered bt-diagrams}.
An ordered bt-diagram for $\Eord$ is a 
usual bt-diagram for $ \E $ except that the $ b $ blocks, say, of the associated set partition $ A $ are decorated 
with the numbers $ \{ 1,2, \ldots, b \} $, each number appearing once. Here are two examples
of ordered bt-diagrams for ${\mathcal E}^{\ord}_5(q)$.

\begin{equation}\label{541}
D=  \raisebox{-.45\height}{\includegraphics[scale=0.7]{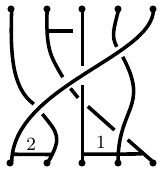}} \raisebox{-20\height}{,} \, \,  \, \,  \, \,  \, \, 
E=   \raisebox{-.45\height}{\includegraphics[scale=0.7]{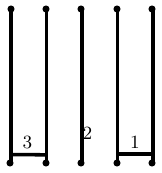}} \raisebox{-20\height}{.}
\end{equation}  
The decorations on $ A $ give rise to an ordered set partition $A^{\ord} $. In the examples
$ D $ and $ E $ of \eqref{541} the associated ordered set partitions are $ \{\{ 3,4,5 \}, \{ 1,2 \}\} $
and $\{ \{ 4,5 \}, \{3\}, \{ 1,2 \}\} $, respectively.

\section{Quiver Hecke algebras}\label{qHeckealgebras}
In this section we recall the quiver Hecke algebra, or KLR-algebra, for the specific
quiver $ \Gamma_{e,d} $ introduced by Rostam, see \cite{Ro}.

\medskip
Suppose that $ \Kdomain $ is a field and that $ q \in \Kdomain^{\times} $. 
Let $ e = {\rm char }_q (\Kdomain) $ and assume for simplicity that $ e \neq  2$. 
Let $ I = \Z/ e  \Z    $, $ J = \Z/d  \Z $ and $ K = I \times J$, as in section
\ref{idempotents in YY}.
For $ i, j \in I $ we write
\begin{equation}
\begin{array}{rl}
   i  \rightarrow   i_1 	    &     \mbox{if  }  i_1 = i+1 \\
   i  \leftarrow  i_1 	    &     \mbox{if  }  i_1 = i-1  \\
   i  \nrelbar  i_1	    &     \mbox{if  }  i_1 \not\in \{i, i +1, i-1 \}. 
\end{array}
\end{equation}

We extend this notation to $ K $ as follows. For $ (i,j), (i_1, j_1 ) \in K $ we write
\begin{equation}
\begin{array}{rl}
   (i,j)  \rightarrow   (i_1, j_1)  	    &     \mbox{if  }  i \rightarrow i_1  \mbox{ and }  j=j_1\\
   (i,j)  \leftarrow  (i_1, j_1) 	    &     \mbox{if  }  i \leftarrow i_1\mbox{ and }  j=j_1  \\
   (i,j)  \nrelbar  (i_1, j_1)	    &     \mbox{if  } i \nrelbar i_1 \mbox{ or if }  j \neq j_1.
\end{array}
\end{equation}

Let $ \Gamma_e $ be the quiver with vertices labeled by the elements of $ I $ and an arrow from
the vertex labeled by $ i $ to the vertex labeled by $ i_1 $ if and only if $ i \rightarrow i_1 $. 
For example
\begin{equation}
\Gamma_{6}= \raisebox{-.45\height}{\includegraphics[scale=0.7]{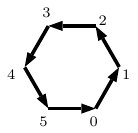}} 
\, \, \, \,  \, \, \, \,  \, \, \, \,  \begin{array}{l} 1 \rightarrow 2 \\ 0 \leftarrow 5\\  1  \nrelbar 3. \end{array}
\end{equation}

Similarly, let $ \Gamma_{e,d} $ be the quiver with vertices labeled by $ K $ and an arrow from the vertex labeled by $ (i,j) $ to
the vertex labeled by $ (i_1,j_1) $
if and only if $ (i, j) \rightarrow (i_1, j_1) $. In other words, $ \Gamma_{e,d} $ is a disjoint union of $d $ copies of $ \Gamma_e $. For example
\begin{equation} \Gamma_{6,3}= \raisebox{-.45\height}{\includegraphics[scale=0.7]{dib52.pdf}}
\raisebox{-.45\height}{\includegraphics[scale=0.7]{dib52.pdf}}
\raisebox{-.45\height}{\includegraphics[scale=0.7]{dib52.pdf}}
\, \, \, \,  \, \, \, \,  \begin{array}{l} (1,1) \rightarrow (2,1)
  \\ (1,1) \nrelbar (2,2). \end{array}
\end{equation}

\medskip
Recall the residue sequences $ \II, \JJn $ and $ \KKn$, introduced in the paragraph following
\eqref{quantum}. 
Recall from \eqref{placeperm}
that $ \Si_n $ acts on $ {\II } $, $ {\JJn } $ and $ {\KKn } $ by place 
permutation.

\begin{definition}\label{cycloHecke}
  \normalfont
Let $ \Kdomain $ be a field. Then the 
cyclotomic quiver Hecke algebra $  \RR $ is defined as the associative $ \Kdomain$-algebra on generators
  $ \{ e( \kk) | \kk \in \KKn    \}  \cup \{ y_1, y_2, \ldots, y_n \} \cup \{ \psi_1, \psi_2, \ldots, \psi_{n-1} \} $ subject to the following relations
\begin{align}
  \label{eq RcycloA} y_1 e(\kk)  & = 0     &  &    \mbox{if }  k_{1,i} = 0  \\
\label{eq RcycloB} e(\kk)  & =  0    &  &    \mbox{if }  k_{1,i} \neq 0  \\
  \label{eq Rone}	\sum_{\kk \in \KKn} e(\kk) 	& =  1   &   &   \\
\label{eq Rtwo}	e(\kk) e(\kk^{\prime})   & = \delta_{\kk, \kk^{\prime}} e(\kk) &  &   \\
\label{eq Rthree} y_a e(\kk)  & = e(\kk) y_a     &  &    \\
\label{eq Rfour} \psi_a e(\kk)  & = e(\sigma_a \kk) \psi_a     &  &       \\
\label{eq Rfive} y_a y_b   & = y_b y_a     &  &      \\
\label{eq Rsix} \psi_a \psi_b   & = \psi_b \psi_a     &  &   \mbox{if }     | a-b | > 1     \\
\label{eq Rseven} \psi_a y_{a+1} e(\kk)    &   =  \left\{\begin{aligned} & (y_a \psi_a +1)e(\kk) & &  \mbox{if }
k_a = k_{a+1} \\
&  y_a \psi_a e(\kk) & & \mbox{if } k_a \neq k_{a+1}  \end{aligned} \right. & &  \\
\label{eq Reight} y_{a+1} \psi_a  e(\kk)    &   =  \left\{\begin{aligned} & (\psi_a y_a  +1)e(\kk) & &  \mbox{if }
k_a = k_{a+1} \\
&  \psi_a y_a  e(\kk) & & \mbox{if } k_a \neq k_{a+1}  \end{aligned} \right. & &  \\
\label{eq Rnine} \psi_a^2  e(\kk)    &   =  \left\{\begin{aligned}
& 0 & &  \mbox{if } k_a = k_{a+1} \\
& e(\kk) & &  \mbox{if } k_a \nrelbar k_{a+1} \\
& (y_{a+1} -y_a)e(\kk) & &  \mbox{if } k_a \rightarrow k_{a+1} \\
& (y_{a} -y_{a+1})e(\kk) & &  \mbox{if } k_a \leftarrow k_{a+1} \\
\end{aligned} \right. & & \\
\label{eq Rten} \psi_{a+1} \psi_{a} \psi_{a+1}   e(\kk)    &   =  \left\{\begin{aligned}
& (\psi_{a} \psi_{a+1} \psi_{a} -1) e(\kk) & &  \mbox{if } k_{a+2} = k_{a} \rightarrow k_{a+1} \\
& (\psi_{a} \psi_{a+1} \psi_{a} +1) e(\kk) & &  \mbox{if } k_{a+2} = k_{a} \leftarrow k_{a+1}       \\
& (\psi_{a} \psi_{a+1} \psi_{a}) e(\kk) & &  \mbox{otherwise.}  \\
\end{aligned} \right. & &   
        \end{align}
  \end{definition}

\medskip
An important feature of $ \RR $ is the fact that it is a $ \Z$-graded algebra. To be
precise, the rules 
\begin{align}
  \label{deg1 }  \deg y_a & = 2     &  &   \\
\label{deg2 }  \deg e(\kk) &= 0     &   &   \\
\label{deg3 }  \deg \psi_a e(\kk) &= 
 \left\{\begin{aligned}
 0  &  & \mbox{if } & k_{a}  \nrelbar k_{a+1}  \\
  1  &  & \mbox{if } & k_a \leftarrow k_{a+1} \mbox{ or } k_a \rightarrow k_{a+1}        \\
  -2  &  & \mbox{if } & k_a = k_{a+1}        \\
\end{aligned} \right. & &  
\end{align}
define a $\Z$-degree function on $ \RR $, as one easily checks on the relations
\eqref{eq RcycloA}--\eqref{eq Rten}.

\medskip
There is a diagram calculus associated with $  \RR $ that we now briefly explain.
A KLR-diagram $ D $ for $ \RR $ is a braid, without over/under information, 
involving $ n$ strands that connect $ n $ northern points 
with $ n $ southern points. Each strand of $ D $ is decorated
by an element of $ K $ and the segments of
each strand 
may be decorated by one or several \lq dots\rq\ {\!\!.}
Here are two examples of KLR-diagrams for $ n = 5, e=3$ and $ d = 2$. 
\begin{equation}
  \raisebox{-.45\height}{\includegraphics[scale=0.7]{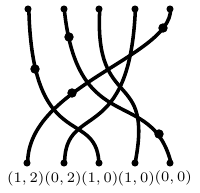}} \raisebox{-20\height}{,} \, \,  \, \,  \, \,  \, \,
 \raisebox{-.45\height}{\includegraphics[scale=0.7]{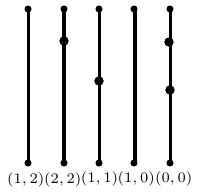}} \raisebox{-20\height}{.}
\end{equation}  
The diagram calculus realization for $ \RR $ is the $\Kdomain $-span of all KLR-diagrams for $\RR$,
with the multiplication $ D D_1 $ of diagrams $ D $ and $D_1 $ given by vertical concatenation with $ D $ on top of $ D_1 $
where $ D D_1 $ is set equal to $ 0 $ if the $\KKn$-decorations on $ D $ and $ D_1 $ do not match.
The elements $ e(\kk), y_a e(\kk)$ and $ \psi_a e(\kk) $ generate $ \RR $, and correspond to the following diagrams 
\begin{equation}
e(\kk) \mapsto  \raisebox{-.45\height}{\includegraphics[scale=0.7]{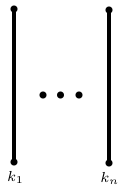}} \raisebox{-20\height}{,} \, \,  \, \,  \, \,  \, \,
y_a e(\kk) \mapsto   \raisebox{-.45\height}{\includegraphics[scale=0.7]{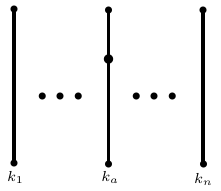}} \raisebox{-20\height}{,}  \, \,  \, \,  \, \,  \, \,
y_a e(\kk) \mapsto   \raisebox{-.45\height}{\includegraphics[scale=0.7]{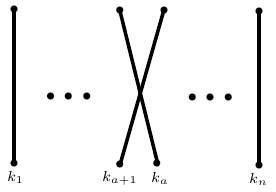}} \raisebox{-20\height}{.}  
\end{equation}  
Using this, the relations \eqref{eq RcycloA}--\eqref{eq Rten} can now be described diagrammatically.
For $ k_a =k_{a+1}$ 
the relations 
\eqref{eq Rseven} and \eqref{eq Reight} look as follows
\begin{equation}\label{REdiag1}
  \raisebox{-.45\height}{\includegraphics[scale=0.7]{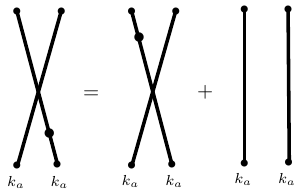}} \raisebox{-20\height}{,} \, \,  \, \,  \, \,  \, \,
   \raisebox{-.45\height}{\includegraphics[scale=0.7]{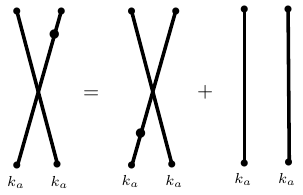}}
\end{equation}  
whereas for $ k_a \neq k_{a+1}$ they look as follows
\begin{equation}
  \raisebox{-.45\height}{\includegraphics[scale=0.7]{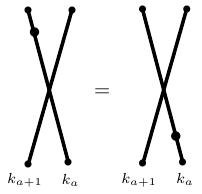}} \raisebox{-20\height}{,} \, \,  \, \,  \, \,  \, \,
   \raisebox{-.45\height}{\includegraphics[scale=0.7]{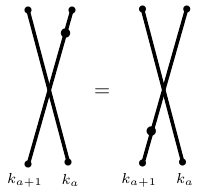}}
\end{equation}

Relation \eqref{eq Rnine} takes the following form
\begin{equation}\label{REdiaglast}
 \raisebox{-.45\height}{\includegraphics[scale=0.7]{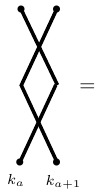}}   
           \left\{\begin{aligned}
  &\, \, \, \, \, \, 0  & &  \mbox{if } k_a = k_{a+1} \\
  & \raisebox{-.45\height}{\includegraphics[scale=0.7]{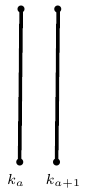}}    & &  \mbox{if } k_a \nrelbar k_{a+1} \\
  & \raisebox{-.45\height}{\includegraphics[scale=0.7]{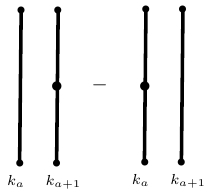}}     & &  \mbox{if } k_a \rightarrow k_{a+1} \\
  & \raisebox{-.45\height}{\includegraphics[scale=0.7]{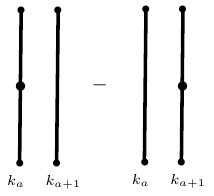}}     & &  \mbox{if } k_a \leftarrow k_{a+1} \\
\end{aligned} \right.
\end{equation}
and similarly \eqref{eq Rten}.

\medskip
For $ d= 1 $, we get that $ \YYdeqone $ is the usual Iwahori-Hecke of type $ A$ and
$ \RRdeqone $ is the quiver-Hecke algebra of the cyclic quiver $ \Gamma_e $.

\medskip
For a general $ d$, Rostam established in \cite{Ro} an isomorphism $ F: \YY \cong
\RR $. His isomorphism is an extension of the famous isomorphism
by Brundan, Kleshchev and Rouquier, see
\cite{brundan-klesc} and \cite{Rouq}, which corresponds to the $ d = 1 $ case of Rostam's
isomorphism.
Via $ F$ we may consider $ \YY $ as a $\Z$-graded
algebra.

\medskip
It follows directly from the construction of $ F $ in \cite{Ro} that 
\begin{equation}\label{second e(k)}
F(e(\kk)) = e(\kk)
\end{equation}
where $ e(\kk) \in \YY $ is as in \eqref{eigenspace}. 
Using $F$, for $ \ii \in \II $, $ \jj \in \JJn $, 
$ A  \in \SetPar_n $ and $ A^{\ord}  \in \OrdSetPar_n $
we define 
elements $ e(\ii)$, $ e(\jj)$, $ e(A) $, $\e(A)$, $ \e(\ii, A)$, 
$\e(A^{\ord})$, $ \e(\ii, A^{\ord}) \in \RR $
via 
\begin{equation}\label{defi e in RR}
\begin{aligned}  
    & e(\ii)  = F(e(\ii) ),&  & e(\jj)  = F(e(\jj) ),&  & e(A)  = F(e(A) ),   &  \e(A)  = F(\e(A) ), & \\
 &  \e(\ii, A)  = F(\e(\ii,  A)), & 
& \e(A^{\ord})  = F(\e( A^{\ord})), &   &  \e(\ii, A^{\ord})  = F(\e(\ii, A^{\ord}) )    
\end{aligned}    
\end{equation}
where $ e(\ii), e(\jj), e(A), \e(A), \e(\ii,A), \e(A^{\ord}),  \e(\ii, A^{\ord}) \in \YY $ are given in
\eqref{def e(i)}, \eqref{description1A}, \eqref{overline e(A)first},
\eqref{firstwethendefine}
and
\eqref{bytheaboveremarks}.
Note however that we don't have elements $ e(A^{\ord}) $ or $ e(\ii, A^{\ord}) $.

\section{The algebras $\RE$ and $ \REord $.}\label{qhecke RE}
Let $ \RR $ be the cyclotomic quiver Hecke algebra
introduced in Definition \ref{cycloHecke}. 
Let $ \REs $ be 
the subalgebra of $\RR $ generated by
\begin{equation} 
 \{ \e(\ii,A)  | \ii \in \II,  A \in \SetPar_n \}
 \cup \{ y_1, y_2, \ldots, y_n \} \cup \{ \psi_1, \psi_2, \ldots, \psi_{n-1} \}
 \end{equation}
where $ \e(\ii, A ) $ is given in \eqref{defi e in RR}.

\medskip 
The following Theorem gives some relations satisfied by the elements of $ \REs$. 
\begin{theorem}\phantomsection\label{teorem 6}
In $ \REs $ the following relations hold.  
\begin{align}
  \label{eq REcycloA} y_1 \e(\ii, A)  & = 0     &   &    \mbox{if }  i_1 = 0   \\
\label{eq REcycloB} \e(\ii,A)  & =  0    &   &   \mbox{if }  i_1 \neq 0   \\
  \label{eq REone}	\sum_{\ii \in \II} \sum_{A \in \SetPar_n} \e(\ii, A) 	& =  1   &   &   \\
  \label{eq REtwo}	\e(\ii, A) \e(\ii^{\prime}, A^{\prime})   & = \delta_{(\ii, A), (\ii^{\prime}, A^{\prime})} \e(\ii, A) &  &   \\
\label{eq REthree} y_a \e(\ii, A)  & = \e(\ii, A) y_a     &  &    \\
\label{eq REfour} \psi_a \e(\ii, A)  & = \e(\sigma_a( \ii, A)) \psi_a     &  &       \\
\label{eq REfive} y_a y_b   & = y_b y_a     &  &      \\
\label{eq REsix} \psi_a \psi_b   & = \psi_b \psi_a     &   &   \mbox{if }     | a-b | > 1     \\
\label{eq REseven} \psi_a y_{a+1} \e(\ii,A)    &   =  \left\{\begin{aligned} & (y_a \psi_a +1)\e(\ii, A) & &  \mbox{if }
i_a = i_{a+1} \mbox{ and } a\sim_A (a+1) \\
&  y_a \psi_a \e(\ii,A) & & \mbox{if } i_a \neq i_{a+1} \mbox{ or } a \not\sim_A (a+1)  \end{aligned} \right. & &  \\
\label{eq REeight} y_{a+1} \psi_a  \e(\ii, A)    &   =  \left\{\begin{aligned} & (\psi_a y_a  +1)\e(\ii, A) & &  \mbox{if }
i_a = i_{a+1} \mbox{ and } a \sim_A (a+1) \\
&  \psi_a y_a  \e(\ii, A) & & \mbox{if } i_a \neq i_{a+1} \mbox{ or } a \not\sim_A (a+1)   \end{aligned} \right. & &  \\
\label{eq REnine} \psi_a^2  \e(\ii, A)    &   =  \left\{\begin{aligned}
& 0 & &  \mbox{if } i_a = i_{a+1} \mbox{ and } a\sim_A (a+1) \\
& \e(\ii, A) & &  \mbox{if } i_a \nrelbar i_{a+1} \mbox{ or } a\not\sim_A (a+1) \\
& (y_{a+1} -y_a)\e(\ii, A) & &  \mbox{if } i_a \rightarrow i_{a+1} \mbox{ and } a\sim_A (a+1) \\
& (y_{a} -y_{a+1})\e(\ii, A) & &  \mbox{if } i_a \leftarrow i_{a+1} \mbox{ and } a\sim_A (a+1) \\
\end{aligned} \right. & & \\
\label{eq REten} \psi_{a+1} \psi_{a} \psi_{a+1}   \e(\ii, A)    &   =  \left\{\begin{aligned}
& (\psi_{a} \psi_{a+1} \psi_{a} -1) \e(\ii, A) & &  \mbox{if }
\left\{\begin{array}{l}i_{a+2} = i_{a} \rightarrow i_{a+1}   \mbox{ and }\\  a\sim_A (a+1) \sim_A( a+2) \end{array}
\right\} \\
& (\psi_{a} \psi_{a+1} \psi_{a} +1) \e(\ii, A) & &  \mbox{if }
\left\{\begin{array}{l}
i_{a+2} = i_{a} \leftarrow i_{a+1}  \mbox{ and }\\  a\sim_A (a+1)  \sim_A (a+2 )
 \end{array}
\right\}
\\
& \psi_{a} \psi_{a+1} \psi_{a} \e(\ii, A) & &  \mbox{otherwise.}  \\
\end{aligned} \right. & &   
\end{align}

\end{theorem}
\begin{dem}
  This is an immediate consequence of the relations for $ \RR $ given in Definition \ref{cycloHecke}
  together with the definition of $ \e(\ii, A) $ in \eqref{defi e in RR}.
  Note that the relations \eqref{eq RcycloA}--\eqref{eq Rnine}, defining $ \RR$, only depend on
  $j_a $ and $j_{a+1}$ through $ j_a = j_{a+1} $ or $ j_a \neq j_{a+1} $ and that for $ \e(\ii, A) $
  this corresponds exactly to $ a \sim_A b $ or $ a\not\sim_A b $. 
\end{dem}
  
\medskip
The Theorem motivates the following definition.
\begin{definition}\label{quiverE}
  \normalfont
Let $ \Kdomain $ be a field. Then $ \RE$ is the associative $\Kdomain$-algebra generated by elements 
\begin{equation} 
  \{ \e(\ii, A)  | \ii \in \II,   A \in \SetPar_n \}
\cup \{ y_1, y_2, \ldots, y_n \} \cup \{ \psi_1, \psi_2, \ldots, \psi_{n-1} \}
 \end{equation}
subject to the relations \eqref{eq REcycloA}--\eqref{eq REten}.

\medskip
For $ \ii \in \II $ and $ A \in \SetPar_n $ we define elements $ \e(\ii) \in \RE $
and $ \e(A) \in \RE $ via 
\begin{equation}\label{define e(A)}
  \e(\ii)= \sum_{ A \in \SetPar_n} \e(\ii, A), \, \, \, \,
  \e(A)= \sum_{ \ii \in \II} \e(\ii, A). 
\end{equation}
\end{definition}

Note that the definition of $ \RE$ neither depends on $ d $ nor on $ \xi$. Note also that the rules 
\begin{align}
  \label{degRE1 }  \deg y_a & = 2     &  &   \\
\label{degRE2 }  \deg \e(\ii, A) &= 0     &   &   \\
\label{degRE3 }  \deg \psi_a \e(\ii, A) &= 
\left\{\begin{aligned}
 0  &  & \mbox{if } & i_{a}  \nrelbar i_{a+1} \mbox{ or } a \not\sim_A (a+1) \\
  1  &  & \mbox{if } & (i_a \leftarrow i_{a+1} \mbox{ or } i_a \rightarrow i_{a+1}) \mbox{ and } a \sim_A (a+1)        \\
  -2  &  & \mbox{if } & i_a = i_{a+1} \mbox{ and } a \sim_A (a+1)        \\
\end{aligned} \right. & &  
\end{align}
define a degree function on $ \RE $, as one easily checks on the relations
\eqref{eq REcycloA}--\eqref{eq REten}. 

\medskip
By construction there exists a homomorphism $ \RE \rightarrow \REs$.
It is surjective but it is not clear if is also injective.

\medskip
There is a diagram calculus for $ \RE$ which is a blend of the diagram calculi for
$ \E $ and $ \RRdeqone$. To be more precise, an $ \RE$-diagram $ D $ is a KLR-diagram for
$ \RRdeqone$ joining 
northern points $ \{1, 2, \ldots, n \} $
with southern points $ \{1^{\prime}, 2^{\prime}, \ldots, n^{\prime} \} $
such that each pairs of points in $ \{1^{\prime}, 2^{\prime}, \ldots, n^{\prime} \} $
may be connected via a \lq tie\rq\ {\!}. Together these ties define a 
set partition denoted $ A_D \in \SetPar_n$.
Here are a couple of examples of $ \RE$-diagrams

\begin{equation}\label{719}
D=  \raisebox{-.45\height}{\includegraphics[scale=0.7]{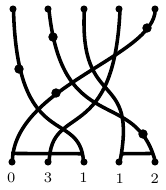}} \raisebox{-20\height}{,} \, \,  \, \,  \, \,  \, \,
E=  \raisebox{-.45\height}{\includegraphics[scale=0.7]{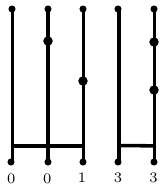}} 
  \end{equation}  
with corresponding 
set partitions $ A_D = \{ \{ 1,2,3 \} , \{ 4,5\} \} $ and also 
$ A_E = \{ \{ 1,2,3 \} , \{ 4,5\} \} $. Just as for $ \E $-diagrams, for $ D $ an $\RE$-diagram
we let $ A_D^{up} \in \SetPar_n $ be the upper set partition obtained from $ D $ by sliding all the ties to the northern points. The diagram calculus for $ \RE $ is the $ \Kdomain $-span of all $\RE$-diagrams
with the multiplication $ D D_1 $ of diagrams $ D $ and $D_1 $
given by vertical concatenation with $ D $ on top of $ D_1 $
where $ D D_1 $ is set equal to $ 0 $ if the $\KKn$-decorations on $ D $ and $ D_1 $ do not match
or if $ A_D \neq  A_{D_1}^{up}$.
The elements $ \e(\ii,A),  y_a e(\ii, A)$ and $ \psi_a e(\ii,A) $ generate $ \RR $, and correspond to the following diagrams 
\begin{equation}\label{720}
  \e(\ii,A) \mapsto  \raisebox{-.5\height}{\includegraphics[scale=0.7]{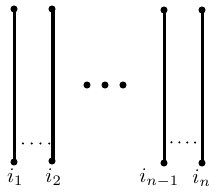}} \raisebox{-20\height}{,}
  \, \,  \, \,   
y_a e(\ii, A) \mapsto   \raisebox{-.5\height}{\includegraphics[scale=0.7]{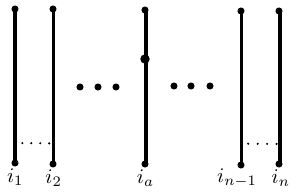}} \raisebox{-20\height}{,}  \, \,  \, \, 
\psi_a e(\ii,A) \mapsto   \raisebox{-.5\height}{\includegraphics[scale=0.7]{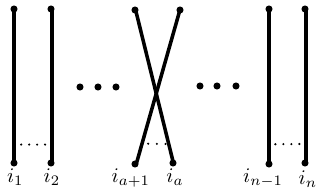}} 
\end{equation}
where for each of
the $ \RE$-diagrams in \eqref{720} we have used dashed lines to indicate that
the associated set partition is $A$.
Using this, the relations
\eqref{eq REcycloA}--\eqref{eq REten}
can now be described diagrammatically.

\medskip
Indeed, 
for $ i_a =i_{a+1}$ and $ a \sim_A (a+1)$
the relations 
\eqref{eq REseven} and \eqref{eq REeight} look as follows
\begin{equation}\label{REdiag1A}
  \raisebox{-.45\height}{\includegraphics[scale=0.7]{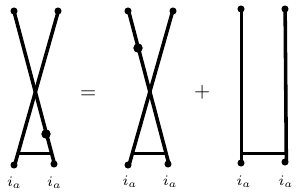}} \raisebox{-20\height}{,} \, \,  \, \,  \, \,  \, \,
   \raisebox{-.45\height}{\includegraphics[scale=0.7]{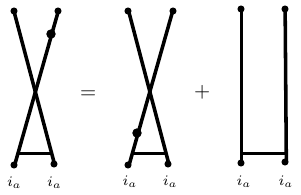}}
\end{equation}  
whereas for $ i_a \neq i_{a+1}$ or $a \not\sim_A (a+1) $ they look as follows
\begin{equation}
  \raisebox{-.45\height}{\includegraphics[scale=0.7]{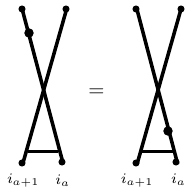}} \raisebox{-20\height}{,} \, \,  \, \,  \, \,  \, \,
  \raisebox{-.45\height}{\includegraphics[scale=0.7]{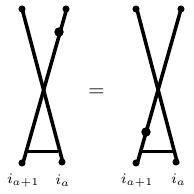}} \raisebox{-20\height}{,} \, \,  \, \,  \, \,  \, \,
  \raisebox{-.45\height}{\includegraphics[scale=0.7]{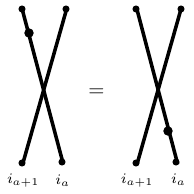}} \raisebox{-20\height}{,} \, \,  \, \,  \, \,  \, \,
     \raisebox{-.45\height}{\includegraphics[scale=0.7]{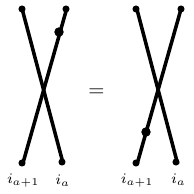}} \raisebox{-20\height}{.}
\end{equation} 

For \eqref{eq REnine} we have the following expressions
\begin{equation}
 \raisebox{-.45\height}{\includegraphics[scale=0.7]{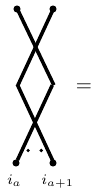}}   
 \left\{\begin{aligned}
& \, \, 0     & &  \mbox{if } i_a = i_{a+1} \mbox{ and } a \sim_A (a+1)\\
            & \raisebox{-.45\height}{\includegraphics[scale=0.7]{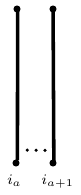}}    & &  \mbox{if }
i_a  \nrelbar i_{a+1} \mbox{ or } a \not\sim_A (a+1) \\
& \raisebox{-.45\height}{\includegraphics[scale=0.7]{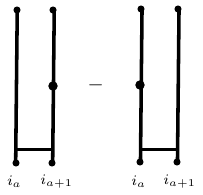}}     & &  \mbox{if }
  i_a \rightarrow i_{a+1} \mbox{ and } a\sim_A (a+1)  \\
& \raisebox{-.45\height}{\includegraphics[scale=0.7]{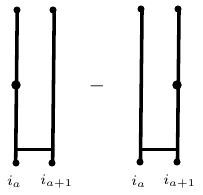}}     & &  \mbox{if }
  i_a \leftarrow i_{a+1} \mbox{ and } a\sim_A (a+1).  \\
\end{aligned} \right.
\end{equation}

Similarly, for \eqref{eq REten} we have the following expressions
\begin{equation}\label{REdiag1Alast}
 \raisebox{-.45\height}{\includegraphics[scale=0.7]{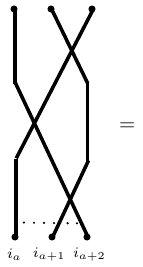}}   
           \left\{\begin{aligned}
           & \raisebox{-.45\height}{\includegraphics[scale=0.7]{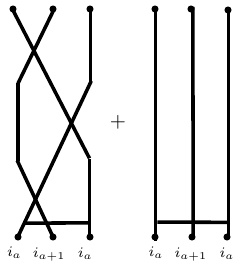}}    & &  \mbox{if }
i_a = i_{a+2} \mbox{ and }            i_a \rightarrow i_{a+1} \mbox{ and } a\sim_A (a+1)  \sim_A (a+2 ) \\
& \raisebox{-.45\height}{\includegraphics[scale=0.7]{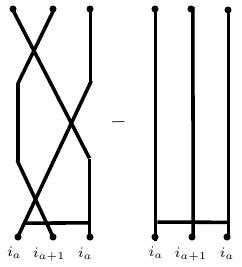}}     & &  \mbox{if }
i_a = i_{a+2} \mbox{ and }  i_a \leftarrow i_{a+1} \mbox{ and } a\sim_A (a+1)  \sim_A (a+2 )  \\
& \raisebox{-.45\height}{\includegraphics[scale=0.7]{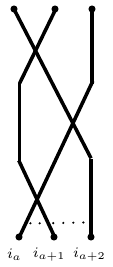}}     & &  \mbox{otherwise.}
\end{aligned} \right.
\end{equation}

\medskip
There is yet another natural subalgebra $  \REords $ 
of $\RR $ generated by

\begin{equation}\label{generators} 
 \{ \e(\ii,A^{\ord})  | \ii \in \II,  A \in \OrdSetPar_n \}
 \cup \{ y_1, y_2, \ldots, y_n \} \cup \{ \psi_1, \psi_2, \ldots, \psi_{n-1} \}
 \end{equation}
where $ \e(\ii, A^{\ord} ) $ is as in \eqref{defi e in RR}.
The elements in \eqref{generators} satisfy relations corresponding to
\eqref{eq REcycloA}--\eqref{eq REten}, with $ A $ replaced by $ A^{\ord} $, 
and so we arrive at the following definition.

\begin{definition}\label{quiverEord}
  \normalfont
Let $ \Kdomain $ be a field. Then $ \REord$ is the associative $\Kdomain$-algebra generated by elements 
\begin{equation} 
  \{ \e(\ii, A^{\ord})  | \ii \in \II,   A^{\ord} \in \OrdSetPar_n \}
\cup \{ y_1, y_2, \ldots, y_n \} \cup \{ \psi_1, \psi_2, \ldots, \psi_{n-1} \}
 \end{equation}
subject to the relations \eqref{eq REcycloA}--\eqref{eq REten}, but with $ A $ replaced by $ A^{\ord}$
throughout.

\medskip
Similarly to what we did for $ \RE$, 
for $ \ii \in \II $ and $ A^{\ord} \in \OrdSetPar_n $ we define elements $ \e(\ii) \in \REord $
and $ \e(A^{\ord}) \in \REord $ via
\begin{equation}\label{eordA}
  \e(\ii)= \sum_{ A^{\ord} \in  \OrdSetPar_n} \e(\ii, A^{\ord}), \, \, \, \,
  \e(A^{\ord})= \sum_{ \ii \in \II} \e(\ii, A^{\ord}).
\end{equation}
\end{definition}

Once again, one checks that $ \REord$ is a graded algebra, with degree function given by the same 
formulas as the ones 
for $ \RE$. There is also a diagram calculus associated with $ \REord $, similar to the one
for $ \RE$, except that this time we associate
with a diagram $D$ for $  \REord $
an ordered set partition $ A_D^{ord} \in \OrdSetPar_n$.

\section{The isomorphisms $ \Eord \cong \REord$ and $ \E \cong \RE$.}\label{isomorphismProof}
In this section we establish our main results, the isomorphisms
$ \E \cong \RE$ and $ \Eord \cong \REord$. 
For the first isomorphism we need to impose the condition that $ q $ admit a square root
$ q^{\frac{1}{2} } $ in $ \Kdomain$, whereas the second
isomorphism holds for any $q\in  \Kdomain$. Our proofs have much in common with Rostam's
proof of $ \YY \cong \RR $ in \cite{Ro}, but there are also several differences throughout.

\medskip
Let $ \E $ and $ \Eord $ be the bt-algebras from Definition \ref{btdefinition} and Definition \ref{orderedBT}. 
We then define elements $  X_1, X_2 , \ldots, X_n  $ in $ \E $ and $ \Eord$ via 
$ X_1 = 1 $ and recursively $ X_{a+1} = q^{-1} g_a X_{a} g_a$. 
These formulas are the same as the formulas for $ X_a $ in \eqref{JM for YY}
and so we have $ \epsilon(X_a ) = X_a $ where $ \iota: \E \rightarrow \YY $ and $ \epsilon: \Eord \rightarrow \YY $ are
the homomorphisms from Theorem \ref{orderedBTinclusion}. 
We set 
\begin{equation}
{\cal L}_n = \{X_1, \ldots, X_n \} 
\end{equation}
and shall refer to $ {\cal L}_n $ as the set of {\it Jucys-Murphy elements} for $ \E$ and $ \Eord$, respectively. 
However, it is important to note that although $ {\cal L}_n $ is a set of $ \JM $-elements
for $ \E $ with respect to the cellular structure introduced in \cite{ERH}, the corresponding content functions
for $\Efractionfield $ 
do not satisfy the separation condition 
of \cite{Mat-So} and therefore the arguments
from section \ref{JM for YYsection} do not carry over to $ \E$, and also not
to $ \Eord$.

\medskip
To remedy this problem, at least partially, we proceed as follows. 
Choose $ d \ge n  $ such that $ d $ is is coprime to $  {\rm char } \,  \Kdomain$. 
Choose moreover a field extension $  \Kdomain^{\prime} $ of $ \Kdomain$ such that $ \Kdomain^{\prime} $ contains
a primitive $ d$'th root of unity $ \xi$. Let $ \Efieldext $ be the bt-algebra defined over $  \Kdomain^{\prime} $
and let $\YYfieldext$ be the Yokonuma-Hecke algebra defined over $  \Kdomain^{\prime} $. 
Then, by Theorem 14 in \cite{ERH}, we have $ \Efieldext \subseteq \YYfieldext $ via the map
given in \eqref{mimicking}
and moreover
$ \E \subseteq \Efieldext$
since the basis for $ \E $ is independent of the ground field. The same arguments also hold for
$ \Eord$ and so we obtain inclusions
of $ \Kdomain$-algebras
\begin{equation}\label{fieldext}
\Efield \subseteq \YYfieldext \, \, \, \mbox{and } \, \, \,  \Eordfield \subseteq \YYfieldext.
\end{equation}

\medskip

Noting that $  {\rm char }_q (\Kdomain)  = {\rm char }_q (\Kdomain^{\prime}) $ we let 
$ \II $, $ \JJn $ and $ \KKn $ be the residue sequences for $\YYfieldext$, as in section
\ref{idempotents in YY}. 

\medskip
By general principles of the Jordan normal form, for $ \ii \in \II $, the idempotent $ e(\ii ) \in \YYfieldext$ introduced in \eqref{eigenspaceIJ}
is a polynomial expression in the $ X_a$'s and therefore $ e(\ii) \in \E \subseteq \Eord$,
$ \e(\ii, A) \subseteq \E$
and $ \e(\ii, A^{\ord}) \subseteq \Eord$
where
\begin{equation}\label{generalprinciples}
\e(\ii, A) = e(\ii ) \e(A)  \, \, \, \mbox{for } \ii \in \II, A \in \SetPar_n  \,  \mbox{ and } \, \,
  \e(\ii, A^{\ord}) = e(\ii ) \e(A^{\ord})  \, \, \, \mbox{for } \ii \in \II, A \in \OrdSetPar_n  
\end{equation}
see \eqref{bytheaboveremarksA}. Once again, although our notation coincides with the notation for 
$ \e(\ii, A)  \in \RE $, see Definition \ref{quiverE}, and $ \e(\ii, A^{\ord}) \in \REord $,
see Definition \ref{quiverEord}, this should not give rise to confusion. 

\medskip
For $ 1\le a \le n $ we define $ y_a \in \E \subseteq \Eord$ via
\begin{equation}\label{defy_a}
y_a = \sum_{\substack{  \ii \in \II }} (1-q^{-i_a} X_a )e(\ii)
\end{equation}
once again using the same notation as for $ y_a \in \RE$ and for $ y_a \in \REord$.

\medskip
For $ \ii \in \II $ we introduce 
$ y_a(\ii) \in  \Kdomain[ y_a ] $ via
\begin{equation}\label{86A}
y_a(\ii) = q^{i_a}(1-y_a)
\end{equation}  
which gives us the following formula
\begin{equation}\label{86}
X_a =  \sum_{\ii \in \II} X_a e(\ii) = \sum_{\ii \in \II} y_a(\ii) e(\ii) 
\end{equation}
valid both in $ \E $ and in $ \Eord$. 
It follows from the definition of the $ e(\ii)$'s in \eqref{eigenspaceIJ} 
that the $ y_a$'s are nilpotent, as elements of $\YYfieldext $ and therefore also
as elements of $ \E $ and $ \Eord$. 

\medskip
On the other hand, repeating the arguments of 
Lemma 2.1 of \cite{brundan-klesc}, where it is shown that the $ y_a$'s are nilpotent in $ \RRdeqone$,
we get that the $ y_a$'s are also nilpotent as elements of $ \RE $ and $ \REord$.

\medskip
For $ 1 \le a < n $ we now introduce the element
$  \displaystyle P_a = \sum_{  \substack{  \ii \in \II \\ A \in \SetPar_n}}     P_a(\ii,A) \e(\ii, A)
\in \E
$
where
$ P_a(\ii,A) $ is given by 
\begin{equation}\label{88}
 P_a(\ii, A)   = \left\{ \begin{aligned}
&   \left\{ \begin{aligned}
&     1 && \mbox{if } i_{a} = i_{a+1} \\
     & (1-q) \left(  1 -y_{a}(\ii) y_{a+1}(\ii)^{-1} \right)^{-1}               && \mbox{if } i_{a} \neq i_{a+1} \\ 
  \end{aligned} 
\right\}  && \mbox{ if }   a \sim_A (a+1)                 \\
& 0 && \mbox{ if }   a \not\sim_A (a+1) .
 \end{aligned}      
 \right. 
\end{equation}
Note that the expressions on the right hand side of \eqref{88} make sense because of \eqref{86}
and the nilpotency of
the $ y_a $'s.
In fact we have
\begin{equation}\label{concreteformula}
 P_a(\ii, A) \e(\ii, A)  = \left\{ \begin{aligned}
&   \left\{ \begin{aligned}
&      \e(\ii, A) && \mbox{if } i_{a} = i_{a+1} \\
     & (1-q) \left(  1 -X_{a} X_{a+1}^{-1} \right)^{-1}   \e(\ii, A)             && \mbox{if } i_{a} \neq i_{a+1} \\ 
  \end{aligned} 
\right\}  && \mbox{ if }   a \sim_A (a+1)                 \\
& 0 && \mbox{ if }   a \not\sim_A (a+1) 
 \end{aligned}      
 \right. 
\end{equation}
see \eqref{inmoredetail} for the concrete formula.

\medskip
We also 
define elements $ P_a $ of $ \Eord $, $ \RE$ and $ \REord$, using formulas corresponding to
the ones 
for $ P_a $ in $\E $, but with $ A $ replaced by $ A^{\ord} $, etc.

\medskip
Recall that in our treatment of $ \E $ we are assuming that $ q $ has a square root
$ q^{\frac{1}{2}} $ in $  \Kdomain $. 
For $1 \le a < n$, we now introduce the 
element $ \displaystyle Q_a =
\sum_{  \substack{  \ii \in \II \\ A \in \SetPar_n}}     Q_a(\ii,A) \e(\ii, A) \in \E$ where
\begin{equation}\label{87}
 Q_a(\ii, A)  = \left\{ \begin{aligned}
&   \left\{ \begin{aligned}
&     1-q- q y_{a+1} -y_a && \mbox{if } i_a = i_{a+1} \\
     & ( y_a(\ii) -q y_{a+1}(\ii) )/( y_a(\ii) - y_{a+1}(\ii))               && \mbox{if } i_a \nrelbar i_{a+1} \\
     & ( y_a(\ii) -q y_{a+1}(\ii) )/( y_a(\ii) - y_{a+1}(\ii))^2               && \mbox{if } i_a \rightarrow i_{a+1}  \\
     & q^{i_a}               && \mbox{if } i_a \leftarrow i_{a+1}  
\end{aligned} 
\right\}  && \mbox{ if } a \sim_A (a+1)  \\
& f_{a, A}  && \mbox{ if }  a \not\sim_A (a+1)
 \end{aligned}      
  \right.
\end{equation}
where $ f_{a, A} $ is the constant function
\begin{equation}\label{bothcases1}
f_{a, A} = q^{\frac{1}{2}}.
\end{equation}

\medskip
We similarly define $ \displaystyle Q_a
= \sum_{  \substack{  \ii \in \II \\ A^{\ord} \in \OrdSetPar_n}}     Q_a(\ii,A^{\ord}) \e(\ii, A^{\ord})
\in \Eord $ using expressions analogous to the ones
in \eqref{87}, except that 
$ f_{a, A^{\ord}} \in \{1, q \} $ is now defined as follows.
Let $ A_i $ be the block of $ A^{\ord} $ containing $ a $ and let $ A_j $ be the block of
$ A^{\ord} $ containing $ a +1$. Then we set
\begin{equation}\label{bothcases2}
f_{a, A^{\ord}} = \left\{ \begin{array}{ll} q & \mbox{if } A_i < A_j  \\
1 & \mbox{otherwise. }   \end{array}\right.
\end{equation}
Note that in both cases \eqref{bothcases1} and \eqref{bothcases2} the function $ f_{a, A} $
has the following property
\begin{equation}\label{wehaveinbothcases}
f_{a, A} f_{ a, \sigma_a(A)} = q  \, \mbox{ when } a \not\sim_A (a+1) 
\end{equation}
since $ A_i < A_j \Longleftrightarrow \sigma_a(A_i) < \sigma_a(A_j) $ by
definition of the $ \Si_n$-action on $ \OrdSetPar_n$.

\medskip
Once again, due to \eqref{86}, the expressions on the right hand side of \eqref{87},
defining $ Q_a $,
are (finite) power 
series 
and $ Q_a $ is invertible. The formulas in \eqref{87} also define
elements of $ \RE$ and $ \REord$ that we shall denote the same way $ Q_a$.

\medskip
For $ 1 \le a  < n $ we use the $ P_a $'s and the $ Q_a $'s to define elements $ g_a $ of $ \RE $ via
\begin{equation}\label{RRgen}
  g_a = \sum_{ \ii \in \II, A \in \SetPar_n} \left( \psi_a Q_a(\ii, A)- P_a(\ii, A) \right)\e(\ii,A) . 
\end{equation}
Similarly we have elements $ g_a \in  \REord $, given by the corresponding formula.

\medskip
Recall the intertwining elements $ \F_a $ for
$  \YYresiduefield $ introduced in
\eqref{intertwinedef}. Using \eqref{intertwineE}, they may also be viewed as elements
of $ \E  $ or of $ \Eord$. Note that by \eqref{concreteformula} we have
\begin{equation}\label{intertwineP}
F_a = g_a + P_a. 
\end{equation}

Using the $ \F_a $'s we can now finally introduce elements $ \psi_a \in \E  $ as follows
\begin{equation}\label{analogouspsi}
  \psi_a=
\sum_{ \ii \in \II, A \in \SetPar_n}  \Phi_a Q_a^{-1}      \e(\ii,A) 
 \end{equation}
and we also introduce $ \psi_a \in \Eord  $ using the same formula 
\eqref{analogouspsi}, but with $ A^{\ord} $ instead of $ A $.

\medskip
The following is the main Theorem of the present article. 
\begin{theorem}\phantomsection\label{twomaintheoremsX}
\begin{description}
\item[a)]   The homomorphism $ F: \RE \rightarrow \E $ given by $ F(\psi_a ) = \psi_a $
  for $ 1 \le a < n $ and $ F(\e(\ii, A ) ) = \e(\ii, A) $ for $ \ii \in \II $ and $ A \in \SetPar_n $
  and $ F(y_a) = y_a $ for $ 1 \le a \le n $
is well-defined and induces an isomorphism $ \E \cong \RE $. 
  \item[b)]   The homomorphism $ F: \REord \rightarrow \Eord $ given by $ F(\psi_a ) = \psi_a $
for $ 1 \le a < n $ and $ F(\e(\ii, A^{\ord} ) ) = \e(\ii, A^{\ord}) $ for $ \ii \in \II $ and $ A^{\ord} \in \OrdSetPar_n $ 
  and $ F(y_a) = y_a $ for $ 1 \le a \le n $
is well-defined and induces an isomorphism $ \Eord \cong \REord $. 
  \end{description}
\end{theorem}  

\medskip
Before embarking on the proof we need the following preparatory Lemmas.
There is a left $ \Si_n$-action on invertible power series $ Q(y_1,  \ldots, y_n) \in \
\Kdomain[[ y_1,  \ldots, y_n ]]^{\times} $ given on simple reflections by 
\begin{equation}
\prescript{\sigma_a \!}{}Q(y_1,  \ldots, y_a, y_{a+1}, \ldots,y_n) = 
Q(y_1,  \ldots, y_{a+1}, y_{a}, \ldots,y_n).  
\end{equation}
Recall the left $ \Si_n$-action on pairs $ (\ii, A ) \in \II \times \SetPar_n $ given by
$ \sigma  (\ii, A) = (\sigma \ii, \sigma A)$, and similarly
the left $ \Si_n$-action on pairs $ (\ii, A^{\ord} ) \in \II \times \OrdSetPar_n $ given by
$ \sigma  (\ii, A^{\ord}) = (\sigma \ii, \sigma A^{\ord})$.
With these actions we have the following Lemma involving $P_a $.

\begin{lemma}\label{lemma11newPPP}
  \phantomsection
  The following \lq braid\rq\ relations hold: 
  \begin{multicols}{2}
    \begin{description}
\item[a)]     $  \prescript{\sigma_{a+1} \sigma_{a}  \!}{}
    P_{a+1}(\ii, A )  = P_{a} (\sigma_{ a+1} \sigma_{a}  (\ii ,A)) $ 
  \item[b)]  $  \prescript{\sigma_{a} \sigma_{a+1}  \!}{}
    P_{a}(\ii, A )  = P_{a+1} (\sigma_{ a} \sigma_{a+1}  (\ii ,A)) $
\item[c)]     $  \prescript{ \sigma_{a}  \!}{}
    P_{a+1}(\sigma_{a}(\ii, A) )  = \prescript{\sigma_{a+1}   \!}{}P_{a} (\sigma_{ a+1}  (\ii ,A)) $ 
  \item[d)] $  \prescript{ \sigma_{a+1}  \!}{}
    P_{a}(\sigma_{a+1}(\ii, A) )  = \prescript{\sigma_{a}   \!}{}P_{a+1} (\sigma_{ a}  (\ii ,A)) .$ 
    \end{description}
  \end{multicols}
There are similar formulas with $ A \in \SetPar_n $ replaced by $ A^{\ord} \in \OrdSetPar_n$.
\end{lemma}
\begin{dem}
The fours parts of the Lemma are formally
equivalent, hence it is enough to show $ {\bf a)} $.
\medskip
If $ (a+1) \not\sim_A (a+2) $ then $P_{a+1}(\ii, A) = 0 $ and therefore also
$\prescript{\sigma_{a+1} \sigma_{a}  \!}{}P_{a+1}(\ii, A) = 0 $. But $ (a+1) \not\sim_A (a+2) $
implies $ a \not\sim_{\sigma_{a+1} \sigma_{a}  A} (a+1) $ and so also
$P_{a}(\sigma_{ a+1} \sigma_{a}  (\ii ,A))  = 0 $, as claimed.

\medskip
If $ (a+1) \sim_A (a+2) $ and $ i_{a+1} \neq i_{a+2} $ we have $ P_{a+1}(\ii, A) = (1-q)y_{a+1}(\ii)/(y_{a+1}(\ii)-y_{a}(\ii)) $
and therefore
\begin{equation}
  \prescript{\sigma_{a+1} \sigma_{a}  \!}{}P_{a+1}(\ii, A) = (1-q)\left( q^{i_{a+1}}(1-y_{a})    \right) /\left(   q^{i_{a+1}}(1-y_{a}) -q^{i_{a+2}}
  (1-y_{a+2})               \right). 
\end{equation}
But $ (a+1) \sim_A (a+2) $
implies $ a \sim_{\sigma_{a+1} \sigma_{a}  A} (a+1) $ and since
$ (\sigma_{a+1} \sigma_{a} \ii)_a  =i_{a+1} \neq  (\sigma_{a+1} \sigma_{a} \ii)_{a+1}= i_{a+2} $ 
we also get
\begin{equation}
  P_{a}(\sigma_{ a+1} \sigma_{a}  (\ii ,A)) = 
(1-q)\left( q^{i_{a+1}}(1-y_{a})    \right) /\left(   q^{i_{a+1}}(1-y_{a}) -q^{i_{a+2}}
  (1-y_{a+2})               \right)
\end{equation}
as claimed. The remaining case $ (a+1) \sim_A (a+2) $ and $ i_{a+1} = i_{a+2} $ is shown the same way. 
\end{dem}

\medskip \medskip
The next Lemma states 
that $Q_a $ has properties similar to the ones for $ P_a $ given in Lemma 
\ref{lemma11newPPP}. 
\begin{lemma}\label{lemma11newPP}
  \phantomsection
  The following \lq braid\rq\ relations hold: 
  \begin{multicols}{2}
  \begin{description}
  \item[a)]     $  \prescript{\sigma_{a+1} \sigma_{a}  \!}{}
    Q_{a+1}(\ii, A )  = Q_{a} (\sigma_{ a+1} \sigma_{a}  (\ii ,A)) $    
  \item[b)]  $  \prescript{\sigma_{a} \sigma_{a+1}  \!}{}
    Q_{a}(\ii, A )  = Q_{a+1} (\sigma_{ a} \sigma_{a+1}  (\ii ,A)) $
\item[c)]     $  \prescript{ \sigma_{a}  \!}{}
    Q_{a+1}(\sigma_{a}(\ii, A) )  = \prescript{\sigma_{a+1}   \!}{}Q_{a} (\sigma_{ a+1}  (\ii ,A)) $ 
  \item[d)] $  \prescript{ \sigma_{a+1}  \!}{}
    Q_{a}(\sigma_{a+1}(\ii, A) )  = \prescript{\sigma_{a}   \!}{}Q_{a+1} (\sigma_{ a}  (\ii ,A)) .$ 
\end{description}
  \end{multicols}
There are similar formulas with $ A \in \SetPar_n $ replaced by $ A^{\ord} \in \OrdSetPar_n$.  
  \end{lemma}
  \begin{dem}
Once again, the four parts of the Lemma are formally equivalent, and so it is enough to show $ {\bf a)} $.
We first consider $ A \in \SetPar_n$. 
\medskip

Suppose that $ (a+1) \sim_A (a+2) $.

\medskip
  If $ i_{a+1} = i_{a+2} $ then
  $ Q_{a+1}(\ii, A) = 1-q - q y_{a+2} -y_{a+1} $ and so 
  \begin{equation}
\prescript{\sigma_{a+1} \sigma_{a}   \!}{}
    Q_{a+1}(\ii, A)=  1-q - q y_{a+1} -y_a. \end{equation}
  On the other hand,  $ \sigma_{a+1} \sigma_{a} \ii = (i_1, \ldots, i_{a+1}, i_{a+2}, i_{a}, \ldots, i_n) $
  and $ a \sim_{\sigma_{a+1} \sigma_{a}  A} (a+1) $ and so also 
  $ Q_{a} (\sigma_{ a+1} \sigma_{a} (\ii ,A))  = 1-q - q y_{a+1} -y_{a} $ 
  which shows $ {\bf a)} $ in this case.

  \medskip
If $ i_{a+1} \nrelbar i_{a+2} $ then 
$ Q_{a+1}(\ii, A) =
( y_{a+1}(\ii) -q y_{a+2}(\ii) )/( y_{a+1}(\ii) - y_{a+2}(\ii)) $
and so
\begin{equation}
\prescript{\sigma_{a+1} \sigma_{a}   \!}{}Q_{a+1}(\ii, A)=  \left( q^{i_{a+1}}(1-y_a ) -q^{i_{a+2}+1}(1-y_{a+1}) \right)/
 \left( q^{i_{a+1}}(1-y_a ) -q^{i_{a+2}}(1-y_{a+1}) \right).
\end{equation}
On the other hand,  we still have 
$ \sigma_{a+1} \sigma_{a} \ii = (i_1, \ldots, i_{a+1}, i_{a+2}, i_{a}, \ldots, i_n) $
  and $ a \sim_{\sigma_{a+1} \sigma_{a} A} (a+1) $ and so also
  \begin{equation}
Q_{a} (\sigma_{ a+1} \sigma_{a}  (\ii ,A))  =
\left( q^{i_{a+1}}(1-y_a ) -q^{i_{a+2}+1}(1-y_{a+1}) \right)/
 \left( q^{i_{a+1}}(1-y_a ) -q^{i_{a+2}}(1-y_{a+1}) \right)
  \end{equation}
  which shows $ {\bf a)} $ in this case. The
  cases $ i_{a+1} \rightarrow i_{a+2} $ and $ i_{a+1} \leftarrow i_{a+2} $ are proved the same way.

  \medskip
  We are only left with the case where $ (a+1) \not\sim_A (a+2) $. Here we have
  $  Q_{a+1}(\ii, A) =q^{ \frac{1}{2}} $ and so we also have
  $ \prescript{\sigma_{a+1} \sigma_{a}   \!}{}Q_{a+1}(\ii, A)  = q^{ \frac{1}{2}} $.
  But $ (a+1) \not\sim_A (a+2) $ implies $ a \not\sim_{\sigma_{a+1} \sigma_{a}  A} (a+1) $ and
  so $ Q_{a} (\sigma_{ a+1} \sigma_{a}  (\ii ,A))  = q^{\frac{1}{2}} $ and therefore
  $  Q_{a} (\sigma_{ a+1} \sigma_{a}  (\ii ,A)) = q^{\frac{1}{2}}$. This shows the
  Lemma in the last
  case.

  \medskip
  We next consider $ A^{\ord} \in \OrdSetPar_n$. Here we can repeat the arguments 
  from the $ A \in \SetPar_n $-situation word by word, except when $ (a+1) \not\sim_A (a+2) $, 
  and we therefore only consider the case $ (a+1) \not\sim_A (a+2) $. 

  \medskip
  Let $ A_i $ be the block of $ A $ containing $ a+1 $ and let $ A_j $ be the block of $ A $ containing $ a+2$. 
  If $ A_i < A_j $ we have
  $  Q_{a+1}(\ii, A) =q $ and so also
  $ \prescript{\sigma_{a+1} \sigma_{a}   \!}{}Q_{a+1}(\ii, A)  = q$.
But $ A_i < A_j  \Longleftrightarrow \sigma_{a+1} \sigma_{a}(A_i) < \sigma_{a+1} \sigma_{a}(A_j)   $
whereas $ a \in \sigma_{a+1} \sigma_{a}(A_i)  $ and $ a+1 \in \sigma_{a+1} \sigma_{a}(A_j) $ 
and 
so also $Q_{a} (\sigma_{ a} \sigma_{a+1}  (\ii ,A)) = q $, as claimed. 
If $ A_i > A_j $ we argue the same way and this completes the proof of the Lemma.
  \end{dem}

\medskip
In the Theorems \ref{commutation intertwiner} and 
\ref{commutation intertwiner setpar}
we have established commutation formulas involving $ \F_a$'s.
The next Lemma gives a few more such commutation formulas.

\begin{lemma}
  \phantomsection\label{lemma41}
  For $ \ii \in \II ,$ 
the $ \Phi_a\! $'s verify the following formulas.
  \begin{description}
  \item[a)] $ \Phi_a y_b = y_b \Phi_a \, \, \mbox{ if } b\neq a, a+1$
  \item[b)] $ \Phi_a Q_b  = Q_b \Phi_a \, \, \mbox{ if } | a-b | >1$
  \end{description}
\end{lemma}
\begin{dem}
The statement $ {\bf a)} $ is immediate from \eqref{intertwine two} in 
Theorem \ref{commutation intertwiner} and \eqref{intertwine two setpar} in Theorem 
\ref{commutation intertwiner setpar}, together with the definition of
$ y_a $ in \eqref{defy_a}. 
From $ {\bf a)} $ we get $ {\bf b)} $ since
$Q_b $ is a polynomial expression in $ y_b $ and $ y_{b+1} $ and
$ \{b, b+1 \} \cap  \{a, a+1 \} = \emptyset $. 
\end{dem}

\medskip
With these preparations we are now in position to prove Theorem \ref{twomaintheoremsX}. 

\medskip
\begin{dem}
  (Theorem \ref{twomaintheoremsX}).
  We start by showing $ {\bf a)} $. 
\medskip
  
We first check that $ F: \RE \rightarrow \E $ is well-defined, in other words that
  the elements 
  $  \psi_a $, $\e(\ii, A) $ and $y_a $ of $ \E $ satisfy the defining relations
\eqref{eq REcycloA}--\eqref{eq REten} for $ \RE$. 

\medskip
\noindent
The relations \eqref{eq REcycloA}--\eqref{eq REthree}
follow immediately from \eqref{bytheaboveremarks} and
\eqref{defy_a} and the definitions. 
The relation \eqref{eq REfour} follows from the definition of $ \psi_a $ in
\eqref{analogouspsi} and $ {\bf a)} $ of Lemma \ref{lemma41}
and the relation \eqref{eq REfive} is clear from \eqref{defy_a} and the commutativity
of the $X_a$'s. Similarly, the relation \eqref{eq REsix} follows from 
$ {\bf d)} $ and $ {\bf e)} $ 
of Lemma \ref{lemma41} together with the definition of $ \psi_a $ in \eqref{analogouspsi}. 

\medskip
For the proof of the relations \eqref{eq REseven}-\eqref{eq REnine} we first make
the following remark in analogy with Remark 2.7 of \cite{Ro}. Suppose that 
$ a \sim_A (a+1) $. Then, in view of Theorem \ref{commutation intertwiner setpar},
Lemma \ref{lemma11newPP}, Lemma \ref{lemma11newPPP}
and the definitions, the 
$ \e(\ii,A) $'s, $ y_a$'s and $ \psi_a$'s are built up of 
ingredients that verify exactly the same commutation relations as the corresponding ingredients
in the Hecke algebra $ \Hecken =\YYdeqone $ and so
the calculations of \cite{brundan-klesc} can be repeated word by word to show that 
the relations \eqref{eq REseven}-\eqref{eq REnine} hold in this case.

\medskip
We are therefore only left with showing \eqref{eq REseven}-\eqref{eq REnine}
in the case $ a \not\sim_A (a+1) $. Here \eqref{eq REseven} is the
equality $ \psi_a y_{a+1} \e(\ii, A) =  y_{a}  \psi_a \e(\ii, A) $. 
But $ y_{a}  \psi_a \e(\ii, A) = y_{a}   \e( \sigma_a(\ii, A)) \psi_a $ 
and so \eqref{eq REseven} is equivalent to the identity 
\begin{equation} \F_a q^{ -\frac{1}{2}} ( 1- q^{-i_{a+1} } X_{a+1}) \e(\ii, A) =
 ( 1- q^{-i_{a+1} } X_{a})  \F_a q^{ -\frac{1}{2}}\e(\ii, A) 
\end{equation}
which follows from \eqref{intertwine four setpar} of 
Theorem \ref{commutation intertwiner setpar}. The relation \eqref{eq REeight} is proved the
same way.

\medskip
When 
$a \not\sim_A (a+1) $
the relation \eqref{eq REnine} is the formula 
$ \psi_a^2 \e(\ii, A) = \e(\ii, A) $, or equivalently
\begin{equation} \F_a q^{ -\frac{1}{2}} \F_a q^{ -\frac{1}{2}} \e(\ii, A)
  = \e(\ii, A)
\end{equation}
which follows directly from
\eqref{intertwine six setpar} of 
Theorem \ref{commutation intertwiner setpar}.

\medskip
Finally, we must must show relation \eqref{eq REten}
in the cases where
$ a , a+1 $ and $ a+2 $ do not all belong to the same block of $ A$, in which cases
\eqref{eq REten} corresponds to the \lq honest braid relation\rq\ 
$ \psi_{a+1} \psi_{a} \psi_{a+1}   \e(\ii, A) =
\psi_{a} \psi_{a+1} \psi_{a}   \e(\ii, A) $.

\medskip
Here, an easy case is 
$ a, a+1 $ and $a+2 $ belonging to three different
blocks of $ A $. In this case \eqref{eq REten} corresponds to the equation
\begin{equation}
\F_{a+1}q^{ -\frac{1}{2}} \F_{a}q^{ -\frac{1}{2}} \F_{a+1} q^{ -\frac{1}{2}}  \e(\ii, A) =
\F_{a}q^{ -\frac{1}{2}} \F_{a+1}q^{ -\frac{1}{2}} \F_{a} q^{ -\frac{1}{2}}  \e(\ii, A)
\end{equation}
which is immediate from Theorem \eqref{intertwine seven setpar} of
Theorem \ref{commutation intertwiner setpar}.

\medskip
Let us next consider the case where $ a \sim_A (a+1) $ but $ (a+2) \not\sim_A a $ and
$ (a+2) \not\sim_A (a+1) $. In this case, the left hand side of \eqref{eq REten} becomes
\begin{equation}\label{823}
\begin{aligned}
\psi_{a+1} \psi_{a} \psi_{a+1}   \e(\ii, A) &=
\F_{a+1} Q_{a+1}^{-1}( \sigma_{a}     \sigma_{a+1}(\ii, A) ) \psi_{a} \psi_{a+1}
\e(\ii, A) \\&=
 \F_{a+1}      \prescript{\sigma_{a} \sigma_{a+1}  \!}{}
    Q_{a}^{-1}(\ii, A ) \psi_a  \psi_{a+1} 
\e(\ii, A) \\&=
 \F_{a+1}      \psi_a  \psi_{a+1} 
 Q_{a}^{-1}(\ii, A ) \e(\ii, A)
\\&=
q^{-1} \F_{a+1}      \F_a  \F_{a+1} 
    Q_{a}^{-1}(\ii, A ) \e(\ii, A)
 \end{aligned}
\end{equation}  
where we used $ {\bf b)} $ of Lemma \ref{lemma11newPP} for
the second equality and \eqref{eq REseven} and \eqref{eq REeight}, which have already been proved,
for the third equality. But the right hand of \eqref{eq REten} is
\begin{equation}
  \psi_{a} \psi_{a+1} \psi_{a}   \e(\ii, A)   =
  q^{-1} \F_{a}      \F_{a+1}  \F_{a} 
    Q_{a}^{-1}(\ii, A ) \e(\ii, A)
\end{equation}
which is equal to the last expression of \eqref{823}, as claimed, 
in view of
\eqref{intertwine seven setpar} of Theorem 
\ref{commutation intertwiner setpar}. The other cases of 
\eqref{eq REten} are proved the same way, and so we have all in all shown that the homomorphism 
$ F: \RE \rightarrow \E  $ indeed is well-defined.

\medskip
In the other direction we now introduce a homomorphism $ G $ via
\begin{equation}
  G:  \E \rightarrow \RE, \,  \e(\ii, A)  \mapsto \e(\ii,A), \, 
  g_a \mapsto g_a \mbox{ for } 1 \le a < n, \ii \in \II, A \in \SetPar_n
\end{equation}
where $  \e(\ii, A) \in \E $ and
$g_a \in \RE $ are the elements defined in
\eqref{generalprinciples}
and 
\eqref{RRgen}, respectively. 
Note that $ \e(A) = \sum_{ \ii \in \II} \e(\ii, A) $ in $ \E $, and also in $ \RE$, 
and so $ G(\e(A) ) = \e(A)$.

\medskip
To show that $ G $ exists
we must check that the $  \e( A) $'s and the $ g_a $'s in $ \RE $ verify the relations
\eqref{eqbt-1}-\eqref{eqbt3} for $ \E$, given in Theorem \ref{ortoBT}.
Here the relations \eqref{eqbt-1} and \eqref{eqbt0} are immediate from the definitions.
To show \eqref{eqbt1}, we first observe that \eqref{define e(A)} and
\eqref{eq REfour} imply 
the relation $ \psi_a \e(A) =  \e(\sigma_aA)  \psi_a $ in $ \RE$. 
Hence, if $ a \sim_A (a+1) $ we get 
\begin{equation}\label{together1}
\begin{aligned}
  g_a \e(A) &=
  \sum_{ \ii \in \II } \left( \psi_a Q_a(\ii, A)- P_a(\ii, A) \right) \e(A) \e(\ii,A) \\
  & = \e(A) \sum_{ \ii \in \II } \left( \psi_a Q_a(\ii, A)- P_a(\ii, A) \right)\e(\ii,A)  \\
&=    \e(A) g_a = \e(\sigma_a A) g_a, 
\end{aligned}
\end{equation}  
as claimed, and if $ a \not\sim_A (a+1) $ we have $ P_a(\ii, A) \e(\ii,A) = 0 $ and so 
\begin{equation}\label{together2}
\begin{aligned}
  g_a \e(A) &=
  \sum_{ \ii \in \II }  \psi_a Q_a(\ii, A)  \e(A) \e(\ii,A) \\
  & = \e(\sigma_a A) \sum_{ \ii \in \II }  \psi_a Q_a(\ii, A) \e(\ii,A)  \\
&=     \e(\sigma_a A) g_a. 
\end{aligned}
\end{equation}  
Together, \eqref{together1} and \eqref{together1} give the proof of \eqref{eqbt1}.

\medskip
We next consider \eqref{eqbt8}. Under the hypothesis of 
\eqref{eqbt8}, that is 
$ a \not\sim_A (a+1) $, we have $ Q_a(\ii, A) \e(\ii, A) = q^{\frac{1}{2} } $ 
and $ P_a(\ii, A) \e(\ii, A) = 0 $ and so 
\begin{equation}
\begin{aligned}
g_a  \e( A)  &= \sum_{ \ii \in \II} \left( \psi_a Q_a(\ii, A)- P_a(\ii, A) \right)\e(\ii,A)
= q^{\frac{1}{2} } \sum_{ \ii \in \II}  \psi_a\e(\ii,A)
= q^{\frac{1}{2} } \sum_{ \ii \in \II} \e(\sigma_a(\ii,A))  \psi_a \implies \\
g_a^2  \e( A)  &=q \sum_{ \ii \in \II} \e(\sigma_a(\ii,A))  \psi_a^2 =
q \sum_{ \ii \in \II} \psi_a^2  \e(\ii,A)
= 
q \sum_{ \ii \in \II}   \e(\ii,A)  = q    \e(A)   
\end{aligned}
\end{equation}
where we used \eqref{eq REfour} and \eqref{eq REnine}. This shows \eqref{eqbt8}.

\medskip
For the proof of \eqref{eqbt9} we may repeat the argument for the Hecke algebra case,
in accordance with the above remark, and so \eqref{eqbt8} is proved.

\medskip
For the proof of the braid relation \eqref{eqbt2} it is enough to show that 
\begin{equation}\label{829}
g_{a+1} g_a g_{a+1} \e(\ii, A) = g_a g_{a+1} g_a \e(\ii, A) 
\end{equation}
for all $ 1 \le a < n $, $ \ii \in \II$ and $ A \in \SetPar_n$. If $a, a+1 $ and $ a+2$ belong to the same
block of $ A $ we get \eqref{829}, according to the above remark once again,
arguing as in the usual Hecke algebra situation, 
and if they belong to different blocks \eqref{829} is equivalent to 
\begin{equation}
  \psi_{a+1} q^{\frac{1}{2} } \psi_a q^{\frac{1}{2} } \psi_{a+1} q^{\frac{1}{2} }\e(\ii, A) =
  \psi_a q^{\frac{1}{2} } \psi_{a+1}q^{\frac{1}{2} } \psi_a q^{\frac{1}{2} } \e(\ii, A) 
\end{equation}
which follows directly from \eqref{eq REten}. Suppose now
that $ a \sim_A (a+1) $ but $ (a+2) \not\sim_A a $ and
$ (a+2) \not\sim_A (a+1) $. In this case the left hand side of \eqref{829} becomes
\begin{equation}
\begin{aligned}
  & \, \, \Big(  \psi_{a+1}Q_{a+1}(\sigma_{a} \sigma_{a+1}(\ii, A)) -P_{a+1}(\sigma_{a} \sigma_{a+1}(\ii, A))
  \Big) \psi_a q^{\frac{1}{2} }\psi_{a+1} q^{\frac{1}{2} } \e(\ii,A)  \\
  =&\, \,q \Big(  \psi_{a+1}
 \prescript{\sigma_{a} \sigma_{a+1}  \!}{}
    Q_{a}(\ii, A )
 - \prescript{\sigma_{a} \sigma_{a+1}  \!}{}
    P_{a}(\ii, A )
    \Big) \psi_a \psi_{a+1}  \e(\ii,A)  \\
=& \, \,q \psi_{a+1} \psi_{a} \Big(   \psi_{a+1}
    Q_{a}(\ii, A )
 -  P_{a}(\ii, A ) \Big)
 \e(\ii,A) =
q \psi_{a+1} \psi_{a}   \psi_{a+1}
    Q_{a}(\ii, A )
     \e(\ii,A)
\end{aligned}
\end{equation}
where we used Lemma \ref{lemma11newPPP} and Lemma
\ref{lemma11newPP} for the first equality and
\eqref{eq REseven} and \eqref{eq REeight} for the second equality. 
But the right hand side of \eqref{829} is equal to 
\begin{equation}
q \psi_{a+1} \psi_{a}   \psi_{a+1}
    Q_{a}(\ii, A )
     \e(\ii,A) = q \psi_{a} \psi_{a+1}   \psi_{a}
    Q_{a}(\ii, A )
     \e(\ii,A) 
\end{equation}  
which shows \eqref{829} in this case. The other cases of \eqref{829} are shown the same way.

\medskip
The only remaining relation is now the commuting braid relation
\eqref{eqbt3} which however follows directly from 
Lemma \ref{lemma41} and the definitions.

\medskip
We must finally verify that $ F $ and $ G $ are inverse isomorphisms,
that is $ G \circ F = id_{\RE} $ and $ F \circ G = id_{\E} $.
Here $ G \circ F ( \e(\ii,A) ) = \e(\ii,A)   $ and $ F \circ G ( \e(A) ) = \e(A)   $ are immediate from 
the definitions.

\medskip
By definition $ F ( y_a ) = y_a $ for $1 \le a < n $. To show
that also $ G ( y_a ) = y_a $ it is 
enough to verify that
\begin{equation}\label{toprove}
  G(X_a ) = X_a  \mbox{ for } 1 \le a < n
\end{equation}
where
$ X_a \in \RE $ is given by
\begin{equation}
  X_a = \sum_{\ii \in \II} y_a(\ii) e(\ii)  =  \sum_{\ii \in \II} q^{i_a}(1-y_a)e(\ii) =
   \sum_{ \substack{\ii \in \II \\ A \in \SetPar_n}} q^{i_a}(1-y_a)\e(\ii,A), 
\end{equation}
see \eqref{86}. We prove \eqref{toprove} by induction on $a$.
We have $ X_1 = 1 $ and so the base case $a=1 $ follows directly from the relations
\eqref{eq REcycloA} and \eqref{eq REcycloB}.
To show the induction step from $ a $ to $ a+1 $ it is enough to prove that
\begin{equation}\label{toshowthe}
q^{-1} g_a X_a g_a  \e(\ii, A) = q^{i_{a+1}}(1-y_{a+1})\e(\ii,A), 
\end{equation}
for all $ \ii \in \II $ and $ A \in \SetPar_n$, 
since $ X_{a+1} = q^{-1} g_{a} X_{a} g_a $.
Suppose first that $ a \not\sim_A (a+1) $. Then, using
\eqref{RRgen} and the inductive hypothesis together with
the relations
\eqref{eq REfour}, 
\eqref{eq REeight} and
\eqref{eq REnine}, the left hand side of 
\eqref{toshowthe} becomes
\begin{equation}
  \psi_a (q^{ i_{a+1} } (1-y_a) \psi_a \e(\ii, A) =
  q^{ i_{a+1} } (1-y_{a+1}) \psi_a^2 \e(\ii, A) =
  q^{ i_{a+1} } (1-y_{a+1}) \e(\ii, A) 
\end{equation}  
which shows \eqref{toshowthe} in this case. For $ a \sim_A (a+1) $ the proof
of \eqref{toshowthe} is carried out the same way as in Hecke algebra case,
and so \eqref{toprove} is proved. From this we deduce that 
$ G \circ F ( y_a ) = y_a  $ and $ F \circ G ( y_a ) = y_a   $ holds for all $ 1 \le a < n$,
as claimed.

\medskip
Finally, using \eqref{analogouspsi}, \eqref{RRgen} and \eqref{intertwineP} we have
\begin{equation}
\begin{aligned}
  G \circ F(\psi_a) & =  G \left( \sum_{ \ii \in \II, A \in \SetPar_n}  \Phi_a Q_a^{-1}      e(\ii,A)  \right) 
   =  G \left( \sum_{ \ii \in \II, A \in \SetPar_n}  (g_a + P_a) Q_a^{-1}      e(\ii,A)  \right) \\
& =    \sum_{ \ii \in \II, A \in \SetPar_n}  (\psi_a Q_a -P_a + P_a) Q_a^{-1}      e(\ii,A)   = \psi_a 
 \end{aligned}
\end{equation}
and 
\begin{equation}
\begin{aligned}
  F \circ G(g_a) & =
  F \left(\sum_{ \ii \in \II, A \in \SetPar_n}  (\psi_a Q_a -P_a )     e(\ii,A)   \right)  =
    \left(\sum_{ \ii \in \II, A \in \SetPar_n}  (\F_a Q_a^{-1} Q_a -P_a )     e(\ii,A)   \right) \\
    & =
 \left(\sum_{ \ii \in \II, A \in \SetPar_n}  (\F_a  -P_a )     e(\ii,A)   \right)  = g_a, 
  \end{aligned}
\end{equation}
as claimed. 
This finishes the proof of $ {\bf a)} $.

\medskip
Let us now turn to the proof of $ {\bf b)} $.
It essentially goes along the same lines
as the proof of $ {\bf a)} $, with the suitable adaptions due to
the differences between the definition of $ Q_a $ in 
\eqref{bothcases1} and \eqref{bothcases2}. 

\medskip
We first check that $  F: \REord \rightarrow \Eord $ is well-defined, in other words
that the elements
  $  \psi_a $, $\e(\ii, A^{\ord}) $ and $y_a $ of $ \Eord $ satisfy the defining relations
\eqref{eq REcycloA}--\eqref{eq REten} for $ \REord$. Here the relations
\eqref{eq REcycloA}--\eqref{eq REsix} are shown the same way as 
in part $ {\bf a)} $ and 
just like in the proof of $ {\bf a)} $, 
in order to show \eqref{eq REseven}-\eqref{eq REnine}
it is enough to consider $ a \not\sim_A (a+1) $.

\medskip
In order to show \eqref{eq REseven} we must then verify the
identity $ \psi_a y_{a+1} \e(\ii, A^{\ord}) =  y_{a}  \psi_a \e(\ii, A^{\ord}) $.
Suppose that $ a $ belongs to the block $ A_i $ of $ A^{\ord} $ and
that $ a+1 $ belongs to the block $ A_j $ of $ A^{\ord} $. If $ A_i < A_j $ then  
\eqref{eq REseven} is equivalent to 
\begin{equation} \F_a q^{-1} ( 1- q^{-i_{a+1} } X_{a+1}) \e(\ii, A^{\ord}) =
 ( 1- q^{-i_{a+1} } X_{a})  \F_a q^{ -1}\e(\ii, A^{\ord}) 
\end{equation}
which holds by \eqref{intertwine four setpar} of 
Theorem \ref{commutation intertwiner setpar}, and otherwise, if $ A_i > A_j $ then  
\eqref{eq REseven} is equivalent to 
\begin{equation} \F_a  ( 1- q^{-i_{a+1} } X_{a+1}) \e(\ii, A^{\ord}) =
 ( 1- q^{-i_{a+1} } X_{a})  \F_a \e(\ii, A^{\ord}) 
\end{equation}
which holds for the same reasons. The relation
\eqref{eq REeight} is shown the same way. 

\medskip
To show \eqref{eq REnine}, that is 
$ \psi_a^2 \e(\ii, A) = \e(\ii, A) $, we must verify 
\begin{equation} \F_a f_{ a, \sigma_a(A)}^{-1}  \F_a f_{a, A}^{-1} \e(\ii, A)
  = \e(\ii, A) .
\end{equation}
But this follows immediately from \eqref{wehaveinbothcases} and
\eqref{intertwine six setpar}.

\medskip
We are then only left with \eqref{eq REnine} which however is shown the same way as in part
$ {\bf a)} $.

\medskip
As in the proof of part $ {\bf a)} $, we now introduce a homomorphism $ G $ via
\begin{equation}
  G:  \Eord \rightarrow \REord, \,  \e(\ii, A^{\ord})  \mapsto \e(\ii,A^{\ord}), \, 
  g_a \mapsto g_a \mbox{ for } 1 \le a < n, \ii \in \II, A \in \OrdSetPar_n
\end{equation}
where $  \e(\ii, A^{\ord}) \in \Eord $ and
$ g_a \in \REord $ are the elements defined in
\eqref{generalprinciples},  
and 
\eqref{RRgen}, respectively. 
Note that $ \e(A^{\ord}) = \sum_{ \ii \in \II} \e(\ii, A^{\ord}) $ in $ \Eord $, and also in $ \REord$, 
and so $ G(\e(A^{\ord}) ) = \e(A^{\ord})$.

\medskip
To show that $ G $ exists we must check that
the $ \e(A^{\ord}) $'s and the $g_a \! $'s in $ \REord $ verify the defining relations for
$ \Eord $ given in Definition \ref{orderedBT}. 

\medskip
Here the relations \eqref{ordeqbt-1}--\eqref{ordeqbt1} are shown the same way as the relations 
\eqref{eqbt-1}--\eqref{eqbt1} are shown in part $ {\bf a)} $. 
To show \eqref{ordeqbt8} we argue as follows.
Under the hypothesis of 
\eqref{ordeqbt8}, that is 
$ a \not\sim_A (a+1) $, we have
$ P_a(\ii, A^{\ord}) \e(\ii, A^{\ord}) = 0 $. Let $ A_i $ be the
block of $ A^{\ord} $ containing $ a $ and let $ A_j $ be the block containing $ a+1 $.
If $ A_i < A_j $ we have
\begin{equation}\label{givesus}
g_a  \e( A^{\ord})  = \sum_{ \ii \in \II} \left( \psi_a Q_a(\ii, A^{\ord})- P_a(\ii, A^{\ord}) \right)\e(\ii,A^{\ord})
= q \sum_{ \ii \in \II}  \psi_a\e(\ii,A^{\ord})
= q \sum_{ \ii \in \II} \e(\sigma_a(\ii,A^{\ord}))  \psi_a 
\end{equation}
where we used \eqref{eq REfour} and \eqref{eq REnine}.
On the other hand, since $ a $ belongs to
the block $ \sigma_a A_j  $ and
$a+1 $
belongs to the block $ \sigma_a A_i  $ of $ \sigma_a A^{\ord} $ and since $ \sigma_a A_j  > \sigma_a A_i  $ 
we have $ Q_a(\sigma_a(\ii, A^{\ord})) = 1 $ and hence \eqref{givesus} gives us
\begin{equation}
  g_a^2\e( A^{\ord})  
  = q \sum_{ \ii \in \II} g_a(\e(\sigma_a(\ii,A^{\ord}))  \psi_a
  = q \sum_{ \ii \in \II} \psi_a(\e(\sigma_a(\ii,A^{\ord}))  \psi_a 
    = q \sum_{ \ii \in \II} \psi_a^2(\e(\ii,A^{\ord})   =q \e( A^{\ord})  
\end{equation}
where we used \eqref{eq REnine} for the last equality. If $ A_i > A_j $ we argue the same way, and so
\eqref{ordeqbt8} is proved.

\medskip
The proof of \eqref{ordeqbt9} is the same as the proof of \eqref{eqbt9} in part $ {\bf a)} $
and we therefore   
consider the braid relation
\eqref{ordeqbt2}, that is
\begin{equation}\label{829AA}
g_{a+1} g_a g_{a+1} \e(\ii, A^{\ord}) = g_a g_{a+1} g_a \e(\ii, A^{\ord}) .
\end{equation}
If $ a, a+1 $ and $ a+2 $ belong to the same block of $ A^{\ord} $, we repeat the argument 
of the proof of \eqref{eqbt2} 
in part $ {\bf a)} $ to arrive at \eqref{829AA}. If they belong to different belongs 
blocks of $A^{\ord} $, then \eqref{829AA} is equivalent to 
\begin{equation}\label{829AAA}
  \psi_{a+1}
f_{a+1, \sigma_{a} \sigma_{a+1} A^{\ord}}
\psi_a   f_{a, \sigma_{a+1}A^{\ord}}  \psi_{a+1} f_{a, A^{\ord}} \e(\ii, A^{\ord}) =
\psi_a  f_{a, \sigma_{a+1}\sigma_a A^{\ord}} \psi_{a+1}     f_{a+1, \sigma_a A^{\ord}}    \psi_a   f_{a, A^{\ord}} \e(\ii, A^{\ord}) 
\end{equation}
which is a consequence of the following identity 
\begin{equation}\label{829AAA}
f_{a+1, \sigma_{a} \sigma_{a+1} A^{\ord}}
  f_{a, \sigma_{a+1}A^{\ord}}  f_{a, A^{\ord}} = 
  f_{a, \sigma_{a+1}\sigma_a A^{\ord}}      f_{a+1, \sigma_a A^{\ord}}      f_{a, A^{\ord}}. 
\end{equation}
Finally, if precisely two of $ a, a+1 $ or $ a+2 $ belong to the same block of
$ A^{\ord} $ we argue the same way as in part $ {\bf a)} $ to arrive at \eqref{829AA}.
This completes the proof of \eqref{ordeqbt2} and the proof of 
the commuting braid relation \eqref{ordeqbt3} is trivial. We have thus proved that
$ G $ is well-defined.

\medskip
Finally, the proof that $ F $ and $ G $ are inverse isomorphisms is carried out
as in $ {\bf a)} $. This completes the proof of the Theorem.
\end{dem}

 \end{document}